\documentclass[11pt,oneside]{amsart}

\usepackage[a4paper, total={6in, 8in}]{geometry}

\usepackage{setspace}

\usepackage{amssymb}   
\usepackage[hidelinks]{hyperref}
\usepackage{amsmath}

\usepackage{mathrsfs}

\usepackage{stmaryrd}

\usepackage{amsthm}
\usepackage{newlfont}
\usepackage{amscd}
\usepackage{mathtools}
\usepackage{bm}
\usepackage{tikz-cd}
\usepackage{graphicx}

\usepackage{hyperref}
\usepackage{framed}
\usepackage{comment}

\theoremstyle{definition}
\newtheorem{definition}{Definition}[section]

\theoremstyle{remark}
\newtheorem{remark}[definition]{Remark}

\theoremstyle{theorem}
\newtheorem{theorem}[definition]{Theorem}

\theoremstyle{corollary}
\newtheorem{corollary}[definition]{Corollary}

\theoremstyle{lemma}
\newtheorem{lemma}[definition]{Lemma}

\theoremstyle{example}
\newtheorem{example}[definition]{Example}

\theoremstyle{prop}
\newtheorem{prop}[definition]{Proposition}

\usepackage{enumitem}

\newcommand {\Hom} {\operatorname{Hom}}
\newcommand {\End} {\operatorname{End}}
\newcommand {\Ext} {\operatorname{Ext}}

\newcommand {\Pro} {\operatorname{Pro}}

\newcommand {\coker} {\operatorname{coker}}
\newcommand {\QCoh} {\mathrm{QCoh}}
\newcommand {\Coh} {\mathrm{Coh}}
\newcommand {\rank} {\operatorname{rank}}
\newcommand {\Stab} {\operatorname{Stab}}

\title[Continuum envelopes on FF curves and elliptic curves]{Continuum envelopes on Fargues--Fontaine curves and elliptic curves}

\begin{document}
\author{Heng Du} 
\address{Yau Mathematical Sciences Center, Tsinghua University, Beijing 100084, China}
\email{hengdu@mail.tsinghua.edu.cn}
\author{Qingyuan Jiang} 
\address{Department of Mathematics,
The Hong Kong University of Science and Technology, Clearwater Bay, Kowloon, Hong Kong.} 
\email{jiangqy@ust.hk}
\author{Yucheng Liu} 
\address{College of Mathematics and Statistics, Center of Mathematics, Chongqing University, Chongqing 401331, China}
\email{noahliu@cqu.edu.cn}
\begin{abstract} 
The theory of stability conditions (see \cite{bridgeland2007stability}), which originated in string theory, provides a natural language in which to reformulate the main theorem of \cite{LeBrasresult}. In this paper, we study its applications to Fargues--Fontaine curves, which play a fundamental role in p-adic Hodge theory and the geometrization of local Langlands program (see \cite{farguesfontaine-courbes} and \cite{fargues2021geometrization}). This leads us to construct the quasi-coherent sheaves $\mathcal{O}(\theta^{\pm})$ from the convergents of an irrational number $\theta$. 

We define the continuum envelope $\QCoh_{\mathbb{R}}(X_{FF})$ to be the smallest abelian subcategory of $\QCoh(X_{FF})$ containing $\Coh(X_{FF})$ and the sheaves $\mathcal{O}(\theta^{\pm})$. We study its homological algebra using Farey diagrams and show that the homological properties of $\mathcal{O}(\theta^{\pm})$ depend strongly on the Diophantine properties of $\theta$. From this perspective, the Fargues--Fontaine curve exhibits strong similarities with elliptic curves.
\end{abstract}

\maketitle

\tableofcontents

\section{Introduction}

	The theory of stability conditions was introduced by Bridgeland in \cite{bridgeland2007stability}, motivated by Douglas's work on D-branes and $\Pi$-stability (see \cite{douglas2002dirichlet}). This theory was further studied by Kontsevich and Soibelman (see \cite{kontsevich2008stability}), and it turns out to be related to many different branches of mathematics including algebraic geometry (\cite{bayer2014projectivity}, \cite{bayer2014mmp}), symplectic geometry (\cite{Quadraticdifferentialsasstabilityconditions}, \cite{Flatsurfacesandstabilityconditions}), representation theory (\cite{kontsevich2008stability}, \cite{Scatteringdiagrams}), and curve counting theories (\cite{CurvecountingtheoriesviastableojectsI}, \cite{CurvecountingtheoriesviastableobjectsII}). 

The Fargues--Fontaine curve, introduced by Fargues and Fontaine in \cite{farguesfontaine-courbes}, has brought a geometric perspectives into the study of $p$-adic Hodge theory, which traditionally focused largely on semilinear algebraic categories: many categories appearing in $p$-adic Hodge theory, such as filtered $\varphi$-modules, $\varphi$-modules or $(\varphi,\Gamma)$-modules, $B$-pairs, Finite Dimension Banach Spaces, can now be understood using the category of coherent sheaves on this curve (\cite{farguesfontaine-courbes}, \cite{LeBrasresult}). A common feature of these semilinear algebraic categories is their slope formalism, Fargues and Fontaine related these slope formalisms to the Harder--Narasimhan formalism of vector bundles on the Fargues--Fontaine curves.

%Additionally, the curve has a notable impact on the geometrization of local class field theory and the local Langlands correspondence, as further explored in references \cite{Fargues-simplyconnected} and \cite{fargues2021geometrization}.

The theory of Bridgeland stability conditions gives a systematic way to extend the concept of slopes defined on an abelian category to the bounded derived category generated by this abelian category. The aim of this paper is to apply the theory of Bridgeland stability conditions to the study of the bounded derived category of coherent sheaves on Fargues--Fontaine curves. 

Let us first recall the definition of a stability condition on a triangulated category $\mathcal{D}$. A stability condition $\sigma$ on $\mathcal{D}$ is a pair $(\mathcal{A},Z)$, where $\mathcal{A}$ is the heart of a bounded $t$-structure on $\mathcal{D}$, and $Z:K_0(\mathcal{A})\rightarrow \mathbb{C}$ is a group homomorphism, satisfying certain properties (see Definition \ref{first WSC} for details). The function $Z$ defines a slope function for all objects in $\mathcal{A}$, and the slope function defines the notion of semistable and stable objects inside $\mathcal{A}$ in the usual way. An important result is that the slope formalism on $\mathcal{A}$ defined using $Z$ can be extended to the triangulated category $\mathcal{D}$ by introducing the concept of \emph{slicing}.

\begin{definition}\label{slicing}
			A \textit{slicing} on a triangulated category $\mathcal{D}$ consists of full subcategories $\mathcal{P}(\phi)\in\mathcal{D}$ for each $\phi\in\mathbb{R}$, satisfying the following axioms:

			\par
			(a) for all $\phi \in \mathbb{R}$, $\mathcal{P}(\phi+1)=\mathcal{P}(\phi)[1]$;
			
			\par
			(b) if $\phi_1>\phi_2$, then $\Hom_{\mathcal{D}}(A_1,A_2)=0$ for any $A_1\in \mathcal{P}(\phi_1)$ and $A_2\in \mathcal{P}(\phi_2)$;
			
			\par

			(c) for every $0\neq E\in\mathcal{D}$ there is a sequence of real numbers
			
			$$\phi_1>\phi_2>\cdots>\phi_m$$and a sequence of morphisms 
			
			$$0=E_0\xrightarrow{f_1}E_1\xrightarrow{f_2} \cdots \xrightarrow{f_m}E_m=E $$such that the cone of $f_j$ is in $\mathcal{P}(\phi_j)$ for all $j$. And the object $Cone(f_i)$ is called the $i$-th Harder--Narasimhan factor of $E$. 
			
		\end{definition}

\begin{remark}\label{remark:different hearts}An object $E\in\mathcal{P}(\phi)$ is called a semistable object of phase $\phi$. Condition (c) can be viewed as a generalization of the Harder--Narasimhan filtration. A slicing can be viewed as a refinement of a t-structure. In fact, for any real number $\phi\in\mathbb{R}$, it defines two t-structures $(\mathcal{P}(>\phi-1),\mathcal{P}(\leq \phi))$ and $(\mathcal{P}(\geq \phi-1),\mathcal{P}(< \phi))$ on $\mathcal{D}$, where $\mathcal{P}(>\phi-1)$ is the full subcategory consisting of the objects whose Harder--Narasimhan factors (as in Definition \ref{slicing}.(c)) have phases greater than $\phi-1$, similar definitions apply to the full subcategories $\mathcal{P}(\leq \phi)), \mathcal{P}(\geq \phi-1),\mathcal{P}(< \phi)$. We use $\mathcal{P}(\phi-1,\phi]$ and $\mathcal{P}[\phi-1,\phi)$ to denote the associated hearts respectively.
\end{remark}
An important fact is that giving a stability condition $\sigma=(\mathcal{A},Z)$ on $\mathcal{D}$ with $\mathcal{A}$ being the heart of a bounded $t$-structure on $\mathcal{D}$, is equivalent to giving a pair $(\mathcal{P},Z)$, where $\mathcal{P}$ is a slicing, and $Z\colon K_0(\mathcal{A})\rightarrow \mathbb{C}$ is as before (see Definition \ref{another definition of BSC} for details). Indeed, the equivalence is given by setting $\mathcal{P}(\phi)$ to be the full subcategory consisting of $\sigma$-semistable objects in $\mathcal{A}$ of slope $-\cot(\pi\phi)$  for $\phi\in(0,1]$, so here $\mathcal{A}=\mathcal{P}(0,1]$. %\textcolor{red}{Heng: Maybe recall how we get the equivalence here very briefly if the idea is clear, or cite a later part.}

We recall the construction of schematic Fargues--Fontaine curves. Recall that given  a complete algebraically closed nonarchimedean field $C^\flat$ over $\bar{\mathbb{F}}_p$, a complete discrete valuation field $E$ with perfect residue field, and a uniformizer $\pi$ of $E$, there is a curve $X_{C^\flat,E}$ defined as the projective scheme associated to the graded algebra $P_{C^\flat,E}$ with 
\[
P_{C^\flat,E}\coloneqq \oplus_{i\geq 0} ({B}_{\mathrm{cris}}^+)^{\varphi=\pi^i}.
\]
In particular, if $E=\mathbb{Q}_p$, we define $X_{FF}=X_{C^\flat,\mathbb{Q}_p}$ to be the absolute Fargues--Fontaine curve.

We will apply the theory of stability conditions to study the bounded derived category of (quasi-)coherent sheaves on the absolute Fargues--Fontaine curve $X_{FF}$. Note that although we will state our results only for the absolute Fargues--Fontaine curve, all the results in this paper hold with only slight changes for the general Fargues--Fontaine curves $X_{C^{\flat},E}$. 

 Bridgeland proved that the space of stability conditions has a complex manifold structure, and it carries a right action of the group $\widetilde{GL}^+(2,\mathbb{R})$, the universal cover of $GL^+(2,\mathbb{R})$. Moreover, the group action of $\widetilde{GL}^+(2,\mathbb{R})$ does not change the class of semistable objects. Using this, we can prove the following theorem, which means that the stability condition considered in \cite{farguesfontaine-courbes} is essentially the only stability condition in $\Stab(X_{FF})$.
\begin{theorem}\label{Theorem: stability conditions on ff curve in introduction section}
The action of $\widetilde{GL}^+(2,\mathbb{R})$ on $\Stab(X_{FF})$ is free and transitive. 
\end{theorem}

\begin{remark}
    If $X$ is a complex smooth projective curve, the Mumford stability condition $(\Coh(X),-\mathrm{deg}+\sqrt{-1}\cdot \rank)$ is a Bridgeland stability condition on $D^b(X)$. Similarly, let $X_{FF}$ denote the absolute Fargues--Fontaine curve, and let $\sigma=(\Coh(X_{FF}), -\mathrm{deg}+\sqrt{-1}\cdot \rank)$ be the stability condition considered in \cite{farguesfontaine-courbes}. Please see \S \ref{prelimenary section} for more details of the construction of the Fargues--Fontaine curve and its associated stability conditions.
\end{remark}

This differs markedly t from the case of the complex projective line $\mathbb{P}^1$, on which there are many different stability conditions (see \cite[Theorem 1.2]{StabilitymanifoldofP1}). In fact, the proof of Theorem \ref{Theorem: stability conditions on ff curve in introduction section} shows that the Fargues--Fontaine curve is rather similar to elliptic curves. This similarity is further reinforced by the following version of Serre duality in $X_{FF}$. 

\begin{theorem}[{\cite[Corollary 3.10 \& Remark 3.11]{anschutz2021fourier}}]
Let $\mathcal{E}, \mathcal{F} \in D^b(\mathrm{Coh}(X_{FF}))$. There is a canonical equivalence in $D_{v}(C; \mathbb{Q}_p)$:
	$$ R \tau_* (R \mathcal{H}\mathrm{om}_{X_{FF}}(\mathcal{E},\mathcal{F})) \simeq R \mathcal{H}\mathrm{om}_v (R \mathcal{H}\mathrm{om}_{X_{FF}}(\mathcal{F}, \mathcal{E}[1]), \underline{\mathbb{Q}_p}[1]).$$
\end{theorem}

Here $\tau$ is the morphism from the topos associated with the big pro\'etale site over the Fargues--Fontaine curve $X_{FF}$ to the topos associated with $\mathrm{Perf}_{C^{\flat},pro\text{\'et}}$ as in \cite[\S 6.2]{LeBrasresult}. Please see \S \ref{subsubsection:Serre duality} for more details of the notation in this theorem. This resembles the Grothendieck--Verdier duality for families of elliptic curves.

In particular, there is a canonical equivalence of Banach--Colmez spaces:
	$$\mathbb{H}\mathrm{om}(\mathcal{E},\mathcal{F}) \simeq \mathbb{E}\mathrm{xt}^{1}(\mathcal{F},\mathcal{E})^{\vee},$$
where we suggestively let $(-)^{\vee}$ denote the functor $\mathcal{H}\mathrm{om}_v(-, \underline{\mathbb{Q}_p}[1])$. This is exactly the form of Serre duality for elliptic curves.

Another example of the connection between the theory of stability conditions and Fargues--Fontaine curves is the following restatement of Le Bras' result. 

 \begin{theorem}[{\cite[Theorem 1.2]{LeBrasresult}}]\label{restatement of Lebras' result}
    Let $\sigma$ be the standard stability condition on $D^b(X_{FF})$, its associated slicing $\mathcal{P}$ gives us the following two abelian categories.\begin{enumerate}
        \item The category $\mathcal{P}(0,1]$ is the category of coherent sheaves on $X_{FF}$.

        \item The category $\mathcal{P}[\frac{1}{2},\frac{3}{2})$ is equivalent to the category of Banach--Colmez spaces.
    \end{enumerate}
    \end{theorem}

\begin{figure}[ht]
		\centering
		\includegraphics[scale=0.8]{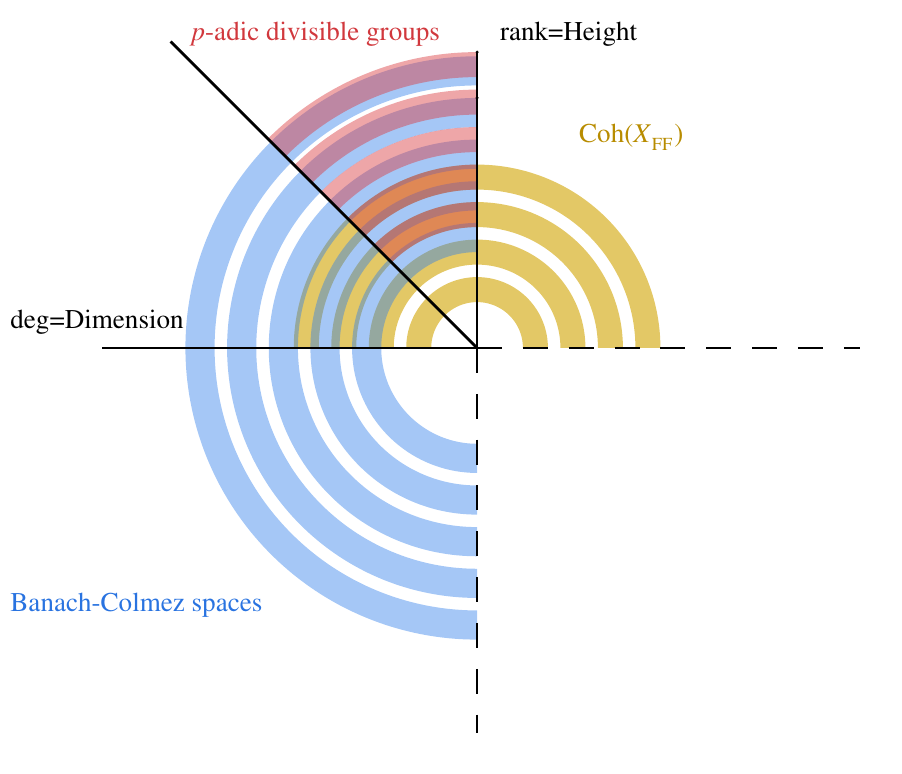}
		\caption{}
        \label{Figure 1}
	\end{figure}

\begin{remark}
    In fact, one can also restate \cite[Theorem 5.1.4.(ii)]{Moduliofp-divisiblegroups} in the following way. The category of p-divisible groups over $\mathcal{O}_C/p$ up to isogeny is equivalent to $\mathcal{P}[\frac{1}{2},\frac{3}{4}]$. Note that this category is only quasi-abelian as the length of the interval $[\frac{1}{2},\frac{3}{4}]$ is less than $1$. See Figure \ref{Figure 1} for a pictorial description of these three categories.
\end{remark}

By a general result in the theory of stability conditions, we can conclude that the category $\Coh(X_{FF})$ is Noetherian, and the category of Banach--Colmez spaces is Artinian.  A natural question is the following: what about $\mathcal{P}(\phi,\phi+1], \mathcal{P}[\phi,\phi+1)$ for other real numbers $\phi$? 

We have the following result, which answers the question. 

\begin{theorem}\label{Hearts of FF curves in introduction}
	Let $\sigma=(\mathcal{A},Z)$ be the standard stability condition on $D^b(X_{FF})$, and $\mathcal{P}$ be its associated slicing.  We have the following results. 
	\begin{enumerate}
		\item If $\cot(\pi \phi)$ is in $\mathbb{Q}$ or ${\infty}$ , the heart $\mathcal{P}(\phi, \phi+1]$ is Noetherian and non-Artinian,  the heart $\mathcal{P}[\phi, \phi+1)$ is Artinian and non-Noetherian.
		\item If $\cot(\pi \phi)$ is an irrational number, the heart $\mathcal{P}(\phi, \phi+1]=\mathcal{P}[\phi, \phi+1)$ is neither Noetherian nor Artinian.
	\end{enumerate}

\end{theorem}
Before discussing this result in more details, we note that there is a similar result for elliptic curves (see \cite{Differentheartsonellipticcurves}). In that case, the heart $\mathcal{P}(\phi,\phi+1]$ is closely related to holomorphic vector bundles on the noncommutative tori $A_{\theta}$. This will be discussed more explicitly in \S \ref{noncommutative tori}.

The above result follows from Theorem \ref{Theorem between minimal triangles and Farey triangles}, which establishes a relation between the homological algebra of $D^b(X_{FF})$ and the Farey diagrams. As a consequence, we construct two infinite sequences $$  \mathcal{O}(\beta_{0})\hookrightarrow \mathcal{O}(\beta_{2})\hookrightarrow\cdots \hookrightarrow \mathcal{O}(\beta_{2i})\hookrightarrow \cdots$$ and $$\cdots\twoheadrightarrow\mathcal{O}(\beta_{2i+1})\twoheadrightarrow \cdots\twoheadrightarrow \mathcal{O}(\beta_{1})\twoheadrightarrow \mathcal{O}(\beta_{-1})$$ in $\mathcal{P}(0,1]=\Coh(X_{FF})$, where $\beta_i=\frac{p_i}{q_i}$ is the $i$-th convergent of the continued fraction of the irrational number $\theta=-\cot(\pi\phi)$, and $\mathcal{O}(\beta_i)$ is the unique stable vector bundle on $X_{FF}$ (see \S \ref{subsection:Fargues--Fontaine curve} for the construction of such a vector bundle). These two infinite sequences imply that the heart $\mathcal{P}(\phi, \phi+1]=\mathcal{P}[\phi, \phi+1)$ is neither Noetherian nor Artinian, since monomorphisms and epimorphisms exchange roles when one pass to the abelian category $\mathcal{P}(\phi, \phi+1]$.

These two infinite sequences also lead us to consider their colimit and limit objects in $\QCoh(X_{FF})$, which we denote by $\mathcal{O}(\theta^-,\{g_{2i}\})$ and $\mathcal{O}(\theta^+,\{f_{2i+1}\})$, where $\{g_{2i}\}$, $\{f_{2i+1}\}$ are sequences of nonzero morphisms $g_{2i}:\mathcal{O}(\beta_{2i})\rightarrow\mathcal{O}(\beta_{2i+2})$, $f_{2i+1}:\mathcal{O}(\beta_{2i+1})\rightarrow \mathcal{O}(\beta_{2i-1})$. As the notation indicates, these objects depend on the choice of these morphisms in a nontrivial way. There are similar limit and colimit objects on an elliptic curve $E$ defined over an algebraically closed field $k$, denoted by $$\mathcal{O}(\theta^{+}, \{\mathcal{L}_{2i+1}\}, \{f_{2i+1}\}),\ \mathcal{O}(\theta^{-},\{\mathcal{L}_{2i}\}, \{g_{2i}\})$$ where $\mathcal{L}_{2i+1}$ ($\mathcal{L}_{2i}$) is a vector bundle on $E$ of rank $q_{2i+1}$ ($q_{2i}$) and degree $p_{2i+1}$ ($p_{2i}$), and $f_{2i+1}:\mathcal{L}_{2i+1}\rightarrow \mathcal{L}_{2i-1}$, $g_{2i}:\mathcal{L}_{2i}\rightarrow \mathcal{L}_{2i+2}$ are nontrivial morphisms. Note that unlike in the Fargues--Fontaine curve case, by Atiyah's classification of vector bundles on elliptic curves (see \cite{Atiyahclassification}), the isomorphism class of $\mathcal{L}_{i}$ lies in a coarse moduli space isomorphic to $\mathrm{Pic}^0(E)\simeq E$.

In \S \ref{section: homological algebra}, we describe the respective coarse moduli space of colimit objects  $\mathcal{O}(\theta^-,\{g_{2i}\})$ and $\mathcal{O}(\theta^{-},\{\mathcal{L}_{2i}\}, \{g_{2i}\})$ on Fargues--Fontaine curves and elliptic curves; they are infinite-dimensional. Essentially, the moduli space of  $\mathcal{O}(\theta^-,\{g_{2i}\})$ comes from the choices of morphisms $$\{g_{2i}:\mathcal{O}(\beta_{2i})\rightarrow\mathcal{O}(\beta_{2i+2})\}$$ while the moduli space of $\mathcal{O}(\theta^{-},\{\mathcal{L}_{2i}\}, \{g_{2i}\})$ comes from the moduli of the bundles $\mathcal{L}_{2i}$. We define the continuum envelope $\QCoh_{\mathbb{R}}(X_{FF})$ ($\QCoh_{\mathbb{R}}(E)$) to be the smallest abelian subcategory of $\QCoh(X_{FF})$ ($\QCoh(E)) $ containing $\Coh(X_{FF})$ ($\Coh(E))$ and these colimit and limit objects. Heuristically, the continuum envelope $\QCoh_{\mathbb{R}}(X_{FF})$ ($\QCoh_{\mathbb{R}}(E)$) is a completion of $\Coh(X_{FF})$ ($\Coh(E)$), analogous to the Archimedean completion of $\mathbb{Q}\hookrightarrow\mathbb{R}$.

As a side remark, one can show that the $K$-class of $\mathcal{O}(\theta^-,\{g_{2i}\})$ is independent of the choice of $\{g_{2i}\}$ (see Lemma \ref{K class lemma}). By this result, we can show there is no Riemann--Roch theorem for the continuum envelope $\QCoh_{\mathbb{R}}(X_{FF})$ (see Example \ref{golden ratio}).

We further study the homological algebra of $\QCoh_{\mathbb{R}}(X_{FF})$ in \S \ref{section: homological algebra}. We have the following theorem.

\begin{theorem}\label{main theorem in introduction}
     Let $\theta=[a_0;a_1,a_2,\cdots]$ be an irrational number, and $\beta_i\coloneqq [a_0;a_1,\cdots,a_i]$ be its $i$-th convergent. We have the following results.
      \begin{enumerate}
        \item For any sequences of nonzero morphisms $\{g_{2i}\}$ and $\{f_{2j+1}\}$, we have the following short exact sequence $$0\rightarrow \Pi_{i=0}^{\infty} \Hom(\mathcal{O}(\beta_i),\mathcal{O}({\beta_i}))^{a_{i+1}}\rightarrow \Hom(\mathcal{O}(\theta^-,\{g_{2i}\}), \mathcal{O}(\theta^+,\{f_{2j+1}\}))\xrightarrow{\vartheta} C\rightarrow 0,$$ where $\vartheta$ can be viewed as Fontaine's $\vartheta$ map. 
        \item The colimit object $\mathcal{O}(\theta^-,\{g_{2i}\})$ is an indecomposable vector bundle of infinite rank for any sequence $\{g_{2i}\}$ of nonzero morphisms. If the set  $\{a_i\}$ is bounded above, then the $\mathbb{Q}_p$-division algebra $$\Hom(\mathcal{O}(\theta^-,\{g_{2i}\}),\mathcal{O}(\theta^-,\{g_{2i}\}))$$ is finite-dimensional for any choice of $\{g_{2i}\}$. Moreover, its $\mathbb{Q}_p$-dimension divides $c(\theta)^2$, where $c(\theta)$ is a finite positive integer depending on $\theta$ (see Definition \ref{defn: c(theta)} for the definition of $c(\theta)$).

        \item For almost every irrational number $\theta$ (with respect to Lebesgue measure), we have $$\Hom(\mathcal{O}(\theta^-,\{g_{2i}\}),\mathcal{O}(\theta^-,\{g_{2i}\}))\simeq \mathbb{Q}_p$$ for any choice of nontrivial morphisms $\{g_{2i}\}$. 
        \item There exists an infinite set $S$ of irrational numbers such that for any irrational number $\theta\in S$, there exists a sequence of nontrivial morphisms $\{g_{2i}\}$ such that $$ \Hom(\mathcal{O}(\theta^-,\{g_{2i}\}),\mathcal{O}(\theta^-,\{g_{2i}\}))$$ is an infinite-dimensional $\mathbb{Q}_p$-division algebra.
    \end{enumerate}
\end{theorem}
By Theorem \ref{main theorem in introduction}.(1), the Hom-space $$\Hom(\mathcal{O}(\theta^-,\{g_{2i}\}), \mathcal{O}(\theta^+,\{f_{2i+1}\}))$$ can be viewed as a $p$-adic Banach space with $C$ dimension $1$ and  uncountable $\mathbb{Q}_p$-dimension.

%     By Langrange's theorem and Theorem \ref{main theorem in introduction}.(2), we know that if $\theta$ is a real quadratic irrational number, then the division algebra $$\Hom(\mathcal{O}(\theta^-,\{g_{2i}\}),\mathcal{O}(\theta^-,\{g_{2i}\})$$ is finite-dimensional. 
    
%  We also have the following example.

% \begin{example}\label{Example:Zaremba conjecture in introduction}
%     If we let
%      $$\mathcal{C}_2\coloneqq \left\{[0;a_1,\cdots,a_i,\cdots] \mid \text{$1\leq a_i\leq 2$ for any $i\geq 1$}\right\}\subset [0,1],$$ 
%     it is well known that $\mathcal{C}_2$ is a Cantor set, with Hausdorff dimension around $0.53128\cdots$ (see \cite{computingthehausdorffdimensionofE_2}).

%     For any irrational number $\theta\in\mathcal{C}_2$, the $\mathbb{Q}_p$-dimension of $$\Hom(\mathcal{O}(\theta^-,\{g_{2i}\}),\mathcal{O}(\theta^-,\{g_{2i}\}))$$ could only be $1$, $2$ or $4$.
% \end{example}

Recall the definition of badly approximable real numbers (see \cite[Definition 1.31]{aigner2015markov}).

\begin{definition}
    A  real number $\theta$ is said to be \textit{badly approximable} if there exists a constant $C>0$ such that $$|\theta-\frac{p}{q}|>\frac{C}{q^2}$$ for all $\frac{p}{q}\in\mathbb{Q}$ with $\frac{p}{q}\neq \theta$.
\end{definition}

Note that a rational number is also badly approximable by definition. Hence, we have the following corollary.

\begin{corollary}
    For any real number $\theta\in\mathbb{R}$, if $\theta$ is badly approximable, then the $\mathbb{Q}_p$-division algebra $$\Hom(\mathcal{O}(\theta^-,\{g_{2i}\}),\mathcal{O}(\theta^-,\{g_{2i}\}))$$ is finite-dimensional for any choice of $\{g_{2i}\}$.
\end{corollary}
\begin{proof}
 By \cite[Proposition 1.32]{aigner2015markov}, we know that a real number $\theta$ is badly approximable if and only if the partial quotients $\{a_i\}$ of its continued fraction expansion are bounded. Hence it follows directly from Theorem \ref{main theorem in introduction}.(2).
\end{proof}

These results show that the homological algebra of $\mathcal{O}(\theta^{-}, \{g_{2i}\})$ relies heavily on the Diophantine properties of $\theta$, and it relates the results from classical number theory to the homological algebra of $\mathcal{O}(\theta^{-},\{\mathcal{L}_{2i}\}, \{g_{2i}\})$.

In the case when $$ \Hom(\mathcal{O}(\theta^-,\{g_{2i}\}),\mathcal{O}(\theta^-,\{g_{2i}\}))$$ is an infinite-dimensional $\mathbb{Q}_p$-division algebra, one can show that its dimension is still countable (see Theorem \ref{thm: infinite-dimensional division algebra}). On the other hand, we can get more complicated division algebras by considering the quotient categories $\mathcal{P}(\phi,\phi+1]/\mathcal{P}(\phi+1)$. These division algebras can be viewed as the $\mathrm{SL}(2,\mathbb{Z})$ variants of Colmez-Fontaine division algebra $\mathscr{C}$ (see \cite[\S 5, \S 9]{BanachColmezspaces} for the definition). 

Indeed, we show that the quotient category $\mathcal{P}(\phi,\phi+1]/\mathcal{P}(\phi+1)$, denoted by $\mathcal{Q}_{\lambda}$, is semi-simple with a single simple object, denoted by $S_{\lambda}$, where we assume that $\lambda=-\mathrm{cot}(\pi\phi)\coloneqq \frac{p}{q}\in\mathbb{Q}_{\infty}$. We have the following result for the division algebras $\End_{\mathcal{Q}_{\lambda}}(S_{\lambda})$.

 \begin{theorem}\label{thm: generalization of clomez-Fontaine corps}
     
  For any $\frac{r}{s}\in \mathbb{Q}$ with $|det\begin{pmatrix}
        r & s \\ p & q
    \end{pmatrix}|=n>0,$ we have a monomorphism of algebras $$D_{\frac{r}{s}}\hookrightarrow \mathrm{M}_n(\mathrm{End}_{\mathcal{Q}_{\lambda}}(S_{\lambda})) ,$$ where $D_{\frac{r}{s}}$ is the central simple division algebra corresponding to $\frac{r}{s}$ in the Brauer group $\mathbb{Q}/\mathbb{Z}=Br(\mathbb{Q}_p)$, and $\mathrm{M}_n(\mathrm{End}_{\mathcal{Q}_{\lambda}}(S_{\lambda}))$ denotes the algebra of $n\times n$ matrices over $\mathrm{End}_{\mathcal{Q}_{\lambda}}(S_{\lambda})$.

For $(r,s)=(1,0)$, if $det\begin{pmatrix}
        1 &0 \\ p & q
    \end{pmatrix}=q>0,$ we have monomorphisms of algebras $$\mathcal{C}\hookrightarrow \mathrm{M}_q(\mathrm{End}_{\mathcal{Q}_{\lambda}}(S_{\lambda})),$$ for any untilt $\mathcal{C}$ of $C^{\flat}$ in characteristic 0.
\end{theorem}
In the special case where $\phi=\frac{1}{2}$, the endomorphism algebra $\End_{\mathcal{Q_0}}(S_{0})$ is Colmez--Fontaine division algebra.

In the last section of this paper, we discuss the space of stability conditions on the twistor projective line and its products.  We also discuss the relation between the hearts $\mathcal{P}(\phi,\phi+1]$ on a complex elliptic curve and holomorphic vector bundles on noncommutative tori $A_{\theta}$. This motivates us to ask the following  question:  just as noncommutative tori can be viewed as noncommutative deformations of elliptic curves, does the Fargues–Fontaine curve similarly admit noncommutative deformations?

Although the method of this paper applies to many interesting stability conditions on smooth projective schemes. We restrict our attention to the Fargues--Fontaine curve and elliptic curves, because of their simplicity, strong similarities, and arithmetic interest.   

\subsection*{Outline of the paper} In \S \ref{prelimenary section}, we briefly recall some preliminary results on Bridgeland stability conditions, 
 Fargues--Fontaine curves, and Banach--Colmez spaces. We also restate some results from \cite{LeBrasresult}, \cite{anschutz2021fourier} from our point of view, and we prove Theorem \ref{Theorem: stability conditions on ff curve in introduction section} in this section. In \S \ref{Section: Farey diagrams}, we recall some classical results about Farey triangles, Farey tessellation, and their relation with continued fractions. In \S \ref{section:limit and colimit objects}, we relate Farey triangle to minimal triangles in $D^b(X_{FF})$, and prove Theorem \ref{Hearts of FF curves in introduction}. We define the continuum envelope in this section. In \S \ref{section: homological algebra}, we study the homological algebra of the continuum envelope $\QCoh_{\mathbb{R}}(X_{FF})$ and relate results in classical number theory to the $\mathbb{Q}_p$-algebras and vector spaces arising from it. In \S \ref{section:final remarks}, we discuss two topics related to our results: the twistor projective line and noncommutative tori. 
	\subsection*{Notation and Conventions} In this paper, we use $\mathbb{H}^n$ to denote the hyperbolic space in dimension n. In the special case when $n=2$, we also identify $\mathbb{H}^2$ with the upper half-plane in $\mathbb{C}$. We denote $\mathbb{R}_{\infty}\coloneqq \mathbb{R}\cup \{\infty\}$, $\mathbb{Q}_{\infty}\coloneqq \mathbb{Q}\cup \{\infty\}$. We adopt the convention $\frac{1}{0}=\infty$. We use $\gcd(p,q)$ to denote the greatest common divisor of two integers $p,q$. For a number $r$ in $\mathbb{Q}_{\infty}$, we always write it in a standard form $r=\frac{p}{q}$, where $\gcd(p,q)=1$ and $q\geq0$.  For such an expression, the complex number $-p+\sqrt{-1}\cdot q$ is a primitive integral vector in $\mathbb{H}^2\cup \mathbb{R}_{< 0} $. We use the numbers in $\mathbb{Q}_{\infty}$ and the primitive integral vectors in $\mathbb{H}^2\cup \mathbb{R}_{< 0} $ interchangeably. An elliptic curve is implicitly assumed to be defined over an algebraically closed field of arbitrary characteristic unless otherwise specified.

 For set-theoretic considerations, we adopt the same strategy as Lurie (\cite[\S 1.2.15]{HTT}) and distinguish ``small” and ``large” objects. For meticulous readers, this means that we will assume the existence of certain uncountable strongly inaccessible cardinal $\kappa$. Then a set is said to be \emph{small} if its cardinality is strictly smaller than $\kappa$, and a simplicial set is \emph{small} if its collection of nondegenerate simplices is small. 
 
 %\textcolor{green}{Maybe we do not need the following examples in this subsection} %In particular, if $Y$ is a Noetherian scheme, both the abelian category $\mathrm{Coh}(Y)$ of coherent sheaves and the stable $\infty$- category $\mathrm{Perf}(Y)$ of perfect complexes are \emph{small}. On the other hand, the abelian category $\mathrm{QCoh}(Y) (\simeq \Ind \mathrm{Coh} (Y))$ of quasi-coherent sheaves and the stable $\infty$-category $\mathrm{D}_{\mathrm{qcoh}}(Y) (\simeq \Ind \mathrm{Perf} (Y))$ of unbounded quasi-coherent complexes are \emph{large}.
 
\medskip
\noindent
\textbf{Acknowledgments.} We would like to thank Sebastian Bartling, Hanfeng Li, Shilei Kong, Koji Shimizu and Zhiyu Tian for helpful discussions. We are grateful to Runlin Zhang for his proof of Lemma \ref{lemma:rare case with bigger dimension}. H.D. is supported by National Key R{\&}D Program of China No. 2023YFA1009703. Y.L. is financially supported by NSFC grant 12201011.

\section{Stability conditions and Fargues--Fontaine curves}\label{prelimenary section}

This section is mainly preliminary. We briefly recall some basic definitions and results concerning Bridgeland stability conditions, Fargues--Fontaine curves, and Banach--Colmez spaces. We also discuss the relation between these topics.
\subsection{Bridgeland stability conditions}

	In this subsection, we review the definition and some basic results of Bridgeland stability conditions (see \cite{beilinson1982faisceaux}, \cite{polishchuk2007constant}, \cite{bridgeland2008stability}, \cite{kontsevich2008stability} and \cite{baye2011bridgeland}).
		
		\begin{definition}
			Let $\mathcal{D}$ be a triangulated category. A \textit{$t$-structure} on $\mathcal{D}$ is a pair of full subcategories $(\mathcal{D}^{\leq 0},\mathcal{D}^{\geq 0})$ satisfying the conditions (i)-(iii) below. We denote $\mathcal{D}^{\leq n}=\mathcal{D}^{\leq 0}[-n]$, $\mathcal{D}^{\geq n}=\mathcal{D}^{\geq 0}[-n]$ for every $n\in\mathbb{Z}$. Then the conditions are:
			
			(i) $\Hom(X,Y)=0$ for every $X\in\mathcal{D}^{\leq 0}$ and $Y\in\mathcal{D}^{\geq 1}$;
			
			(ii) $\mathcal{D}^{\leq -1}\subset \mathcal{D}^{\leq 0}$ and  $\mathcal{D}^{\geq 1}\subset \mathcal{D}^{\geq 0}$;
			
			(iii) Every object $X\in\mathcal{D}$ fits into an exact triangle 
			
			$$\tau^{\leq 0}X\rightarrow X\rightarrow \tau^{\geq 1}X\rightarrow \cdots$$ with $\tau^{\leq 0}X\in\mathcal{D}^{\leq 0}$, $\tau^{\geq 1}X\in\mathcal{D}^{\geq 1}$.
			
		It is not hard to see that the objects $\tau^{\leq 0}X$ and $\tau^{\geq 0}X$ are defined functorially. The functors $\tau^{\leq 0}$ and $\tau^{\geq 1}\coloneqq \tau^{\geq 0}[-1]$ are called truncation functors (see e.g. \cite[Proposition 8.1.8]{hotta2007d} for the basic properties of truncation functors). 	The heart of a $t$-structure is $\mathcal{A}=\mathcal{D}^{\leq 0}\cap\mathcal{D}^{\geq 0}$. It is an abelian category (see e.g. \cite[Theorem 8.1.9]{hotta2007d}). We have the associated cohomology functors defined by $H^0=\tau^{\leq 0}\tau^{\geq 0}:\mathcal{D}\rightarrow\mathcal{A}$, $H^i(X)\coloneqq H^0(X[i])$. 
		\end{definition}
		Motivated by work in string theory (see e.g. \cite{douglas2002dirichlet}), the notion of t-structures was refined by Bridgeland to the notion of slicings (see Definition \ref{slicing}) in \cite{bridgeland2007stability}. Combining a slicing and the notion of a central charge from physics, Bridgeland coined the following definition of stability conditions.

		\begin{definition}\label{first WSC}
			A \textit{stability condition} on $\mathcal{D}$ consists of a pair $(\mathcal{P},Z)$, where $\mathcal{P}$ is a slicing and $Z:K_0(\mathcal{D})\rightarrow\mathbb{C}$ is a group homomorphism such that the following conditions are satisfied:
			\par

			(a) If $0\neq E\in\mathcal{P}(\phi)$ then $Z(E)=m(E)exp(i\pi\phi)$ for some $m(E)\in \mathbb{R}_{> 0}$.
			
			\par 
			
			(b) (Support property) The group homomorphism $Z$ factors as $$K_0(\mathcal{D})\xrightarrow{v} \Lambda\xrightarrow{g} \mathbb{C},$$ where $\Lambda$ is a finite rank lattice, $v$ and $g$ are both group homomorphisms, and there exists a quadratic form $Q$ on $\Lambda_\mathbb{R}\coloneqq\Lambda\otimes\mathbb{R}$ such that $Q|_{\ker(g)}$ is negative definite, and $Q(v(E))\geq 0$, for any semistable object $E$. 
			
		\end{definition}

		\begin{remark}\label{rotation}
			The group homomorphism $Z$ is usually called the central charge in physics literature, we also adopt this terminology in this paper.
			
			If we only require $m(E)$ to be nonnegative in (a), the pair $(\mathcal{P},Z)$ is called a weak stability condition (see \cite[Definition 14.1]{Stabilityconditionsinfamilies}). 
   \end{remark}
   	\begin{remark}
   	    Bridgeland proved that the space of stability conditions $\Stab(\mathcal{D})$ has the structure of a complex manifold. It carries a right action of the group $\widetilde{GL}^+(2,\mathbb{R})$, the universal cover of $GL^+(2,\mathbb{R})$, and a left action of the group $Aut(\mathcal{D})$ of exact autoequivalences of $\mathcal{D}$; these actions commute (see \cite[Lemma 2.2]{bridgeland2008stability}). Note that the group action of $\widetilde{GL}^+(2,\mathbb{R})$ does not change the class of semistable objects.
   	\end{remark}

		There is an equivalent way of defining a stability condition, which is more closely related to the classical Mumford stability condition for vector bundles on curves. To introduce such a definition, we need to define a stability function $Z$ on an abelian category $\mathcal{A}$.
		\begin{definition}\label{Definition of stability function} Let $\mathcal{A}$ be an abelian category. We call a group homomorphism $Z:K_0(\mathcal{A})\rightarrow \mathbb{C}$ a \textit{stability function} on $\mathcal{A}$ if, for any nonzero object $E\in \mathcal{A}$, we have $Im(Z(E))\geq0$, with $Im(Z(E))=0 \implies  Re(Z(E))<0$. 
			
		\end{definition}
		
		\begin{definition}\label{another definition of BSC}
			A \textit{stability condition} on $\mathcal{D}$ is a pair $\sigma=(\mathcal{A},Z)$ consisting of the heart $\mathcal{A}\subset\mathcal{D}$ of a bounded t-structure and a stability function $Z:K_0(\mathcal{A})\rightarrow \mathbb{C}$ such that (a) and (b) below are satisfied:
			
			\par
			(a) (HN-filtration) The function $Z$ allows us to define a slope for any object $E$ in the heart $\mathcal{A}$ by
			
			$$\mu_{\sigma}(E):=\begin{cases} -\frac{Re(Z(E))}{Im(Z(E))}\ &\text{if} \  Im(Z(E))> 0,\\ +\infty &\text{otherwise.} \end{cases}$$

			The slope function gives a notion of stability: A nonzero object $E\in \mathcal{A}$ is $\sigma$-semistable if for every proper subobject $F$, we have $\mu_{\sigma}(F)\leq  \mu_{\sigma}(E)$. Moreover, if $\mu_{\sigma}(F)<  \mu_{\sigma}(E)$ holds for any proper subobject $F\subset E$, we say that $E$ is $\sigma$-stable.
			
			We require any nonzero object $E$ of $\mathcal{A}$ to have a Harder--Narasimhan filtration whose factors are $\sigma$-semistable, i.e., there exists a unique filtration $$0=E_0\subset E_1 \subset E_2\subset \cdots \subset E_{m-1} \subset E_m=E$$ such that $E_i/E_{i-1}$ is  $\sigma$-semistable and $\mu_{\sigma}(E_i/E_{i-1})>\mu_{\sigma}(E_{i+1}/E_i)$ for any $1\leq i< m$. The quotient $E_{i}/E_{i-1}$ is called the $i$-th HN factor of $E$.
			
			(b) (Support property) Equivalently, as in Definition \ref{first WSC}, the central charge $Z$ factors as $K_0(\mathcal{A})\xrightarrow{v} \Lambda\xrightarrow{g} \mathbb{C}$. There also exists a quadratic form $Q$ on $\Lambda_{\mathbb{R}}$ such that $Q|_{\ker(g)}$ is negative definite and $Q(v(E))\geq 0$ for any $\sigma$-semistable object $E\in\mathcal{A}$.
		\end{definition}
		
		\begin{remark}\label{remark: equivalence of definitions}
			  The equivalence of these two definitions is given by setting $\mathcal{P}(\phi)$ to be the full subcategory consisting of $\sigma$-semistable objects in $\mathcal{A}$ of slope $-\cot(\pi\phi)$  for $\phi\in(0,1]$, so here $\mathcal{A}=\mathcal{P}(0,1]$. There is another equivalent way to define a stability condition if we take $\mathcal{A}$ to be $\mathcal{P}[0,1)$, we leave this definition to the readers.
     
     If the imaginary part of $Z$ is discrete in $\mathbb{R}$, and $\mathcal{A}$ is Noetherian, then condition (a) is satisfied automatically (see e.g. \cite[Proposition 2.8]{Filtrationsandtorsionpairs}). Note that if the induced $\mathbb{R}$-linear morphism $g_{\mathbb{R}}:\Lambda_{\mathbb{R}}\rightarrow \mathbb{C}$ is injective, the support property in \ref{another definition of BSC}.(b) would hold automatically.
			
		\end{remark}

  Given a Bridgeland stability condition $\sigma=(\mathcal{P},Z)$ on a triangulated category $\mathcal{D}$, we have the following proposition.
	
	\begin{prop}\label{proposition of noetherian and aritinian}
		Assume that the central charge $Z$ of $\sigma$ satisfies the following two conditions: 
		\begin{enumerate}
			\item the image of $Z$ is discrete; 
			\item the image of the imaginary part of $Z$ is discrete.
		\end{enumerate}
		Then the heart $\mathcal{P}(0,1]$ is Noetherian, and the heart $\mathcal{P}[0,1)$ is Artinian.
	\end{prop}
	
	\begin{proof}
		The Noetherianity of $\mathcal{P}(0,1]$ is \cite[Proposition 5.0.1]{APsheaves}, and the case of $\mathcal{P}[0,1)$ is a completely dual statement. We include its proof for readers' convenience. 
		
		We consider a decreasing sequence $$\cdots \subset F_2\subset F_1\subset E$$ of nonzero objects in the heart $\mathcal{P}[0,1)$.
		
		Firstly, if $E\in\mathcal{P}(0)$, then the sequence stabilizes after finitely many step. Indeed, since $$Im(Z(F_i))=-ImZ(E/F_i)$$ and both terms are nonnegative, we have $Im(Z(F_i))=0$ and thus $F_i\in\mathcal{P}(0)$. By the discreteness assumption (1), we have that there exists $N>0$ such that $Z(F_i)$ is constant for all $i\geq N$, hence $F_i$ becomes constant. 
		
		By the discreteness assumption (2), we may assume that $Im(Z(F_i))$ becomes constant for $i$ large enough. By the discussion of the case when $E\in\mathcal{P}(0)$, we may assume that $Im(Z(F_i))$ is a positive constant for $i$ large enough. Assume that sequence does not stabilize after finitely many steps. Then for any $\phi \in[0,1)$, we can find $N$ large enough such that for any $n>N$, we have $\mu(F_n)>-\cot(\pi\phi)$ (this also uses the discreteness assumption (1)).  If we take $\phi=\phi_{max}(E)$, i.e. the phase of the first Harder--Narasimhan factor of $E$, then we get that $F_n$ can not be a subobject of $E$ as the first Harder--Narasimhan factor of $F_n$ maps trivially to $E$ by phase inequality, which is a contradiction.
	\end{proof}

\subsection{Fargues--Fontaine curves}
\label{subsection:Fargues--Fontaine curve} Let $\mathbb{Q}_p$ be the $p$-adic completion of $\mathbb{Q}$, and we fix $C$ a complete algebraically closed field extension of $\mathbb{Q}_p$. $C$ is a non-Archimedean field with ring of integers $\mathcal{O}_C$ and residue field $k$. $C$ is a perfectoid field in the sense of \cite[Definition 1.2]{scholze-perfectoidspacesIHES}, we let $C^\flat$ be the tilt. In particular, $C^\flat$ is a non-Archimedean over $\mathbb{F}_p$ with ring of integers $\mathcal{O}_C^\flat$. We write $A_{\mathrm{inf}}=A_{\mathrm{inf}}(\mathcal{O}_C)$ for the infinitesimal period ring of Fontaine. There is a canonical surjection $\vartheta\colon A_{\mathrm{inf}} \to \mathcal{O}_C$ such that $\ker(\vartheta)$ is principal, generated by $\xi$. This map will be referred to Fontaine's $\vartheta$ map. Define $A_{\mathrm{cris}}$ as the $p$-complete PD envelope of $A_{\mathrm{inf}}$ with respect to $\ker(\vartheta)$, and define $B_{\mathrm{cris}}^+=A_{\mathrm{cris}}[1/p]$. The map $\vartheta$ extends to a map $\vartheta: B_{\mathrm{cris}}^+ \to C$. We define $B_{\mathrm{dR}}^+$ to be $\ker(\vartheta)$-adic completion of $B_{\mathrm{cris}}^+$. Let $\Breve{\mathbb{Q}}_p=W(\overline{\mathbb{F}}_p)[1/p]$, which carries a Frobenius endomorphism $\varphi$ lifting the $p$-th power map modulo $p$. Since $k$ is algebraically closed, there is a canonical inclusion $\overline{\mathbb{F}}_p \to k$, which induces a map $\Breve{\mathbb{Q}}_p \to B_{\mathrm{cris}}^+$.

This subsection will review the definition of (absolute) Fargues--Fontaine curve in \cite{farguesfontaine-courbes} and the construction of vector bundles over the Fargues--Fontaine curve from isocrystals. 

Recall in \cite[Chapter 6]{farguesfontaine-courbes}, given $C^\flat$ and a complete discrete valuation field $E$ with perfect residue field, and $\pi$ being a uniformizer of $E$, there is a curve $X_{C^\flat,E}$ defined as the projective scheme associated to the graded algebra $P_{C^\flat,E}$ with 
\[
P_{C^\flat,E}\coloneqq \oplus_{i\geq 0} ({B}_{\mathrm{cris}}^+)^{\varphi=\pi^i}.
\]
In particular, if we let $E=\mathbb{Q}_p$, and we define $X_{FF}=X_{C^\flat,\mathbb{Q}_p}$ to be the absolute Fargues--Fontaine curve. Although all the results in this paper hold with slight changes for any $X_{C^{\flat},E}$, we will only state our results for $X_{FF}$ for simplicity. We summarize the properties of $X_{FF}$.
\begin{theorem}[{\cite[\S 6.5]{farguesfontaine-courbes}}]\label{thm:FFmain}
    \begin{enumerate}
        \item $X_{FF}$ is a complete curve with field of definition equal to $\mathbb{Q}_p$ in the sense of \cite[\S 5.1]{farguesfontaine-courbes};
        \item the map $\theta: B_{\mathrm{cris}}^+ \to C$ defines a closed point $\infty$ of $X_{FF}$, whose residue field is isomorphic to $C$;
        \item $\widehat{\mathcal{O}}_{X_{FF},\infty}$ is isomorphic to $B_{\mathrm{dR}}^+$;
        \item $B_e\coloneqq H^0(X_{FF}\backslash\{\infty\},\mathcal{O}_{X_{FF}})$ is a principal ideal domain; 
        \item Residue fields at closed points of $X_{FF}$ are untilts of $C^\flat$ in characteristic $0$. This defines a bijection from closed points of $X_{FF}$ to isomorphism classes of untilts of $C^\flat$ in characteristic $0$.
    \end{enumerate}
\end{theorem}
Thus, $X_{FF}$ behaves like a Riemann surface (so-called generalized Riemann sphere). In particular, one can define the notion of semistable vector bundles, and any vector bundle on $X_{FF}$ admit a  unique Harder--Narasimhan filtration. After recalling the ``Dieudonn\'e-Manin classification'' of vector bundles, we will come back to this.

Let $\mathrm{Bun}_{X_{FF}}$ denote the category of vector bundles over $X_{FF}$. We recall the following classification theorem of vector bundles over $X_{FF}$ after Fargues--Fontaine. Let $\varphi\text{-}\mathrm{Mod}_{\Breve{\mathbb{Q}}_p}$ be the category of isocrystals over $\overline{\mathbb{F}_p}$, i.e., the category of pairs $(D,\varphi_D)$ with $D$ a finite-dimensional vector space over $\Breve{\mathbb{Q}}_p$ and $\varphi_D: D\to D$ a $\varphi$-semilinear endomorphism of $D$. There is a well-defined functor $\mathcal{E}\colon \varphi\text{-}\mathrm{Mod}_{\Breve{\mathbb{Q}}_p} \to \mathrm{Bun}_{X_{FF}}$ defined as
\[
\mathcal{E}((D,\varphi_D))\coloneqq \widetilde{\oplus_{i\geq 0} (D\otimes_{\Breve{\mathbb{Q}}_p}{B}_{\mathrm{cris}}^+)^{\varphi=p^i}}
\]
i.e., the vector bundle associated with the graded module $M((D,\varphi_D))\coloneqq\oplus_{i\geq 0} (D\otimes_{\Breve{\mathbb{Q}}_p}{B}_{\mathrm{cris}}^+)^{\varphi=p^i}$. 

\begin{theorem}[{\cite[Chapter 8]{farguesfontaine-courbes}}]
The well-defined functor $\mathcal{E}$ is faithful, essentially surjective, and not full. This map induces a bijection between the isomorphism classes in both categories.    
\end{theorem}

Moreover, using the above theorem, one can interpret the slope of a vector bundle in the following way.
\begin{prop}
\begin{enumerate}
    \item A vector bundle $\mathcal{V}$ is semistable of slope $\lambda\in \mathbb{Q}$ if and only if $\mathcal{V}$ is in the essential image of the restriction of $\mathcal{E}$ to the subcategory of isocrystals of slope $-\lambda$;
    \item $\mathcal{E}$ induces an equivalence of the category of semisimple isocrystals over $\overline{\mathbb{F}_p}$ of slope $-\lambda$ and the category of semistable vector bundles over $X_{FF}$ of slope $\lambda$;
    \item Let $\mathcal{O}(\lambda)$ be the image of the simple isocrystal of slope $-\lambda$ under $\mathcal{E}$, assume that $\lambda=\frac{d}{h}$ is in lowest terms with $h>0$, then $\mathcal{O}(\lambda)$ is of rank $h$ and degree $d$;
    \item Every vector bundle is (non-canonically) isomorphic to a finite direct sum of $\mathcal{O}(\lambda_i)$ with $\lambda_i\in \mathbb{Q}$.
\end{enumerate}
\end{prop}

The functor $\mathcal{E}$ is not full. The category of isocrystals is semi-simple, but the $\Hom$ and $\Ext^1$ spaces between (semistable) vector bundles over $X_{FF}$ are usually quite large. 

\begin{prop}\label{prop:homandext1}
\begin{enumerate}
    \item We have $\Hom(\mathcal{O}(\lambda),\mathcal{O}(\mu))=0$ if $\lambda > \mu$;
    \item $\Ext^1(\mathcal{O}(\lambda),\mathcal{O}(\mu))=0$ if $\lambda \leq \mu$;
    \item $\Ext^i$ vanishes for $i \neq 0, 1$ since $X_{FF}$ is of dimension $1$.
\end{enumerate}
\end{prop}

We can extend the slope function of the vector bundles over $X_{FF}$ to the category of coherent sheaves by letting torsion sheaves have slope $\infty$. This can be modified to the following definition.

\begin{definition}
Let $\mathcal{E}$ be a coherent sheaf on $X_{FF}$, define $v(\mathcal{E}) \in \mathbb{Z}^2$ to be the \emph{character class} of $\mathcal{E}$ which corresponds to the integral vector $(\deg(\mathcal{E}), \rank(\mathcal{E}))$. When $\mathcal{E}$ is torsion, we define the degree of $\mathcal{E}$ to be the sum of $\mathrm{dim}_{\kappa_x} \mathcal{E}_x$ over closed points of $X_{FF}$. Note that $\mathcal{E}$ is supported on finitely many closed points of $X_{FF}$, whose residue fields $\kappa_x$ are untilts of $C^\flat$ by Theorem~\ref{thm:FFmain}.
\end{definition}

Let $x$ be a closed point of $X_{FF}$, let $\iota_x\colon \mathrm{Spec}(\kappa_x) \to X_{FF}$ be the closed immersion. For any $k>0$, the torsion sheaf $\iota_{x,\ast}B_{\mathrm{dR}}^+(\kappa_x)/t^k B_{\mathrm{dR}}^+(\kappa_x)$ has a resolution
$$
0\to \mathcal{O} \to \mathcal{O}(k) \to \iota_{x,\ast}B_{\mathrm{dR}}^+/t^k B_{\mathrm{dR}}^+ \to 0
$$
which implies that the degree function defined above is additive. The exactness of the above sequence can be found in \cite[Lemma 7.3]{LeBrasresult}. Recall that we have fixed a characteristic $0$ untilt $C$ of $C^\flat$, it defines the closed point $\infty$ on $X_{FF}$. We will use $\mathcal{O}(\frac{k}{0})$ to denote the torsion sheaf $\iota_{\infty,\ast}B_{\mathrm{dR}}^+/t^k B_{\mathrm{dR}}^+$ unless otherwise specified.

\subsection{Banach--Colmez spaces}
%\textcolor{blue}{Qingyuan: I temporarily added the following subsection; to revise and add more details and references later. Feel free to delete/add/correct.}

\subsubsection{Colmez's definition}
We begin by reviewing Colmez's original definition \cite{BanachColmezspaces} of Banach--Colmez spaces, which he refers to as ``finite-dimensional Banach Spaces." 
%First, we recall several terminologies.  
A \emph{sympathetic $C$-algebra} $\Lambda$ is a Banach $C$-algebra equipped with the spectral norm $\Vert \cdot \Vert$ such that $\lambda \mapsto \lambda^p$ is surjective on $\{\lambda \mid \Vert \lambda -1 \Vert <1\}$ (\cite[\S 4]{BanachColmezspaces}). By definition, sympathetic algebras are perfectoid algebras. 
A \emph{Banach space} is a functor
%	$$(\text{Sympathetic $C$-Algebras} ) \to (\text{$\mathbb{Q}_p$-Banach spaces}).$$	
from the category of sympathetic $C$-algebras to the category of $\mathbb{Q}_p$-Banach spaces. A morphism of Banach spaces is a natural transformation of functors. A sequence of Banach spaces is called \emph{exact} if it remains exact when evaluated on any sympathetic algebra $\Lambda$. 
				
\begin{example} 
\label{example.BC.Q_p.C}
\begin{enumerate}
	\item For any finite-dimensional $\mathbb{Q}_p$-vector space $V$, the constant sheaf $\Lambda \mapsto V$ defines a Banach space, denoted by $\underline{V}$.
	\item For any finite-dimensional $C$-vector space $W$, the functor $\Lambda \mapsto W\otimes \Lambda$ defines a Banach space, denoted by $W \otimes \mathcal{O}_C$. If $W = C$, we simply denote $W \otimes \mathcal{O}_C$ by $\mathcal{O}_C$.
\end{enumerate}
\end{example}	

A Banach--Colmez space is a Banach space that is ``equal to $W \otimes \mathcal{O}_C$ up to finite-dimensional $\mathbb{Q}_p$-vector spaces". More precisely:

\begin{definition}
An \emph{effective Banach--Colmez space} is a Banach space $\mathbb{Y}$ that is an extension of $W \otimes \mathcal{O}_C$ by some $\underline{V}$, where $V$ (resp. $W$) is a finite-dimensional $\mathbb{Q}_p$-vector space (resp. $C$-vector space). A {\em Banach--Colmez space} is a Banach space $\mathbb{W}$ that is the quotient of an effective Banach--Colmez space $\mathbb{Y}$ by some $\underline{V'}$, where $V'$ is a finite-dimensional $\mathbb{Q}_p$-vector space. In particular, we have exact sequences
	\begin{equation}\label{eqn:presentation:BC}
	 0 \to \underline{V} \to \mathbb{Y} \to W \otimes \mathcal{O}_C \to 0 \qquad 0 \to \underline{V'} \to \mathbb{Y} \to \mathbb{W} \to 0.
\end{equation}
Therefore, we can regard $\mathbb{W}$ as obtained from $W \otimes \mathcal{O}_C$ by ``adding $\underline{V}$ and quotient out $\underline{V'}$".
We will refer to \eqref{eqn:presentation:BC} as \emph{a presentation of $\mathbb{W}$}. We let $\mathcal{BC}$ denote the category of Banach--Colmez spaces.
\end{definition}

\begin{theorem}[Colmez {\cite[\S 6]{BanachColmezspaces}}]
\label{thm:Colmez}
\begin{enumerate}
	\item 
	\label{thm:Colmez-1}
	The category $\mathcal{BC}$ of Banach--Colmez spaces is an abelian category. 
	\item
	\label{thm:Colmez-2}
	 For any $\mathbb{W} \in \mathcal{BC}$ with a presentation as \eqref{eqn:presentation:BC}, the pair of integers 
	$$\operatorname{Dim} \mathbb{W} \coloneqq (\dim \mathbb{W}, \operatorname{ht} \, \mathbb{W}) \coloneqq (\dim_C (W), \dim_{\mathbb{Q}_p} (V) - \dim_{\mathbb{Q}_p} (V'))
	$$
does not depend on the presentation of $\mathbb{W}$ and defines an additive functor
	$$\operatorname{Dim}=(\dim, \operatorname{ht}) \colon \mathcal{BC} \to \mathbb{N} \times \mathbb{Z}.$$
	\item 
	\label{thm:Colmez-3}
	 If $f\colon \mathbb{W}_1 \to \mathbb{W}_2$ is a morphism of Banach--Colmez spaces, then 
			\begin{align*}
\operatorname{Dim} (\mathbb{W}_1) &= \operatorname{Dim} (\ker(f)) + \operatorname{Dim} (\mathrm{im}(f)) \\
\operatorname{Dim} (\mathbb{W}_2) &= \operatorname{Dim} (\operatorname{coker}(f)) + \operatorname{Dim} (\mathrm{im}(f)).
\end{align*}
	\item
	\label{thm:Colmez-4}
	 If $\dim \mathbb{W} = 0$, then $\operatorname{ht} \, \mathbb{W} \ge 0$.
\end{enumerate}
\end{theorem}

\begin{proof}
Claims \eqref{thm:Colmez-1} and \eqref{thm:Colmez-3} follow from \cite[Theorem 0.4, Propositions 6.16 and 6.22]{BanachColmezspaces}. Claim \eqref{thm:Colmez-2} is \cite[Corollary 6.17]{BanachColmezspaces}. Claim \eqref{thm:Colmez-4} follows from the definition. 
\end{proof}

\begin{remark}
The functor of evaluation at $C$,
	$$\mathcal{BC} \to (\text{$\mathbb{Q}_p$-Banach Spaces}) \quad \mathbb{W} \mapsto \mathbb{W}(C),$$
is exact and conservative (exactness follows from the definition of exact sequences of Banach spaces, and conservativeness follows from checking on effective Banach--Colmez spaces). 
Therefore, exact sequences and isomorphisms of Banach--Colmez spaces can be detected on their underlying $C$-points, and we may think of a Banach--Colmez space $\mathbb{W}$ as the $\mathbb{Q}_p$-Banach space $\mathbb{W}(C)$ with additional ``analytical" structure. 
\end{remark}

\subsubsection{Alternative characterizations}
This subsection reviews Le Bras' characterization (\cite{LeBrasresult}) of Banach--Colmez spaces as pro-{\'e}tale sheaves or $v$-sheaves on the category $\mathrm{Perf}_C$ of perfectoid spaces over $\mathrm{Spa}(C)$.

We refer to \cite{scholze2017etale} or \cite[\S 2.1]{LeBrasresult} for the definitions of pro-\'{e}tale and $v$-topology; here, we only review the basic notations. Let $K = C$ or $C^\flat$. We denote the big pro-\'{e}tale site of $K$ by $\mathrm{Perf}_{K, \text{pro-\'{e}t}}$ and the $v$-site of $K$ by $\mathrm{Perf}_{K, v}$. A theorem of Scholze states that there are canonical equivalences $\mathrm{Perf}_{C, \text{pro-\'{e}t}} \simeq \mathrm{Perf}_{C^\flat, \text{pro-\'{e}t}}$ and $\mathrm{Perf}_{C, v} \simeq \mathrm{Perf}_{C^\flat, v}$ (see \cite{scholze2017etale}  or \cite[Theorem 2.7]{LeBrasresult}). Let $\mathcal{S}\mathrm{hv}_{\text{pro-\'{e}t}}(C; \mathbb{Q}_p)$ (resp. $\mathcal{S}\mathrm{hv}_{v}(C; \mathbb{Q}_p)$) denote the abelian category of sheaves of $\mathbb{Q}_p$-vector spaces on $\mathrm{Perf}_{C,\text{pro-\'{e}t}}$ (resp. on $\mathrm{Perf}_{C,v}$) and let $D_{\text{pro-\'{e}t}}(C; \mathbb{Q}_p)$ (resp. $D_{v}(C; \mathbb{Q}_p)$) denote its $\infty$-categorical derived category.

\begin{example}\label{example.BC.Q_p.C:v-sheaves}
\begin{enumerate}
	\item The constant presheaf $\underline{\mathbb{Q}_p}$ is a sheaf for the $v$-topology (and hence also for the pro-{\'e}tale topology).
	\item The presheaf $\mathcal{O}_C$ which associates each $S \in \mathrm{Perf}_{C}$ to $\mathcal{O}_S(S)$ is a sheaf for the $v$-topology (and hence also for the pro-{\'e}tale topology); see \cite[Theorem 8.7]{scholze2017etale}.
\end{enumerate}
\end{example}

\begin{definition}[{\cite[Definition 2.11]{LeBrasresult}}]
Let $\mathcal{BC}' \subseteq \mathcal{S}\mathrm{hv}_{\text{pro-{\'e}t}}(C; \mathbb{Q}_p)$ (resp. $\mathcal{BC}'' \subseteq \mathcal{S}\mathrm{hv}_{v}(C; \mathbb{Q}_p)$) denote the smallest abelian subcategory that is stable under extensions and contains the sheaves $\underline{\mathbb{Q}_p}$ and $\mathbb{G}_a$.
\end{definition}

\begin{prop}[{\cite[Proposition 7.11 \& Remark 2.12]{LeBrasresult}}]
\label{prop:BC:3equivalent.definitions}
There are canonical equivalences of abelian categories:
    $$
    \mathcal{BC} \xleftarrow{\sim} \mathcal{BC}' \xleftarrow{\sim} \mathcal{BC}''.
    $$
Here, the first equivalence is induced by the forgetful functor sending an object in $\mathcal{BC}'$ to the underlying functor on the category of sympathetic $C$-algebras. The second equivalence is induced by the pullback functor $\nu^*$, where $\nu \colon \mathrm{Perf}_{C, v} \to \mathrm{Perf}_{C, \text{pro-\'{e}t}}$ is the morphism of sites.  
\end{prop}

In light of this proposition, we can identify $\mathcal{BC}$ with the full subcategory $\mathcal{BC}''$ of $\mathcal{S}\mathrm{hv}_{v}(C; \mathbb{Q}_p)$ and regard Banach--Colmez spaces as $v$-sheaves of $\mathbb{Q}_p$-vector spaces on $\mathrm{Perf}_C$. Note that under this equivalence, the Banach spaces $\underline{\mathbb{Q}_p}$ and $\mathcal{O}_C$ of Example \ref{example.BC.Q_p.C} correspond canonically to the $v$-sheaves $\underline{\mathbb{Q}_p}$ and $\mathcal{O}_C$ of Example \ref{example.BC.Q_p.C:v-sheaves}, respectively.

\subsubsection{Cohomologies of coherent sheaves on Fargues--Fontaine curves}
%\textcolor{blue}{Qingyuan: I typed a preliminary presentation of the BC version of cohomologies and Hom-spaces, in terms of the adic version of relative FF curves and the $v$-topology (see the stuff before Prop 2.16). However, I am not entirely familiar with the references for this part, so please make any necessary corrections.}

The construction of the Fargues--Fontaine curve $X_{FF}=X_{C^\flat, \mathbb{Q}_p}$ presented in \S \ref{subsection:Fargues--Fontaine curve} depends on a perfectoid field $C^\flat$ (and choosing an untilt $C$ of $C^\flat$ corresponds to specifying a closed point $\infty \in |X_{FF}|$). In general, the perfectoid field $C^\flat$ can be replaced by any perfectoid space $S \in \mathrm{Perf}_{\mathbb{F}_p}$. There is a functorial construction of a ``relative Fargues--Fontaine curve" $\mathcal{X}_S$ over $S$, which is an adic space over $\mathrm{Spa}(\mathbb{Q}_p)$; see \cite[\S II.1.2]{fargues2021geometrization}. If $S=\mathrm{Spa}(C^\flat)$, then the adic space $\mathcal{X}_{C^\flat}$ is an ``analytification" of the Fargues--Fontaine curve $X_{FF}$. There is a \emph{GAGA} theorem that implies that the natural morphism $\mathcal{X}_{C^\flat} \to X_{FF}$ of locally ringed topological spaces induces an equivalence of categories  $\mathrm{Bun}_{\mathcal{X}_{C^\flat}} \simeq \mathrm{Bun}_{X_{FF}}$ (\cite[Theorem 6.3.9]{kedlaya2015relative} and \cite[Proposition II.2.7]{fargues2021geometrization}). Moreover, the natural map $|\mathcal{X}_{C^\flat}| \to |X_{FF}|$ induces a bijection between classical points of $\mathcal{X}_{C^\flat}$ and the closed points of $X_{FF}$ (\cite[Proposition II.2.9]{fargues2021geometrization}).

%\textcolor{blue}{Qingyuan: Update: I just realized we might not need the GAGA for Coh, as long as the GAGA equivalence for Bun is exact.
%Original questions: I think we can argue $\mathrm{Coh}(\mathcal{X}_{C^\flat}) \simeq \mathrm{Coh}(X_{FF})$ as follows? -- for vector bundles, this is proved in \cite[Theorem 6.3.9]{kedlaya2015relative} and \cite[Proposition II.2.7]{fargues2021geometrization}; for torsion sheaves, it follows from \cite[Proposition II.2.9]{fargues2021geometrization}. And is there a direct reference for this?}
%\textcolor{red}{Heng: there is a version of GAGA for Coh in Kedlaya--Liu II, which I could find the exact reference number if needed. On the other hand, there is no good theory of coherent sheaves on adic spaces without Noetherian assumptions(the adic FF curve is not noetherian). That is the main point that in KLII, they define something called fpd sheaves and pseudo-coherent sheaves.}
%\textcolor{blue}{Qingyuan: I see, thanks!}

\begin{prop}
For any perfectoid space $S$ of characteristic $p$ and any perfect complex $\mathcal{E}$ on $\mathcal{X}_S$,  the functor 
			$$(T \in \mathrm{Perf}_{S}) \mapsto R\Gamma(\mathcal{X}_T; \mathcal{E}|_{\mathcal{X}_T})$$
is a $v$-sheaf of complexes.
\end{prop}

\begin{proof}
In the case where $\mathcal{E}$ is a vector bundle, this is proved in \cite[Lemma 17.1.8, Proposition 19.5.3]{scholze2020berkeley}; see also \cite[Proposition II.2.1]{fargues2021geometrization}. The general case of perfect complexes is proved using a similar strategy in \cite[Proposition 2.4]{anschutz2021fourier}. %\textcolor{red}{Heng: should add more details.}
\end{proof}

Consequently, taking $S=C^\flat$ and combining with the GAGA equivalence $\mathrm{Bun}_{\mathcal{X}_{C^\flat}} \simeq \mathrm{Bun}_{X_{FF}}$, the fact that every bounded coherent complex on $X_{FF}$ is a perfect complex, and Scholze's equivalence $\mathrm{Perf}_{C;v} \simeq \mathrm{Perf}_{C^\flat, v}$ (implying $D_{v}(C; \mathbb{Q}_p) \simeq D_{v}(C^\flat; \mathbb{Q}_p)$), we obtain a functor 
	$$R \tau_* \colon D^b(\mathrm{Coh}(X_{FF})) \to D_{v}(C; \mathbb{Q}_p).$$
In particular, for any bounded coherent complex $\mathcal{E} \in  D^b(\mathrm{Coh}(X_{FF}))$ and $i \in \mathbb{Z}$, $R^i \tau_*(\mathcal{E})$ is a $v$-sheaf of $\mathbb{Q}_p$-vector spaces on $\mathrm{Perf}_C$.

\begin{definition}
For any coherent complex $\mathcal{E}$ on $X_{FF}$ and $i \in \mathbb{Z}$, let 
	$$\mathbb{H}^i(X_{FF}; \mathcal{E}) \coloneqq R^i \tau_*(\mathcal{E}) \in \mathcal{S}\mathrm{hv}_{v}(C;\mathbb{Q}_p)$$
be the $v$-sheaf of $\mathbb{Q}_p$-vector space which sends each $T \in \mathrm{Perf}_{C^\flat}$ to $R^i \Gamma(\mathcal{X}_T, \mathcal{E}|_T)$.
More generally, for any pair of coherent complexes $\mathcal{E}, \mathcal{F}$ on $X_{FF}$ and any $i \in \mathbb{Z}$, let 
    $$\mathbb{E}\mathrm{xt}^i(\mathcal{E},\mathcal{F}) \coloneqq \mathbb{H}^i(X_{FF}; \mathcal{E}^\vee \otimes^{\mathbb{L}} \mathcal{F}) = R^i \tau_*(R \mathcal{H}\mathrm{om}_{X_{FF}}(\mathcal{E},\mathcal{F})).$$
In the case where $i=0$, we also write $\mathbb{H}\mathrm{om}(\mathcal{E},\mathcal{F}) \coloneqq \mathbb{E}\mathrm{xt}^0(\mathcal{E},\mathcal{F})$.
\end{definition}

We present an upgraded version of the cohomologies of coherent sheaves on Fargues–Fontaine curves in terms of Banach--Colmez spaces. Let $\mathbb{U}_{d,h} =(\mathbb{B}_{\mathrm{cris}}^+)^{\varphi^h=\pi^d}$ (if $d \ge 0$), $\mathbb{V}_{(-d,h)} = \mathbb{B}_{\mathrm{dR}}^+/(t^{-d} \mathbb{B}_{\mathrm{dR}}^+ + \underline{\mathbb{Q}_{p^h}})$ (if $d<0$) and $\mathbb{B}_m = \mathbb{B}_{\mathrm{dR}}^+/t^{m} \mathbb{B}_{\mathrm{dR}}^+$ be the Banach--Colmez spaces defined in \cite[Example 2.3]{Modularvariants}, regarded as $v$-sheaves of $\mathbb{Q}_p$-vector spaces via the identification $\mathcal{BC} \simeq \mathcal{BC}''$ of Proposition \ref{prop:BC:3equivalent.definitions}.

\begin{prop} 
\label{prop:cohomologies.O(lambda).Fargues--Fontaine}
Let $\lambda = d/h \in \mathbb{Q}$, where $d \in \mathbb{Z}$, $h \in \mathbb{Z}_{>0}$, $(d,h)=1$, let $m >0$ be an integer and let $i \in \mathbb{Z}$. 
\begin{enumerate}
	\item ({\cite[Proposition 8.2.3]{farguesfontaine-courbes}}) 
	We have
		\begin{itemize}
			\item 
				$\mathbb{H}^0(X_{FF}; \mathcal{O}(\lambda)) = 
					\begin{cases}
					\mathbb{U}_{d,h} & \text{if $\lambda \ge $0}. \\
					0 & \text{if $\lambda<0$.}
					\end{cases}$
			\item
				$\mathbb{H}^1(X_{FF}; \mathcal{O}(\lambda)) = 
					\begin{cases}
					0 & \text{if $\lambda \ge $0}. \\
					\mathbb{V}_{(-d,h)} & \text{if $\lambda<0$.}
					\end{cases}.$
			\item 
			$\mathbb{H}^i(X_{FF}; \mathcal{O}(\lambda)) = 0$ if $i \ge 2$.
		\end{itemize}
	\item ({\cite[Proposition 8.2.8]{farguesfontaine-courbes}}) We have $\mathrm{End}(\mathcal{O}_X(\lambda)) = D_\lambda$.
	\item We have 
		\begin{itemize}
			\item $\mathbb{H}^0(X_{FF}, \mathcal{O}(\frac{m}{0})) = \mathbb{B}_m$.
			\item $\mathbb{H}^i(X_{FF}, \mathcal{O}(\frac{m}{0})) =0$ for $i \neq 0$.
		\end{itemize}
More generally, for any torsion coherent sheaf $\mathcal{E}$ on $X_{FF}$ of degree $m>0$, we have
		\begin{itemize}
			\item $\mathbb{H}^0(X_{FF}; \mathcal{E})$ is a Banach--Colmez space with 
				$$\operatorname{Dim} \mathbb{H}^0(X_{FF}; \mathcal{E}) = (m, 0).$$
			\item  $\mathbb{H}^i(X_{FF}; \mathcal{E}) = 0$ if $i \ne 0$.
		\end{itemize}
	\item 
	For any $\mathcal{E} \in \mathrm{Coh}(X_{FF})$ and $i \in \mathbb{Z}$, $\mathbb{H}^i(X_{FF}; \mathcal{E})$ is a Banach--Colmez space. Moreover, $\mathbb{H}^i(X_{FF}; \mathcal{E})=0$ if $i \neq 0, 1$, and  
			$$\mathrm{Dim} \, \mathbb{H}^0(X_{FF}; \mathcal{E}) - \mathrm{Dim} \, \mathbb{H}^1(X_{FF}; \mathcal{E}) = (\deg(\mathcal{E}), \operatorname{rank}(\mathcal{E})).$$ 
\end{enumerate}
\end{prop}
%\textcolor{blue}{Qingyuan: this proposition is slightly wordy, we might simplify it if needed. Also I will add the details for the proof.}
\begin{proof}
Claim (3): for any closed point $x$, let $\mathcal{E} = \iota_{x,\ast}B_{\mathrm{dR}}^+/t^m B_{\mathrm{dR}}^+$, then from the short exact sequence
	$$
0\to \mathcal{O} \to \mathcal{O}(m) \to \mathcal{E} \to 0,	$$
we obtain a long exact sequence of $v$-sheaves:
	\begin{equation*}
\begin{split}
	0 \to &\mathbb{H}^0(X_{FF}; \mathcal{O}) \to \mathbb{H}^0(X_{FF}; \mathcal{O}(m)) \to  \mathbb{H}^0(X_{FF}; \mathcal{E}) \to \\
	 & \mathbb{H}^1(X_{FF}; \mathcal{O}) \to \mathbb{H}^1(X_{FF}; \mathcal{O}(m)) \to  \mathbb{H}^1(X_{FF}; \mathcal{E}) \to \cdots. 
	\end{split}
\end{equation*}
Hence claim (3) follows from claim (1). Since $X_{FF}$ is a Dedekind scheme, every coherent sheaf splits as the direct sum of torsion coherent sheaves and vector bundles. Therefore, claim (4) follows from claims (1) and (3).
\end{proof}

We can view the formula in claim (4) as the Riemann--Roch formula for Fargues--Fontaine curves. More generally, we have the following:

 \begin{prop}\label{Riemann-Roch for Fargues-Fontaine curve.}
     Let $\mathcal{E},  \mathcal{F}\in \operatorname{Coh}(X_{FF})$ be two coherent sheaves with $v(\mathcal{E})=(p,q), v(\mathcal{F})=(r,s)$, we have $$\chi(\mathcal{E},\mathcal{F})\coloneqq \mathrm{Dim} \,\mathbb{H}\mathrm{om}(\mathcal{E},\mathcal{F})-\mathrm{Dim} \,\mathbb{E}\mathrm{xt}^1(\mathcal{E},\mathcal{F})=(qr-ps,qs). $$
 \end{prop}\begin{proof}
     Since this formula holds when $\mathcal{E},\mathcal{F}$ are vector bundles, hence it holds in general by the additivity of the dimension vector.
 \end{proof}
 
 %\textcolor{blue}{Qingyuan: That sounds great.
%We may need to provide additional arguments for torsion sheaves, as their behaviour differs slightly from the vector bundle cases. -- Update: Now the previous proposition has been modified and might be used directly here.}
 
 %\textcolor{green}{one more thing, we use $\mathcal{E}$ here for an object in $\Coh(X_{FF})$, but we use $E$ in Definition 2.13, maybe we should make the notation consistent by using $\mathcal{E}$ all the time, as we also use $E$ for an elliptic curve sometimes.}

\begin{remark}
\label{remark:homandext1}
The vanishing results in Proposition \ref{prop:homandext1} can be upgraded to the corresponding $v$-sheaf version.
%\begin{itemize}
    %\item $\mathbb{H}\mathrm{om}(\mathcal{O}(\lambda),\mathcal{O}(\mu))=0$ if $\lambda > \mu$;
    %\item $\mathbb{E}\mathrm{xt}^1(\mathcal{O}(\lambda),\mathcal{O}(\mu))=0$ if $\lambda \leq \mu$;
    %\item $\mathbb{E}\mathrm{xt}^i(\mathcal{E},\mathcal{F})=0$ for $i \neq 0, 1$ and any $\mathcal{E},  \mathcal{F}\in \operatorname{Coh}(X_{FF})$. 
%\end{itemize}
\end{remark}
%\textcolor{green}{This overlaps with proposition 2.4, maybe we just keep one?}
%\textcolor{blue}{Qingyuan: we can delete this one. Maybe mention in the previous proposition something like "there is a BC version such that these are C-points of the BC spaces..."}

\subsubsection{Serre duality}\label{subsubsection:Serre duality}
This subsection presents Serre duality for Fargues--Fontaine curves. While the results of this subsection are not directly needed for the subsequent computations in the paper, they provide a perspective that partially explains the striking similarity between the homological properties of Fargues--Fontaine curves and elliptic curves.

%\begin{prop}[{\cite[Theorem 4.1]{LeBrasresult}, \cite[Theorem 3.8]{anschutz2021fourier}}]
%Let $\mathbb{R} \mathcal{H}om_{v}$ denote the internal derived Homs in $D_{v}(C; %\mathbb{Q}_p)$. Then we have:
%\begin{align*}
	%\mathcal{R} \mathcal{H}om_v(\mathcal{O}_C, \mathcal{O}_C) &\simeq \mathcal{O}_C \oplus \mathcal{O}_C [-1] 
	%&\mathcal{R}  \mathcal{H}om_v(\mathcal{O}_C,  \underline{\mathbb{Q}_p}) &\simeq \mathcal{O}_C [-1] \\
	%\mathcal{R}  \mathcal{H}om_v(\underline{\mathbb{Q}_p}, \mathcal{O}_C) &\simeq \mathcal{O}_C
	%&\mathcal{R}  \mathcal{H}om_v(\underline{\mathbb{Q}_p}, \underline{\mathbb{Q}_p}) &\simeq \underline{\mathbb{Q}_p}.
%\end{align*}
%\end{prop}
%Applying global section functors, we obtain the following computations of derived Hom-spaces between objects in $\mathcal{S}\mathrm{hv}_{v}(C;\mathbb{Q}_p)$:\begin{align*}
	%R\Hom_v(\mathcal{O}_C, \mathcal{O}_C) &\simeq C \oplus C [-1] 
	%& R \Hom_v(\mathcal{O}_C,\underline{\mathbb{Q}_p}) &\simeq C[-1] \\
	%R\Hom_v(\underline{\mathbb{Q}_p}, \mathcal{O}_C) &\simeq C
	%& R \Hom_v(\underline{\mathbb{Q}_p}, \underline{\mathbb{Q}_p}) &\simeq \mathbb{Q}_p.
%\end{align*}
%\textcolor{blue}{While the above proposition is not necessary for stating the Serre duality, it plays a crucial role in the proof, and having this insight is psychologically beneficial for the readers. So I temporarily included it.}

\begin{theorem}[Serre Duality]
\label{thm:Serre.duality}
Let $\mathcal{E}, \mathcal{F} \in D^b(\mathrm{Coh}(X_{FF}))$. There is a canonical equivalence in $D_{v}(C; \mathbb{Q}_p)$:
	$$ R \tau_* (R \mathcal{H}\mathrm{om}_{X_{FF}}(\mathcal{E},\mathcal{F})) \simeq R \mathcal{H}\mathrm{om}_v (R \mathcal{H}\mathrm{om}_{X_{FF}}(\mathcal{F}, \mathcal{E}[1]), \underline{\mathbb{Q}_p}[1]).$$
\end{theorem}
\begin{proof}
This follows from \cite[Corollary 3.10 \& Remark 3.11]{anschutz2021fourier}.
\end{proof}

Let $\mathcal{E}, \mathcal{F} \in \mathrm{Coh}(X_{FF})$, and assume that $\mathcal{E}$ only has slopes $< \phi$ and $\mathcal{F}$ only has slopes $\ge \phi$, where $\phi \in \mathbb{R}$. From Proposition \ref{prop:homandext1}, we have $\mathbb{E}\mathrm{xt}^1(\mathcal{E},\mathcal{F}) \simeq 0 \simeq \mathbb{H}\mathrm{om}(\mathcal{F},\mathcal{E}) \simeq 0.$ The above theorem implies that there is a canonical equivalence of Banach--Colmez spaces:
	$$\mathbb{H}\mathrm{om}(\mathcal{E},\mathcal{F}) \simeq \mathbb{E}\mathrm{xt}^{1}(\mathcal{F},\mathcal{E})^{\vee},$$
where we suggestively let $(-)^{\vee}$ denote the functor $\mathcal{H}\mathrm{om}_v(-, \underline{\mathbb{Q}_p}[1])$. Note that there are canonical equivalences $\mathbb{H}\mathrm{om}(\mathcal{E},\mathcal{F}) \simeq \mathbb{H}\mathrm{om}(\mathcal{E},\mathcal{F})^{\vee \vee}$ and $\mathbb{E}\mathrm{xt}(\mathcal{F},\mathcal{E}) \simeq \mathbb{E}\mathrm{xt}(\mathcal{F},\mathcal{E})^{\vee \vee}$. In particular, we also have the equivalence $\mathbb{E}\mathrm{xt}^{1}(\mathcal{F},\mathcal{E}) \simeq \mathbb{H}\mathrm{om}(\mathcal{E},\mathcal{F})^\vee$. Consequently, if $\operatorname{Dim} \mathbb{H}\mathrm{om}(\mathcal{E},\mathcal{F}) = (d, h)$, then $\operatorname{Dim} \mathbb{E}\mathrm{xt}^{1}(\mathcal{F},\mathcal{E}) = (d,-h)$ and vice versa. 

\begin{example} 
Let $\mathcal{E} = \mathcal{O}( \frac{p}{q})$ and $\mathcal{F}  =\mathcal{O}( \frac{r}{s})$, where $p,q,r,s \in \mathbb{Z}$, $\frac{p}{q} < \frac{r}{s}$, $\gcd(p,q)=1$, $\gcd(r,s) = 1$. Then 
	$$R \mathcal{H}\mathrm{om}_{X_{FF}}( \mathcal{E}, \mathcal{F}) \simeq \mathcal{E}^\vee \otimes \mathcal{F}  \simeq \mathcal{O}\left(\frac{qr-ps}{qs}\right)^{\oplus \gcd(qr-ps,qs)}.$$
From Proposition \ref{prop:cohomologies.O(lambda).Fargues--Fontaine}, we obtain that
    \begin{align*}
&\mathbb{H}\mathrm{om}(\mathcal{E},\mathcal{F}) = \mathbb{U}_{(qr-ps,qs)} \quad  \operatorname{Dim} \mathbb{H}\mathrm{om}(\mathcal{E},\mathcal{F})= (qr-ps,qs)\\
&\mathbb{E}\mathrm{xt}^{1}(\mathcal{F},\mathcal{E}) = \mathbb{V}_{(qr-ps, qs)} \quad \operatorname{Dim} \mathbb{E}\mathrm{xt}^{1}(\mathcal{F},\mathcal{E}) = (qr-ps, -qs).
\end{align*}
The above theorem further implies that $\mathbb{U}_{(qr-ps,qs)} \simeq \mathbb{V}_{(qr-ps, qs)}^\vee.$ 
\end{example}

Notably, the formulas above exhibit a remarkable resemblance to \emph{Serre duality} for elliptic curves:

\begin{remark}[Serre duality for elliptic curves]
Recall that if $\pi \colon E_S \to S$ is a family of elliptic curves over a Gorenstein  scheme $S$, then for any $\mathcal{E}, \mathcal{F} \in D^b(\mathrm{Coh}(E_S))$, there is a canonical equivalence in the category of perfect complexes on $S$:
	$$R \pi_* (R \mathcal{H}\mathrm{om}_{E_S}(\mathcal{E},\mathcal{F})) \simeq R \mathcal{H}\mathrm{om}_{\mathcal{O}_S} (R\pi_* R \mathcal{H}\mathrm{om}_{E_S}(\mathcal{F}, \mathcal{E}[1]),  \mathcal{O}_S).$$
If $S = \operatorname{Spec}(\kappa)$ for a field $\kappa$, denote $E = E_{\kappa}$ and let $\mathcal{E}, \mathcal{F} \in \mathrm{Coh}(E)$. The above equivalence then implies that there are canonical equivalences of $\kappa$-vector spaces for all $i \in \{0,1\}$:
	$$\Ext^i_{E}(\mathcal{E}, \mathcal{F}) \simeq \Ext_E^{1-i}(\mathcal{F},\mathcal{E})^*,$$
where $(-)^* = \Hom_\kappa(-, \kappa)$ denotes the $\kappa$-linear dual functor. 
\end{remark}

%\textcolor{blue}{Qingyuan: Up to this point, I have roughly completed the material on Serre duality, except for some modifications.}

\subsection{Stability conditions on Fargues--Fontaine curves}\label{subsection:stability conditions on ff curves}
%\textcolor{red}{Heng: This could be introduced after next subsection or move this subsection after the next one.}
%\textcolor{blue}{Qingyuan: Yes, le bras results can be formulated more naturally using tilted hearts.}
By the classification theorem of vector bundles on the absolute Fargues--Fontaine curve $X_{FF}$, the formula $Z(\mathcal{O}(\frac{p}{q}))=-p+\sqrt{-1}\cdot q\in\mathbb{C}$ determines a group homomorphism $$K_0(\Coh(X_{FF}))\rightarrow \mathbb{C},$$  and $\sigma=(\Coh(X_{FF}), Z)$ is a Bridgeland stability condition on $D^b(X_{FF})$. In fact, this is essentially the only stability condition on $D^b(X_{FF})$.

\begin{prop}\label{Prop: stabilit conditions on FF curves}
    The action of $\widetilde{GL}^+(2,\mathbb{R})$ on $\Stab(X_{FF})$ is free and transitive.
\end{prop}
\begin{proof}
  The argument parallels that of \cite[Theorem 9.1]{bridgeland2007stability}, reflecting the close analog between $X_{FF}$ and elliptic curves. We claim that any indecomposable coherent sheaf $\mathcal{E}$ must be semistable in any stability condition $\sigma\in \Stab(X_{FF})$. This includes every indecomposable torsion sheaf supported at any closed point $x$ in $X_{FF}$.

  Suppose, to the contrary, that $\mathcal{E}$ is not semistable. Then there is a nontrivial triangle $$A\rightarrow \mathcal{E}\rightarrow B$$ with $\Hom(A,B)=0.$ We will show that we have $$\Ext^1(B,A)=0.$$ It suffices to prove this vanishing result for any indecomposable factor of $A$ and $B$. We let $A_1$ and $B_1$ be two indecomposable factors of $A$ and $B$ respectively, and we can assume that there are nontrivial morphisms $f:A_1\rightarrow\mathcal{E}$ and $g:\mathcal{E}\rightarrow B_1$. Otherwise, we could split off the triangle $A_1\rightarrow A_1$ or $B_1\rightarrow B_1$ from the previous distinguished triangle. As $A_1,B_1$ are indecomposable objects, they are semistable with respect to the standard stability condition $\sigma_{st}=(\Coh(X_{FF}), -\mathrm{deg}+\sqrt{-1}\cdot \rank)$. By Proposition \ref{prop:homandext1} and Proposition \ref{prop:cohomologies.O(lambda).Fargues--Fontaine}, we know that $$\phi_{st}(\mathcal{E})-1\leq \phi_{st}(A_1)\leq \phi_{st}(\mathcal{E}),$$ $$\phi_{st}(\mathcal{E})\leq \phi_{st}(B_1)\leq \phi_{st}(\mathcal{E})+1,$$ and $\Hom(A_1,B_1)=0$ implies that $$\phi_{st}(A_1)\geq \phi_{st}(B_1), \text{ or $\phi_{st}(A_1)\leq \phi_{st}(B_1)-1$},$$ and these two equalities hold only when they are torsion sheaves. Here $\phi_{st}$ is the phase of semistable objects with respect to $\sigma_{st}$. 
  
  The case $\phi_{st}(A_1)\geq \phi_{st}(B_1)$ implies that $\phi_{st}(A_1)= \phi_{st}(B_1)=\phi_{st}(\mathcal{E})$ and they are indecomposable torsion sheaves supported at the same closed point $x$ (since there are nontrivial morphisms $f:A_1\rightarrow\mathcal{E}$ and $g:\mathcal{E}\rightarrow B_1$).  This contradicts the assumption $\Hom(A_1,B_1)=0$.

  In the latter case $\phi(A_1)\leq \phi(B_1)-1$, we have $\Ext^1(B_1,A_1)=\Hom(B_1,A_1[1])$. This vanishes unless $\phi(A_1)=\phi(B_1)-1$, which implies that $B_1\simeq A_1[1]$ if $\phi(A_1)\notin \mathbb{Z}$. The isomorphism $B_1\simeq A_1[1]$ contradicts the definition of Harder--Narasimhan filtration of $\mathcal{E}$ with respect to $\sigma$ (see Definition \ref{slicing}.(c)). If $\phi(A_1)\in\mathbb{Z}$, we can assume that $A_1\simeq T_1$ and $B_1\simeq T_2[1]$, where $T_1,T_2$ are two indecomposable torsion sheaves on $X_{FF}$. The condition $\Hom(A_1,B_1)=0$ implies that $T_1,T_2$ are supported on different closed points, which implies that $\Ext^1(B_1,A_1)=\Hom(T_2,T_1)=0.$  Hence $\mathcal{E}\simeq A\oplus B$, contradicting the fact that $\mathcal{E}$ is indecomposable. The claim is proved. 

  Take an element $\sigma=(\mathcal{A},Z)\in \Stab(X_{FF})$. Suppose for a contradiction that the image of $Z$ is contained in a real line in $\mathbb{C}$. It is easy to show that it contradicts the support property in Definition \ref{first WSC}. Indeed, the support property requires a quadratic form $Q$ on $\mathbb{R}^2$, which is negative definite on a one dimensional-subspace in $\mathbb{R}^2$ and nonnegative on all integral vectors in $\mathbb{R}^2$. This is impossible.

  Then we can argue as in \cite[Theorem 9.1]{bridgeland2007stability} to conclude.
\end{proof}

From now on, fix $\sigma=(\Coh(X_{FF}),Z)$ to be the standard stability condition on $D^b(X_{FF})$. This stability condition determines a slicing $\mathcal{P}$ on $D^b(X_{FF})$. We can restate Le Bras's result (see \cite{LeBrasresult}) in terms of this slicing.

Indeed, let $\mathcal{P}[\frac{1}{2}, \frac{3}{2})$ be as in Remark \ref{remark:different hearts}, i.e. objects are in $D^b(X_{FF})$ such that $H^i(C) \neq 0$ if and only if $i=-1,0$ and $H^{-1}(C)$ (resp. $H^0(C)$) has negative (resp. non-negative) slopes. %In other words, objects in $D^b(X_{FF})$ are of the form 
%\[
%[\mathcal{E}_{0} \to \mathcal{E}_1]
%\]
%as a complex of vector bundles over $X_{FF}$ sits in \emph{homological} degree $0,1$ with $\mathcal{E}_{-1}$ (resp. $\mathcal{E}_0$) has negative (resp. non-negative) slope. Note that the condition on $\mathcal{E}_{-1}$ (resp. $\mathcal{E}_{0}$) is equivalent to the condition 
%\[
%\mathrm{H}_0(X_{FF,T},\mathcal{E}_0)=0 (\text{ resp. } %\mathrm{H}_1(X_{FF,T},\mathcal{E}_1)=0)
%\]
%for all sympathetic algebra $T$ over $C^\flat$ by \cite[Prop. 6.8]{LeBrasresult}. We will define a functor $\mathcal{BC}([\mathcal{E}_{0} \to \mathcal{E}_1])$ from the category of perfectoid affinoids over $C^\flat$ to the category of abelian groups as follow. For any perfectoid affinoid $S$ over $C^\flat$, define
%\[
%\mathcal{BC}([\mathcal{E}_{0} \to \mathcal{E}_1])(T)= %\mathbb{H}^0(X_{FF,T}, [\mathcal{E}_{0} \to \mathcal{E}_1])
%\]
%where $\mathbb{H}^\bullet$ meaning taking hypercohomology. It can be deduced from \cite[Prop. II.2.1]{fargues2021geometrization}, that $\mathcal{BC}([\mathcal{E}_{0} \to \mathcal{E}_1])$ defines a sheaf over $\mathrm{AffPerfd}_{C^\flat,v}$, in particular, it defines a pro-\'etale sheave. 

%$\mathcal{BC}([\mathcal{E}_{0} \to \mathcal{E}_1])$ is usually called the relative Banach--Colmez spaces. The name comes from the following main theorem in \cite{LeBrasresult}.

\begin{theorem}
 The functor $R^0\tau_*$ induces an equivalence between $\mathcal{P}[\frac{1}{2}, \frac{3}{2})$ and the category $\mathcal{BC}$. %of abelian sheaves on $\mathrm{AffPerfd}_{C^\flat,v}$ (or $\mathrm{AffPerfd}_{C^\flat,\text{pro\'et}}$). %Moreover, the essential image of $\mathcal{BC}$ is equal to the smallest abelian subcategory ofside the category of abelian sheaves over $\mathrm{AffPerfd}_{C^\flat,v}$ that is stable under extensions and contains the sheaves $\underline{\mathbb{Q}_p}$ and $\mathbb{G}_a$.
\end{theorem}
\begin{proof}
This is a restatement of \cite[Cor. 6.10 and Thm. 7.1]{LeBrasresult}.
\end{proof}

As a corollary of Proposition \ref{proposition of noetherian and aritinian}, we have the following result.

\begin{corollary}
    The category of coherent sheaves on $X_{FF}$ is Noetherian, whereas the category of Banach--Colmez spaces is Artinian.
\end{corollary}

%\textcolor{blue}{Qingyuan: I added the following parts:}

Regarding the hearts $\mathcal{P}(\phi,\phi+1]$ and $\mathcal{P}[\phi, \phi+1)$, we can establish the following results:

\begin{prop} 
\label{prop:tilted.heart.hd=1}
For any $\phi \in \mathbb{R}$, the following statements are true:
\begin{enumerate}
	\item 
	\label{prop:tilted.heart.hd=1-1}
	The inclusions $\mathcal{P}(\phi, \phi+1] \subseteq D^b(X_{FF})$ and $\mathcal{P}[\phi, \phi+1) \subseteq D^b(X_{FF})$ induce canonical equivalences of derived categories:
		\[D^b(\mathcal{P}(\phi, \phi+1]) \simeq D^b(X_{FF}) \quad \text{and} \quad D^b(\mathcal{P}[\phi, \phi+1)) \simeq D^b(X_{FF}).\]
	\item
	\label{prop:tilted.heart.hd=1-2}
	The abelian categories $\mathcal{P}(\phi, \phi+1]$ and $\mathcal{P}[\phi, \phi+1)$ both have homological dimension $\leq 1$.
\end{enumerate}
\end{prop}

\begin{proof}
We will only provide the proof for the abelian category $\mathcal{P}(\phi, \phi+1]$, as the proof for $\mathcal{P}[\phi, \phi+1)$ is identical. 

Claim (1) follows directly from the proof of  \cite[Proposition 2.6]{Modularvariants}. From claim \eqref{prop:tilted.heart.hd=1-1}, for all $A, B \in \mathcal{P}(\phi, \phi+1]$, we have equivalences
\[\mathrm{Ext}_{\mathcal{P}(\phi, \phi+1]}^{i}(A, B ) \simeq \mathrm{Hom}_{D^b(X_{FF})}(A, B[i]) \quad \text{for all $i$}.\]
Therefore, to prove claim \eqref{prop:tilted.heart.hd=1-2}, it suffices to show that
\[\mathrm{Hom}_{D^b(X_{FF})}(A, B[i]) = 0 \quad \text{for all $i \ge 2$ and $A, B \in \mathcal{P}(\phi, \phi+1]$}.\]
We consider the following cases:
\begin{itemize}
	\item If $A, B \in \mathcal{T}$, then $\mathrm{Hom}_{D^b(X_{FF})}(A, B[i]) = \mathrm{Ext}_{X_{FF}}^i(A,B)=0$ for all $i>1$ by Proposition \ref{prop:homandext1}.(3). The same holds for $A, B \in \mathcal{F}[1]$.
	\item If $A \in \mathcal{T}$ and $B \in \mathcal{F}[1]$, then $B[-1] \in \mathcal{F}$, and we have
	\[\mathrm{Hom}_{D^b(X_{FF})}(A, B[i]) = \mathrm{Ext}_{X_{FF}}^{i+1}(A,B[-1]) = 0\]
	for all $i \ge 2$; see Proposition \ref{prop:homandext1}.(3).
	\item If $A \in \mathcal{F}[1]$ and $B \in \mathcal{T}$, then $A[-1] \in \mathcal{F}$. By Serre duality (Theorem \ref{thm:Serre.duality}), we have the following equivalence of Banach--Colmez spaces:
    \[
    \mathbb{E}\mathrm{xt}^{i}(A,B) \simeq \mathbb{E}\mathrm{xt}^{2-i}(B, A[-1])^\vee \simeq 0
    \]
    for all $i \ge 2$, where $(-)^\vee = \mathcal{H}\mathrm{om}_v(-, \underline{\mathbb{Q}_p}[1])$. The second equivalence holds when $i = 2$ due to the properties of slicing, and when $i > 2$ due to Remark \ref{remark:homandext1}.(3). Taking the underlying $C$-points, we obtain $\mathrm{Hom}_{D^b(X_{FF})}(A, B[i]) = 0$ for $i \ge 2$.
\end{itemize}
Since every object in $\mathcal{P}(\phi, \phi+1]$ is an extension of objects in $\mathcal{T}$ and $\mathcal{F}[1]$, claim \eqref{prop:tilted.heart.hd=1-2} follows.
\end{proof}

\begin{remark}
\begin{enumerate}
	\item The same result as in Proposition \ref{prop:tilted.heart.hd=1} holds for any elliptic curve, with the same proof.
	\item An alternative proof of Proposition \ref{prop:tilted.heart.hd=1} can be provided as follows. The calculation in the proof of claim \eqref{prop:tilted.heart.hd=1-2} shows that $\mathrm{Hom}_{D^b(X_{FF})}(A, B[i]) = 0$ for all $A, B \in \mathcal{P}(\phi, \phi+1]$ and $i > 1$ (i.e., $\bigoplus_{i} \mathrm{Hom}_{D^b(X_{FF})}(A, B[i])$ is ``generated in degree $i=1$"). By combining this observation with \cite[Prop. 3.1.6]{beilinson1982faisceaux}, we establish the desired equivalence of \eqref{prop:tilted.heart.hd=1-1}, and then \eqref{prop:tilted.heart.hd=1-2} follows.
\end{enumerate}
\end{remark}

\section{Farey diagrams and continued fractions}\label{Section: Farey diagrams}

\subsection{Farey tessellation} Farey tessellation is a well-studied subject, and it is related to many different branches of mathematics, including number theory, dynamical systems, and Teichm\"uller spaces (see, e.g. \cite{themodularsurfaceandcontinuedfractions}, \cite{Themodulargroup}, \cite{curvatureandrank}, \cite{Circlehomeomorphisms}). In this subsection, we review some basic constructions and results concerning the Farey graph.

First, recall the modular group and the extended modular group, which act on the complex plane $\mathbb{C}$.

\begin{definition}
    The modular group $\Gamma$ is defined as $\mathrm{PSL(2,\mathbb{Z})}$.$$\Gamma=\{z\mapsto \frac{az+b}{cz+d}|a,b,c,d\in\mathbb{Z}, ad-bc=1\}.$$

     The extended modular group $\widetilde{\Gamma}$ is defined as $\mathrm{PGL(2,\mathbb{Z})}.$
     $$\widetilde{\Gamma}=\{z\mapsto \frac{az+b}{cz+d}|a,b,c,d\in\mathbb{Z}, |ad-bc|=1\}.$$
\end{definition}

It is easy to see that $\Gamma$ is generated by $\begin{pmatrix}
1 & 1\\
0 & 1 
\end{pmatrix}$ and $\begin{pmatrix}
0 & 1\\
-1 & 0 
\end{pmatrix}$, while $\widetilde{\Gamma}$ is generated by $\begin{pmatrix}
1 & 1\\
0 & 1 
\end{pmatrix}$ and $\begin{pmatrix}
0 & 1\\
1 & 0 
\end{pmatrix}$.  Moreover, $\Gamma$ is a normal subgroup of $\widetilde{\Gamma}$ of index $2$.

\begin{figure}[ht]
		\centering
		\includegraphics[scale=1]{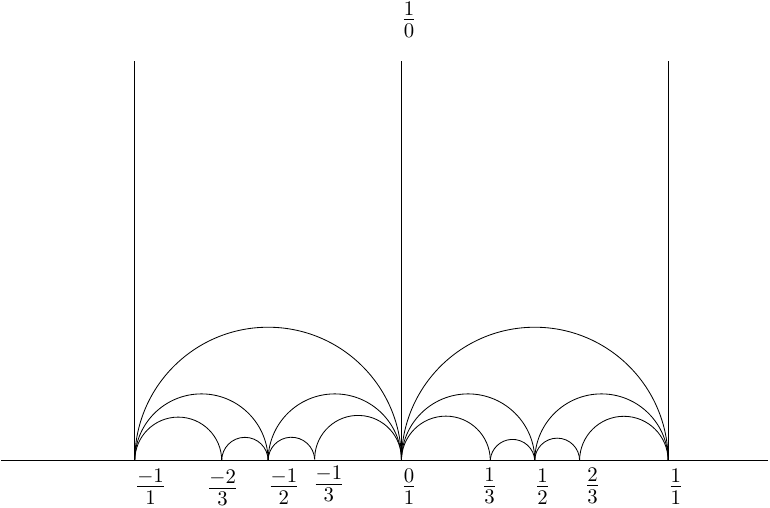}
		\caption{}
        \label{Figure 2}
	\end{figure}

Note that $\Gamma$ and $\widetilde{\Gamma}$ act on the extended complex plane $\mathbb{C}\cup \{\infty\}$ by linear fractional transformations, and $\Gamma$ preserves the upper half-plane $\mathbb{H}^2$, while $\widetilde{\Gamma}$ does not. However, it is a simple fact that the Poincar\'e extension of the action of $\widetilde{\Gamma}$ on $\mathbb{H}^3$, the three-dimensional hyperbolic space, preserves a vertical hyperbolic plane $\mathbb{H}^{\perp}$ in $\mathbb{H}^3$. The resulting collection of fundamental domains gives the \textit{Farey tessellation} (see Figure \ref{Figure 2} ).

\begin{figure}[ht]
		\centering
		\includegraphics[scale=0.8]{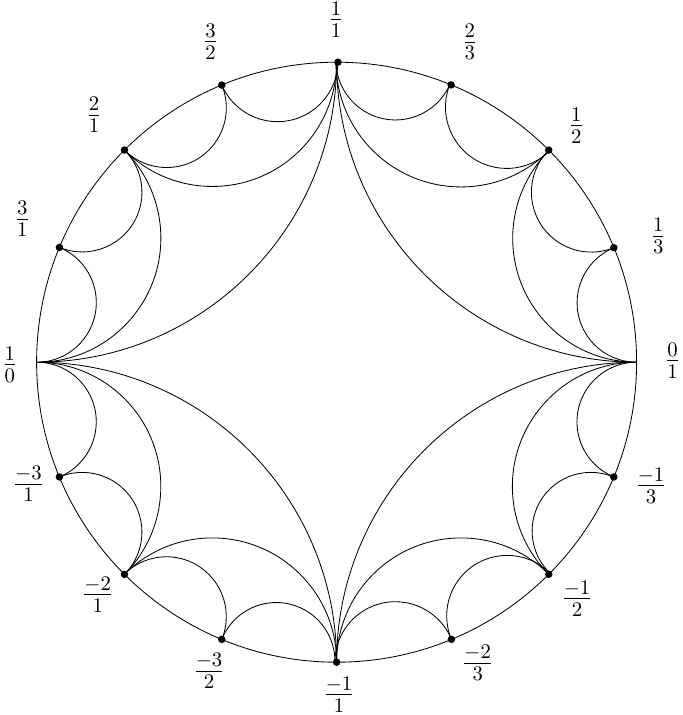}
		\caption{}
        \label{Figure 3}
	\end{figure}

 In the Poincar\'e disc, the Farey tessellation can be drawn as in Figure \ref{Figure 3}. A triangle in this figure is called a \textit{Farey triangle}, and a geodesic in it is called a \textit{Farey geodesic}. We have the following well-known basic results about the Farey triangles.
\begin{theorem}\label{Classcial rsults on Farey graphs}
    The following statements are true.

    \begin{enumerate}
        \item No two Farey geodesics intersect in their interiors.

        \item Two numbers $\frac{p}{q}$, $\frac{r}{s}\in\mathbb{Q}_{\infty}$ are connected by a Farey geodesic if and only if $|ps-rq|=1$.

        \item The extended modular group $\widetilde{\Gamma}$ acts transitively on the set of Farey geodesics.

        \item If $\frac{p}{q}<\frac{x}{y}<\frac{r}{s}\in\mathbb{Q}_{\infty}$ are three vertices of a Farey triangle, then we have $\frac{x}{y}=\frac{r+p}{q+s}$.
    \end{enumerate}
\end{theorem}
\begin{proof}
    These results are well known (see \cite{themodularsurfaceandcontinuedfractions}).
\end{proof}

 % \subsection{Farey diagrams, binary trees and homogeneous trees}

 % In this subsection, we will introduce the notion of Farey graphs, the binary trees, and homogeneous trees associated to Farey graphs.

Given a fixed irrational number $\theta$, for any other number $r\in\mathbb{R}_{\infty}$ different from $\theta$, we have a geodesic connecting these two points $\theta$ and $r$ in the upper half-plane $\mathbb{H}^2$, denoted by $G_{\theta,r}$. We orient this geodesic from $r$ to $\theta$.

The Farey triangles through which $G_{\theta,r}$ passes form a sequence. The sequence is infinite only at one end (if $r$ is rational) or infinite at both ends (if $r$ is irrational). 

\begin{definition}
    Given an irrational number $\theta$,  for any number $r\in\mathbb{R}_{\infty}$ different from $\theta$, the \textit{Farey diagram} associated with $\theta,r$ is defined as a sub-diagram of the Farey tessellation: its vertices are the vertices of Farey triangles through which $G_{\theta,r}$ passes; two such vertices are joined whenever a Farey geodesic connects them. We denote it by $F_{\theta,r}$. 
\end{definition}
\begin{figure}[h]
		\centering
		\includegraphics[scale=0.8]{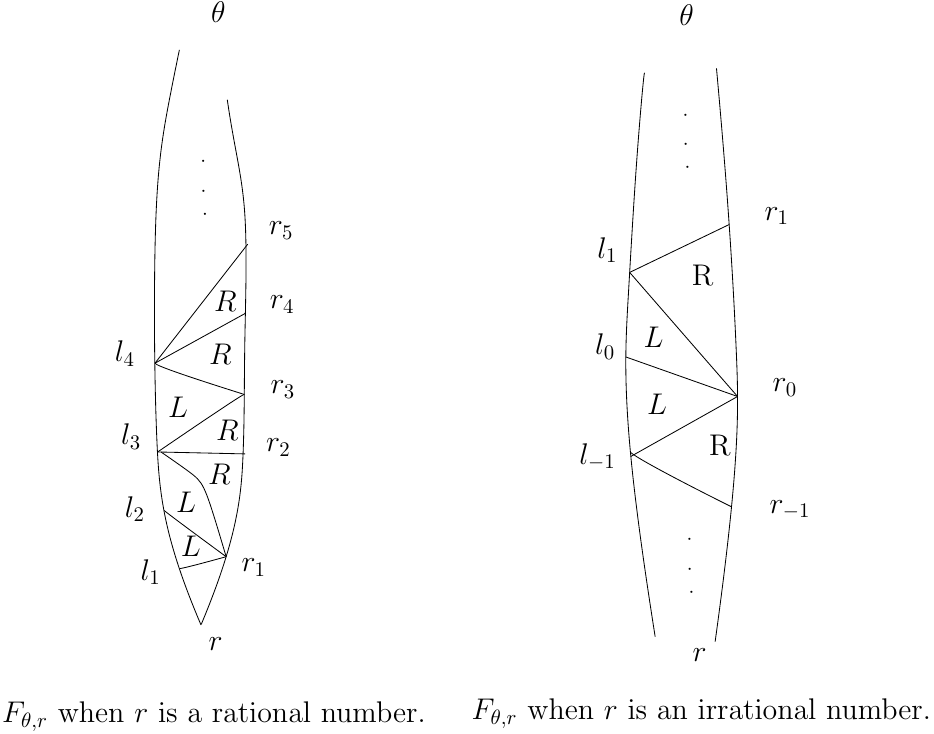}
		\caption{}
        \label{Figure 4}
	\end{figure}
We orient the circle in Figure \ref{Figure 3} anticlockwise.  This orientation provides us with the following definition. 

\begin{definition}\label{Types of triangles}

    A triangle in $F_{\theta,r}$ is called of \textit{type} $L$ if  two vertices of this triangle are in the interior of the arc  from $\theta$ to $r$ in the boundary of the Poincar\'e disc.  A triangle in $F_{\theta,r}$ is called of \textit{type} $R$ if two vertices of this triangle are in the interior of the arc  from $r$ to $\theta$. 

    If $r\in\mathbb{Q}_{\infty}$, we label the vertices of the Farey triangles in $F_{\theta,r}$ as in Figure \ref{Figure 4}, that is, a vertex will be denoted by $l_i$ ($r_j$) if it is on the left arc (right arc), and its subscript $i$ equals the number of vertices on the same arc that are closer to $r$.

    If $r$ is an irrational number and two base points $l_0$ and $r_0$ have been chosen on the left arc and right arc respectively, we label the vertices as indicated in Figure \ref{Figure 4}.
\end{definition} 
 \subsection{Relation with continued fractions}\label{subsection of continued fractions}
In this subsection, we will review some basic facts about continued fractions and their well-known relationship with Farey diagrams (see \cite{themodularsurfaceandcontinuedfractions}).

Let $\theta$ be an irrational number and let $$[a_0;a_1,a_2,\cdots]$$ be the continued fraction representation of $\theta$. We call the rational numbers 

    $$\beta_i\coloneqq [a_0; a_1,\cdots, a_i]=\frac{p_i}{q_i}, \ \text{where $p_i, q_i$ are coprime and $q_i>0$},$$
the convergents of $\theta$. We also have the notion of semi-convergents, i.e. any rational number of the form $$\beta_{i,m}=\frac{p_{i,m}}{q_{i,m}}\coloneqq \frac{p_{i}+mp_{i+1}}{q_{i}+mq_{i+1}}$$ where $m$ is an integer such that $0\leq m\leq a_{i+2}$, is called a semi-convergent of $\theta$. We know that $\beta_{i,0}=\beta_i, \ \beta_{i,a_{i+2}}=\beta_{i+2}$, and $p_iq_{i+1}-p_{i+1}q_i=(-1)^{i+1}$ (see \cite[Chapter 1]{Langbookondiophantineapproximations} for these basic facts). Also note the following common convention: we use $p_{-1}=1$ and $q_{-1}=0$ to denote the starting convergent $\beta_{-1}$ of any continued fractions. This corresponds to the infinite point $\frac{1}{0}$ in the Poincar\'e disc. As we will see in \S \ref{section:limit and colimit objects}, this matches perfectly with other conventions if we associate it with a skyscraper sheaf.

The continued fraction for $\theta\in\mathbb{R}_{>0}$ can be read from the Farey tessellation as follows. Join any point (except the origin) on the imaginary axis to $\theta$ by a hyperbolic geodesic $\gamma$. This arc $\gamma$ cuts a sequence of Farey triangles. Label these triangles by $L$ and $R$ as in Definition \ref{Types of triangles}. Thus we get the so-called \textit{cutting sequence} $L^{n_0}R^{n_1}L^{n_2}\cdots$ of $\theta$, where $n_i\in\mathbb{Z}_{>0}$ for any $i$. 

\begin{prop}\label{continued fractions and Farey diagrams}
  Let $L^{n_0}R^{n_1}L^{n_2}\cdots$ be the cutting sequence of a real number $\theta>1$, then $\theta=[n_0;n_1,n_2,\cdots]$. Likewise, if $R^{n_1}L^{n_2}\cdots$ is the cutting sequence of a real number $0<\theta<1$, then $\theta=[0;n_1,n_2,\cdots]$. 
\end{prop}
\begin{proof}
    See \cite[Theorem A]{themodularsurfaceandcontinuedfractions}.
\end{proof}
\section{Limit and colimit objects $\mathcal{O}(\theta^{\pm})$}\label{section:limit and colimit objects}
\subsection{Minimal triangles and Farey triangles}
In this subsection, we introduce the notion of minimal triangles in $D^b(X_{FF})$, and study its relation with Farey triangles.
\begin{definition}\label{Definition of minimal triangles}
    A \textit{minimal triangle} in the derived category of coherent sheaves on $X_{FF}$ is a triangle $$E\rightarrow F\rightarrow G\rightarrow E[1]$$ satisfying the following conditions.

    \begin{enumerate}

    \item The objects $E,F,G$ are coherent sheaves on $X_{FF}$, and the classes $v(E),v(F),v(G)$ are primitive vectors in $\mathbb{Z}^2$;
        \item for any nontrivial morphism $f: E\rightarrow F$, we have $Cone(f)\simeq G$, or $Cone(f)$ is a torsion sheaf;
        \item  for any nontrivial morphism $g: F\rightarrow G$, we have $Cone(g)\simeq E[1]$;
        \item for any nontrivial morphism $h: G\rightarrow E[1]$, we have $Cone(h)\simeq F[1]$.
 \end{enumerate}
 \end{definition}

\begin{remark}\label{remark about coherence}
    In fact, the requirement $E,F,G\in \Coh(X_{FF})$  is not essential in the definition, and we include this condition only for simplicity.
\end{remark}

In the following theorem,  when we say that a morphism $f$ is an injection or a surjection,  we implicitly mean that $f$ is an injection or a surjection in $\Coh(X_{FF})$.
\begin{theorem}\label{Theorem between minimal triangles and Farey triangles}
    The set of character classes of minimal triangles corresponds bijectively to the set of Farey triangles. 
\end{theorem}
\begin{proof}
   Let us construct the correspondence between these two sets. 

   First, consider a Farey triangle with vertices $s_0,s_1,s_2$. We claim that the following distinguished triangle  $$\mathcal{O}(s_0)\rightarrow\mathcal{O}(s_1)\rightarrow \mathcal{O}(s_2)\rightarrow \mathcal{O}(s_0)[1]$$ is a minimal triangle as in Definition \ref{Definition of minimal triangles}. The proof of this claim is similar to that of \cite[Lemma 3.7]{Modularvariants}, we include it here for the readers' convenience.
   
   Condition (1) is obvious because the rational numbers $s_0,s_1,s_2$ are reduced. We only deal with condition (2), as the remaining conditions can be proved similarly.

We claim that for any nonzero morphism $$f\in \Hom(\mathcal{O}(s_0), \mathcal{O}(s_1)),$$ $f$ is an injection and $\coker(f)\simeq \mathcal{O}(s_2)$, if $s_2\in\mathbb{Q}$; or a torsion sheaf of degree $1$ if $s_2=\frac{1}{0}$. 

We claim that $f$ is an injection. Indeed, if $f$ is not an injection, there is a nontrivial kernel $\ker(f)\neq 0$. Hence we have $$s_0\leq \mu(\operatorname{im}(f))\leq s_1$$ by semistability of $\mathcal{O}(s_0)$, $\mathcal{O}(s_1)$) and $$\rank(\operatorname{im}(f))<\rank (\mathcal{O}(s_0))$$ by nontriviality of $\ker(f)$. However, as we see in Figure \ref{Figure 5}, there is no integral point in the interior of the triangle $A$ by Pick's theorem. This is a contradiction; therefore, $f$ is an injection in $\Coh(X_{FF})$.

\begin{figure}[h]
		\centering
		\includegraphics[scale=0.8]{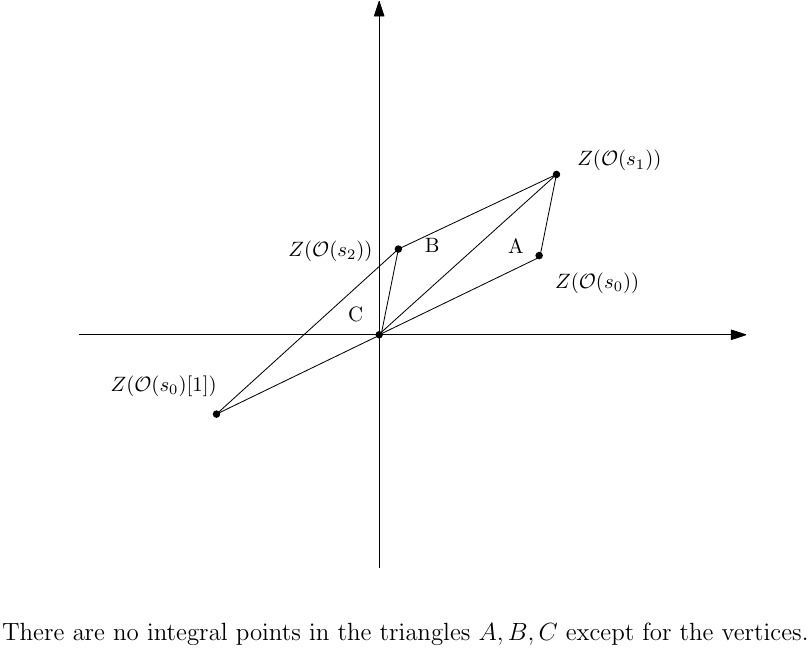}
		\caption{}
        \label{Figure 5}
	\end{figure}
   By a similar argument, we can show that $\coker(f)$ is semistable. Indeed, if $\coker(f)$ is not semistable, consider the last Harder--Narasimhan factor $Q$ of $\coker(f)$. It is easy to see that $s_1\leq \mu(Q)\leq s_2$ and $\rank(Q)<\rank(\mathcal{O}(s_2))$, which contradicts the fact that the triangle $B$ has no integral point in its interior in Figure \ref{Figure 5}. Therefore, we have proved that $\coker(f)$ is semistable, hence $\coker(f)$ is isomorphic to $\mathcal{O}(s_2)$ if $s_2\in\mathbb{Q}$, or a torsion sheaf of degree 1 if $s_2=\frac{1}{0}$.

   For the opposite direction, given a minimal triangle $$E\rightarrow F\rightarrow G\rightarrow E[1],$$ we need to show that there exists a Farey triangle $F$ with vertices $s_0,s_1,s_2$ such that $E\simeq \mathcal{O}(s_0)$, $F\simeq\mathcal{O}(s_1)$ and $G\simeq \mathcal{O}(s_2)$ or a torsion sheaf of degree $1$.

   We first prove that $E,F,G$ are three indecomposable objects, and we treat only $E$, since the other two cases are similar. Suppose, for contradiction, that $E$ decomposes as a direct sum of two objects $E_0 \oplus E_1$, we can further assume that $E_0$ is not isomorphic to $E_1$ by the assumption that $v(E)$ is a primitive vector in $\mathbb{Z}^2$. Hence, a nonzero morphism $f:E_0\rightarrow F$ can be viewed as a morphism $E\rightarrow F$, which is not an injection in $\Coh(X_{FF})$. This contradicts Definition \ref{Definition of minimal triangles}, as $E_1[1]$ is a direct summand of $\coker(f)$, hence $\coker(f)$ can not be isomorphic to an object in $\Coh(X_{FF})$. 

   Hence, we can write $E\simeq \mathcal{O}(s_0)$, $F\simeq\mathcal{O}(s_1)$, $G\simeq \mathcal{O}(s_2)$ for $s_i\in\mathbb{Q}_{\infty}$ and $s_0<s_1<s_2$. Here, by abuse of notation, $\mathcal{O}(\frac{1}{0})$ denotes the skyscraper sheaf on a closed point (which is the support of $G$).  Next, we need to prove that $s_0,s_1,s_2$ are vertices of a Farey triangle. It suffices to prove that $s_0$ and $s_1$ are connected by a Farey geodesic. 

   Assume to the contrary that $s_0$ and $s_1$ are not connected by a Farey geodesic. Let $T$ be the triangle in the complex plane with $0,s_0,s_1$ as its vertices. The Euclidean area of $T$ is at least $1$.  Let $C$ be the convex hull of all integral points in $T$ other than the origin. By Pick's lemma, $C$ contains other integral points beside its vertices. If $C$ has interior points, we let $\lambda_0, \cdots, \lambda_n$ denote the integral points on the edges of $C$ that are below the edge connecting $s_0$ and $s_1$, otherwise,  we let $\lambda_0, \cdots, \lambda_n$ denote the integral points on the edge connecting $s_0$ and $s_1$ (see Figure \ref{Figure 6}), where $n\geq 2$ by assumption. If we identify $\lambda_0, \cdots, \lambda_n$ with their slopes $\frac{p_i}{q_i}$, we order them so that $s_0=\lambda_0<\cdots <\lambda_n=s_1$. It is easy to show that $\lambda_i$ and $\lambda_{i+1}$ are connected by a Farey geodesic by Pick's lemma for any $0\leq i\leq n-1$. 

   \begin{figure}[h]
		\centering
		\includegraphics[scale=0.8]{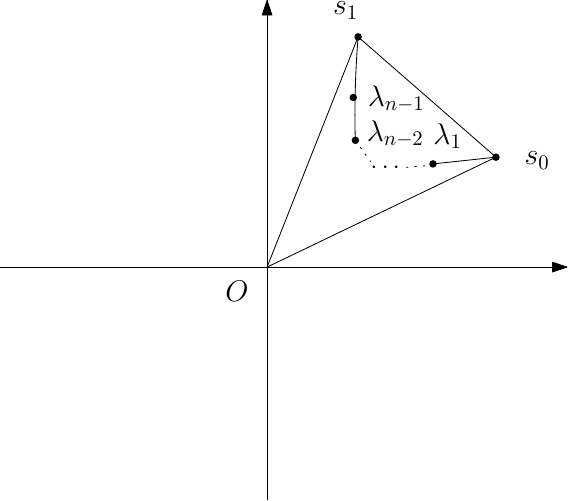}
		\caption{}
        \label{Figure 6}
	\end{figure}

   According to the first half of the proof, any nonzero morphism $f_i\in \Hom(\mathcal{O}(\lambda_i),\mathcal{O}(\lambda_{i+1}))$ is either an injection in $\Coh(X_{FF})$ if $q_i\leq q_{i+1}$, or a surjection in $\Coh(X_{FF})$ if $q_{i}>q_{i+1}$. 

   By the convexity of $C$, one can see that if $f_i$ is an injection, $f_j$ is also an injection for all $n-1\geq j\geq i$. We let 
   
   $$i_0\coloneqq \operatorname{min}\left\{0\leq i\leq n-1 \mid \text{$f_i$ is an injection}\right\},$$ 
   then if $i_0\geq 1$, we can take $$f=f_{n_1}\circ \cdots\circ f_0\in \Hom(\mathcal{O}(s_0),\mathcal{O}(s_1)). $$ As %$$\mathcal{O}(s_0)\stackrel{f_{i_0-1}\circ\cdots\circ f_0}{\twoheadrightarrow} \mathcal{O}(\lambda_{i_0})\xhookrightarrow{f_{n-1}\circ\cdots\circ f_{i_0}}\mathcal{O}(\lambda_n),$$
    \[\begin{tikzcd}[column sep=huge]
    \mathcal{O}(s_0)\arrow[r,twoheadrightarrow,"f_{i_0-1}\circ\cdots\circ f_0"] & \mathcal{O}(\lambda_{i_0})\arrow[r,hookrightarrow,"f_{n-1}\circ\cdots\circ f_{i_0}"] & \mathcal{O}(\lambda_n),
    \end{tikzcd}\]
   we know that $f\neq 0$ and $\operatorname{im}(f)\simeq \mathcal{O}(\lambda_{i_0})$, which contradicts Definition \ref{Definition of minimal triangles}.(2) as $i_0\geq 1$.

   Now we suppose that $i_0=0$, hence all the maps $f_i$ are injections, take $$f=f_{n-1}\circ \cdots\circ f_0\in \Hom(\mathcal{O}(s_0),\mathcal{O}(s_1)). $$ Let $r_i$ be the reduced rational number in $\mathbb{Q}_{\infty}$ such that we have the following minimal triangle $$\mathcal{O}(\lambda_i)\xrightarrow[]{f_i}\mathcal{O}(\lambda_{i+1})\rightarrow\mathcal{O}(r_i)\rightarrow\mathcal{O}(\lambda_i)[1].$$
By convexity of $C$, we know that $r_0\geq r_1\geq\cdots\geq r_{n-1}$. Let 

    $$j_0\coloneqq \operatorname{min}\{1\leq j\leq n-1 \mid \text{$r_{j}$ is not $\infty$} \}.$$ 
As $v(\mathcal{O}(s_2))$ is primitive, we know that not all $r_j$ are $\frac{1}{0}=\infty$, hence $j_0$ is well defined. By Proposition \ref{prop:homandext1} and convexity of $C$, it is straightforward to check that $$\coker(f)\simeq \bigoplus_{j=j_0}^{n-1}\mathcal{O}(r_j)\oplus \coker(f_{j_0-1}\circ \cdots \circ f_0).$$
  
 This contradicts the indecomposability of $\coker(f)$ as $n\geq 2$. Hence, $s_0$ and $s_1$ are connected by a Farey geodesic, and furthermore, a minimal triangle corresponds to a Farey triangle.  
   One can easily check that these two maps are inverse to each other.
\end{proof}

\begin{remark}
    This theorem provides the converse direction of \cite[Theorem 1.2]{Modularvariants}.
\end{remark}
\begin{remark}
    As we mentioned in Remark \ref{remark about coherence}, the condition $E,F,G\in \Coh(X_{FF})$ is not essential. If we drop this condition, the only ambiguity in Theorem \ref{Theorem between minimal triangles and Farey triangles} arises from the equality $\frac{p}{q}=\frac{-p}{-q}$, where $\frac{p}{q}$ corresponds to $\mathcal{O}(\frac{p}{q})$, while $\frac{-p}{-q}$ corresponds to $\mathcal{O}(\frac{p}{q})[1]$ if $q>0$.

    It is a useful terminology coincidence that the minimal triangle corresponding to the fundamental Farey triangle is the fundamental short exact sequence in Banach--Colmez spaces (see Figure \ref{Figure 7}).

      \begin{figure}[h]
		\centering
		\includegraphics[scale=0.66]{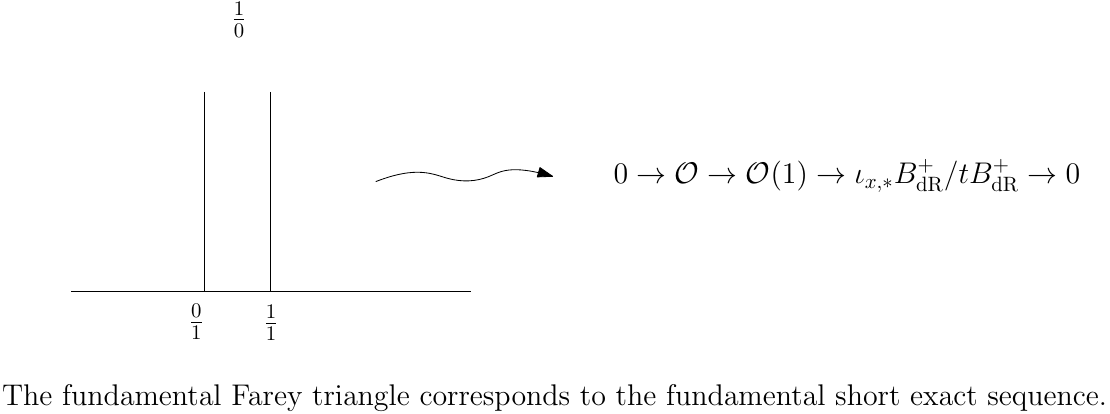}
		\caption{}
        \label{Figure 7}
	\end{figure}

\end{remark}

By the previous theorem, we can further study the abelian categories $\mathcal{P}(\phi,\phi+1]$ and $\mathcal{P}[\phi,\phi+1)$. 
\begin{corollary}\label{Hearts of FF curves}
	Let $\sigma=(\mathcal{A},Z)$ be the standard Bridgeland stability condition on $D^b(X_{FF})$ as in \S \ref{subsection:stability conditions on ff curves},   we have the following results. 
	\begin{enumerate}
		\item If $\cot(\pi \phi)$ is in $\mathbb{Q}_{\infty}$ , the heart $\mathcal{P}(\phi, \phi+1]$ is Noetherian and non-Artinian, and the heart $\mathcal{P}[\phi, \phi+1)$ is Artinian and non-Noetherian.
		\item If $\cot(\pi \phi)$ is an irrational number, the heart $\mathcal{P}(\phi, \phi+1]=\mathcal{P}[\phi, \phi+1)$ is neither Noetherian nor Artinian.
	\end{enumerate}

\end{corollary}

\begin{proof}
    The proof is essentially the same as the proof of the main theorem in \cite{Differentheartsonellipticcurves}, we include the proof here for the reader's convenience.

    Half of part (1) follows directly from Proposition \ref{proposition of noetherian and aritinian}. Indeed, as in Remark \ref{rotation}, we can rotate the stability condition $\sigma$ by $e^{-\sqrt{-1}\phi}$. If $\cot(\pi \phi)$ is a rational number or undefined, the stability condition $e^{-\sqrt{-1}\phi}\sigma$ satisfies the assumption in Proposition \ref{proposition of noetherian and aritinian}, hence $(e^{-\sqrt{-1}\phi}\mathcal{P})(0,1]=\mathcal{P}(\phi,\phi+1]$ is Noetherian, and $(e^{-\sqrt{-1}\phi}\mathcal{P})[0,1)=\mathcal{P}[\phi,\phi+1)$ is Artinian. 

For the other half of part (1), there exist two rational numbers $\frac{p_0}{q_0}>\frac{p_{-1}}{q_{-1}}$ that are larger than $-\cot{\pi\phi}$ and the numbers $\frac{p_0}{q_0},\frac{p_{-1}}{q_{-1}}, -\cot{\pi\phi}=\frac{r}{s}$ are the vertices of a Farey triangle. If we denote $\frac{p_{-i}}{q_{-i}}\coloneqq \frac{p_0-ir}{q_0-is}$ for any $i\in\mathbb{Z}$. By Theorem \ref{Theorem between minimal triangles and Farey triangles}, we have the following sequence $$\cdots \hookrightarrow \mathcal{O}(\frac{p_{-i}}{q_{-i}})\hookrightarrow \mathcal{O}(\frac{p_{-1}}{q_{-1}})\hookrightarrow\mathcal{O}(\frac{p_0}{q_0})$$ in $\mathcal{P}(\phi,\phi+1]$. Hence $\mathcal{P}(\phi,\phi+1]$ is not Artinian; similarly, the sequence $$\mathcal{O}(\frac{p_0}{q_0})\twoheadrightarrow \mathcal{O}(\frac{p_1}{q_1})\twoheadrightarrow\cdots$$ implies that $\mathcal{P}[\phi,\phi+1)$ is not Noetherian.

    For part (2), without loss of generality, we may assume that $\phi\in (0, 1)$. Let $$[a_0;a_1,a_2,\cdots]$$ be the continued fraction representation of $\theta\coloneqq -\cot(\pi\phi)$, let $$\beta_i\coloneqq [a_0; a_1,\cdots, a_i]=\frac{p_i}{q_i}, \ q_i>0$$ and  $$\beta_{i,m}=\frac{p_{i,m}}{q_{i,m}}\coloneqq \frac{p_{i}+mp_{i+1}}{q_{i}+mq_{i+1}}$$ be the associated convergents and semi-convergents as in \S \ref{subsection of continued fractions}, where $m$ is an integer such that $0\leq m\leq a_{i+2}$.

    By the results in \S \ref{subsection of continued fractions}, we know that $$\beta_{i,m}, \beta_{i,m+1}, \beta_{i+1}$$ form a Farey triangle for any integer $i\geq -1$ and integer $m$ with $0\leq m\leq a_{i+2}-1$ (see Figure \ref{Figure 8} for a schematic graph. We call the associated graph a \textit{roller coaster}. It has no integral points in its interior and, as a labelled graph, agrees with $F_{\theta, \frac{1}{0}}$).

 \begin{figure}[h]
		\centering
		\includegraphics[scale=0.8]{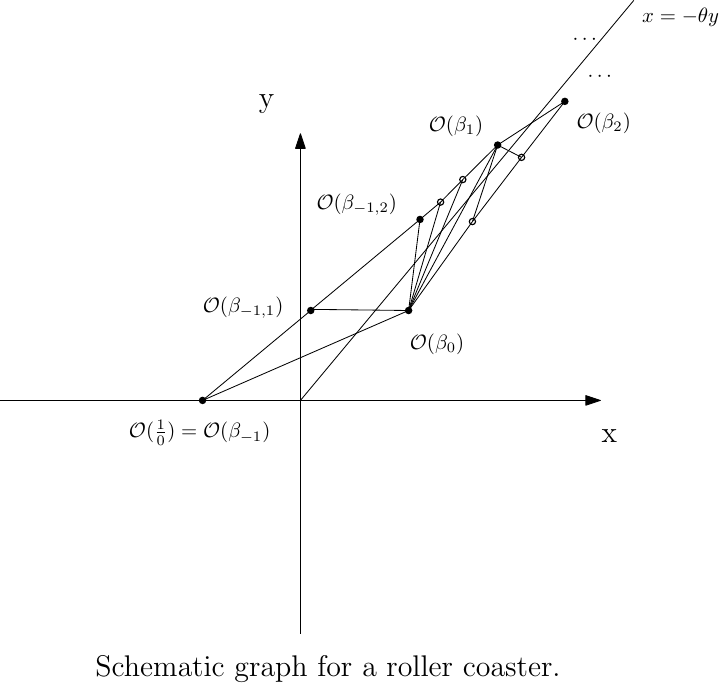}
		\caption{}
        \label{Figure 8}
	\end{figure}

    Hence, by Theorem \ref{Theorem between minimal triangles and Farey triangles}, we know that any nontrivial morphism $$f\in \Hom(\mathcal{O}(\beta_{2i,m}),\mathcal{O}(\beta_{2i,m+1}))$$ is an injection in $\mathcal{P}(0,1]$ with $\coker(f)\simeq \mathcal{O}(\beta_{2i+1})$. Similarly, any nontrivial morphism $$g\in \Hom(\mathcal{O}(\beta_{2i+1,m+1}),\mathcal{O}(\beta_{2i+1,m}))$$ is a surjection in $\mathcal{P}(0,1]$ with $\ker(g)\simeq\mathcal{O}(\beta_{2i+2})$. 

    Therefore, any nontrivial morphism $$f[1]\in \Hom(\mathcal{O}(\beta_{2i,m})[1],\mathcal{O}(\beta_{2i,m+1})[1])$$ is a surjection in $\mathcal{P}(\phi,\phi+1]$ with $\ker(f[1])\simeq \mathcal{O}(\beta_{2i+1})$. Similarly, any nontrivial morphism $$g\in \Hom(\mathcal{O}(\beta_{2i+1,m+1}),\mathcal{O}(\beta_{2i+1,m}))$$ is an injection in $\mathcal{P}(\phi,\phi+1]$ with $\coker(g)\simeq\mathcal{O}(\beta_{2i+2})[1]$. 

    We get the following two infinite sequences in $\mathcal{P}(\phi,\phi+1]$

    $$ \mathcal{O}(\beta_{0})[1]\twoheadrightarrow \mathcal{O}(\beta_{2})[1]\twoheadrightarrow\cdots \twoheadrightarrow \mathcal{O}(\beta_{2i})[1] \twoheadrightarrow\cdots ,$$ and $$\cdots\hookrightarrow\mathcal{O}(\beta_{2i+1})\hookrightarrow\cdots\hookrightarrow \mathcal{O}(\beta_{1})\hookrightarrow \mathcal{O}(\beta_{-1}).$$

    These two infinite sequences imply that $\mathcal{P}(\phi,\phi+1]=\mathcal{P}[\phi,\phi+1)$ is neither Noetherian nor Artinian.
\end{proof}

\begin{remark}
    Similar results hold for elliptic curves (see \cite{Differentheartsonellipticcurves}). 
\end{remark}

\subsection{Continuum envelopes} As we saw in the proof of the previous corollary, we have two infinite sequences $$  \mathcal{O}(\beta_{0})\hookrightarrow \mathcal{O}(\beta_{2})\hookrightarrow\cdots \hookrightarrow \mathcal{O}(\beta_{2i})\hookrightarrow \cdots$$ and $$\cdots\twoheadrightarrow\mathcal{O}(\beta_{2i+1})\twoheadrightarrow \cdots\twoheadrightarrow \mathcal{O}(\beta_{1})\twoheadrightarrow \mathcal{O}(\beta_{-1})$$ in $\mathcal{P}(0,1]=\Coh(X_{FF})$. These two infinite sequences lead us to consider the associated limit and colimit objects in $\QCoh(X_{FF})$, which is a complete and cocomplete abelian category. 

\begin{definition}\label{definition of limit and colimit objects}
    Let $\theta$ be an irrational number, $\beta_i=\frac{p_i}{q_i},\beta_{i,m}=\frac{p_{i,m}}{q_{i,m}}$ be its associated convergents and semi-convergents  of the continued fraction $$\theta= [a_0;a_1,a_2,\cdots],$$ where $i\geq -1$ is an integer and $0\leq m\leq a_{i+2}$. We use $$\mathcal{O}(\theta^{+}, \{f_{2i+1}\})\coloneqq \varprojlim_{f_{2i+1}}\mathcal{O}(\beta_{2i+1})$$ to denote the limit of 
   \[ \begin{tikzcd}
        \cdots \arrow[r,twoheadrightarrow,"f_{2i+1}"] & \mathcal{O}(\beta_{2i-1}) \arrow[r,twoheadrightarrow,"f_{2i-1}"] & \mathcal{O}(\beta_{2i-3}) \arrow[r,twoheadrightarrow,"f_{2i-3}"] & \cdots\arrow[r,twoheadrightarrow,"f_{1}"] & \mathcal{O}(\beta_{-1}), 
    \end{tikzcd}\]
where $\{f_{2i+1}\}$ is a sequence of nontrivial morphisms, and use $$\mathcal{O}(\theta^{-},\{g_{2i}\})\coloneqq \varinjlim_{g_{2i}}\mathcal{O}(\beta_{2i})$$ 
  to denote the colimit of 
    \[\begin{tikzcd}
     \mathcal{O}(\beta_0) \arrow[r, hookrightarrow,  "g_{0}"] &    \cdots \arrow[r, hookrightarrow,  "g_{2i-4}"] & \mathcal{O}(\beta_{2i-2}) \arrow[r, hookrightarrow,  "g_{2i-2}"] & \mathcal{O}(\beta_{2i}) \arrow[r, hookrightarrow,  "g_{2i}"] & \cdots, 
    \end{tikzcd}\]
    where $\{g_{2i}\}$ is a sequence of nontrivial morphisms. 

    We use $\QCoh_{\mathbb{R}}(X_{FF})$ to denote the smallest full abelian  subcategory of $\QCoh(X_{FF})$ that contains $\Coh(X_{FF})$ and  $\mathcal{O}(\theta^{+}, \{f_{2i+1}\})$, $\mathcal{O}(\theta^{-},\{g_{2i}\})$ for all irrational numbers $\theta $ and sequences of nontrivial morphisms $\{f_{2i+1}\}$ and $\{g_{2i}\}$. 
    
    We call $\QCoh_{\mathbb{R}}(X_{FF})$ the \textit{continuum envelope} of $\Coh(X_{FF})$.
\end{definition}

\begin{definition}
    
    For an elliptic curve $E$, define $\QCoh_{\mathbb{R}}(E)$ similarly: it is the minimal full abelian subcategory of $\QCoh(E)$ that contains $\Coh(E)$ and $$\mathcal{O}(\theta^{+}, \{\mathcal{L}_{2i+1}\}, \{f_{2i+1}\}),\ \mathcal{O}(\theta^{-},\{\mathcal{L}_{2i}\}, \{g_{2i}\})$$ where $\mathcal{L}_{2i+1}$ (respectively, $\mathcal{L}_{2i}$) is a vector bundle on $E$ with rank $q_{2i+1}$ ($q_{2i}$) and degree $p_{2i+1}$ ($p_{2i}$), and  $$\mathcal{O}(\theta^{+}, \{\mathcal{L}_{2i+1}\}, \{f_{2i+1}\})\coloneqq \varprojlim_{f_{2i+1}}\mathcal{L}_{2i+1}$$ is the limit object of  \[ \begin{tikzcd}
        \cdots \arrow[r,twoheadrightarrow,"f_{2i+1}"] & \mathcal{L}_{2i-1} \arrow[r,twoheadrightarrow,"f_{2i-1}"] & \mathcal{L}_{2i-3} \arrow[r,twoheadrightarrow,"f_{2i-3}"] & \cdots\arrow[r,twoheadrightarrow,"f_{1}"] & \mathcal{L}_{-1}, 
    \end{tikzcd}\]
and $$\mathcal{O}(\theta^{-},\{\mathcal{L}_{2i}\}, \{g_{2i}\})\coloneqq \varinjlim_{g_{2i}}\mathcal{L}_{2i}$$ is the colimit object of   \[\begin{tikzcd}
     \mathcal{L}_0 \arrow[r, hookrightarrow,  "g_{0}"] &    \cdots \arrow[r, hookrightarrow,  "g_{2i-2}"] & \mathcal{L}_{2i} \arrow[r, hookrightarrow,  "g_{2i}"] & \mathcal{L}_{2i+2} \arrow[r, hookrightarrow,  "g_{2i+2}"] & \cdots. 
    \end{tikzcd}\]  Here $\{f_{2i+1}\}$ and $\{g_{2i}\}$ are two sequences of nontrivial morphisms.

     We call $\QCoh_{\mathbb{R}}(E)$ the \textit{continuum envelope} of $\Coh(E)$.
  \end{definition}

  \begin{remark}
 As we will see in \S \ref{section: homological algebra}, the limit and colimits objects in $\QCoh_{\mathbb{R}}(X_{FF})$ depend on  the choices of morphisms $\{f_{2i+1}\}, \{g_{2i}\}$, and form a very large moduli space. The limit and colimits objects in $\QCoh_{\mathbb{R}}(E)$ depend on  the choices of objects $\mathcal{L}_{i}$ and the choices of morphisms $\{f_{2i+1}\}, \{g_{2i}\}$, and form an infinite-dimensional moduli space.

  The main difference between $X_{FF}$ and $E$ is that the corresponding objects $\mathcal{L}_{i}$ on $E$ have a moduli space isomorphic to $\mathrm{Pic}^0(E)$. On the other hand, for $X_{FF}$, the objects $\mathcal{O}(\beta_{i})$ are rigid, but the parameter spaces of sequences of nontrivial morphisms $\{f_{2i+1}\}$ and $\{g_{2i}\}$ are very large.

  Note that we only consider convergents rather than all semi-convergents in the definitions of $\QCoh_{\mathbb{R}}(X_{FF})$ and $\QCoh_{\mathbb{R}}(E)$. The next proposition implies that we would get the same continuum envelope for $X_{FF}$. 
\end{remark}

From now on, we will state the results mainly for Fargues--Fontaine curves, and comment on the corresponding results for elliptic curves. We also fix the following setup.

\textbf{Set-up}:  Let $\theta$, $\beta_i=\frac{p_i}{q_i},\beta_{i,m}=\frac{p_{i,m}}{q_{i,m}}$ be as in Definition \ref{definition of limit and colimit objects}, $\mathcal{O}(\beta_{i})$ and $\mathcal{O}(\beta_{i,m})$ be the associated vector bundles on $X_{FF}$. For each $i$, choose nonzero morphisms $g_{2i}$ in $\Hom(\mathcal{O}(\beta_{2i}), \mathcal{O}(\beta_{2i+2}))$, nonzero morphisms $f_{2i+1}$ in $\Hom(\mathcal{O}(\beta_{2i+1}), \mathcal{O}(\beta_{2i-1}))$, and $$f_{2i+1,n}\in \Hom(\mathcal{O}(\beta_{2i+1,n}),\mathcal{O}(\beta_{2i+1,n-1})),$$ $$g_{2i,m}\in \Hom(\mathcal{O}(\beta_{2i,m}),\mathcal{O}(\beta_{2i,m+1}))$$ be nonzero morphisms.

\begin{prop}\label{prop:focus on convergents}
 With the notation of Definition \ref{definition of limit and colimit objects}, the following map is surjective: 
  \[\begin{tikzcd}
      \mathrm{Pr}(\otimes_{m=0}^{a_{2i+2}-1} \Hom(\mathcal{O}(\beta_{2i,m}),\mathcal{O}(\beta_{2i,m+1})))  \rightarrow \Hom(\mathcal{O}(\beta_{2i}),\mathcal{O}(\beta_{2i+2})).
  \end{tikzcd} \] 
  
 Here, $\mathrm{Pr}$ denotes the subset of pure tensors in the tensor product. A similar result holds for the other parity of semi-convergents. 
\end{prop}
\begin{proof}
    In fact, we will prove a more general result: the natural composition map 
    
   \[ \begin{tikzcd}
      \mathrm{Pr}(\otimes_{m=n_0}^{n_1-1} \Hom(\mathcal{O}(\beta_{2i,m}),\mathcal{O}(\beta_{2i,m+1})))  \rightarrow \Hom(\mathcal{O}(\beta_{2i,n_0}),\mathcal{O}(\beta_{2i,n_1}))
  \end{tikzcd} \]
  is surjective for any $0\leq n_0<n_1\leq a_{2i+2}$. We prove this by induction on $n_1-n_0$. The case $n_1-n_0=1$ is immediate. 
  
  For the induction step, assume the result holds when  $n_1-n_0=k$, and now let $n_1-n_0=k+1$. We claim that for any nontrivial morphism $f\in \Hom(\mathcal{O}(\beta_{2i,n_0}),\mathcal{O}(\beta_{2i,n_1}))$, we have the following short exact sequence in $\Coh(X_{FF})$. $$0\rightarrow \mathcal{O}(\beta_{2i,n_0})\xrightarrow{f} \mathcal{O}(\beta_{2i,n_1}) \xrightarrow{h} \mathcal{O}(\beta_{2i+1})^{k+1}\rightarrow 0$$
The injectivity of $f$ follows in exactly the same way as in the proof of Theorem \ref{Theorem between minimal triangles and Farey triangles}. 

It remains to prove that $\coker(f)\simeq \mathcal{O}(\beta_{2i+1})^{k+1}$. Suppose otherwise, and let $Q$ denote the last Harder--Narasimhan factor of $\coker(f)$, we know that $$\beta_{2i+1}>\mu(Q)> \beta_{2i,n_1}\quad \text{and} \quad \rank (Q)<(k+1) \ \rank(\mathcal{O}(\beta_{2i+1})).$$ 
However, there is no such integral point by Pick's theorem (see Figure \ref{Figure 9} for a pictorial illustration when $k+1=3$). Hence we have proved our claim. 
\begin{figure}[h]
		\centering
		\includegraphics[scale=0.8]{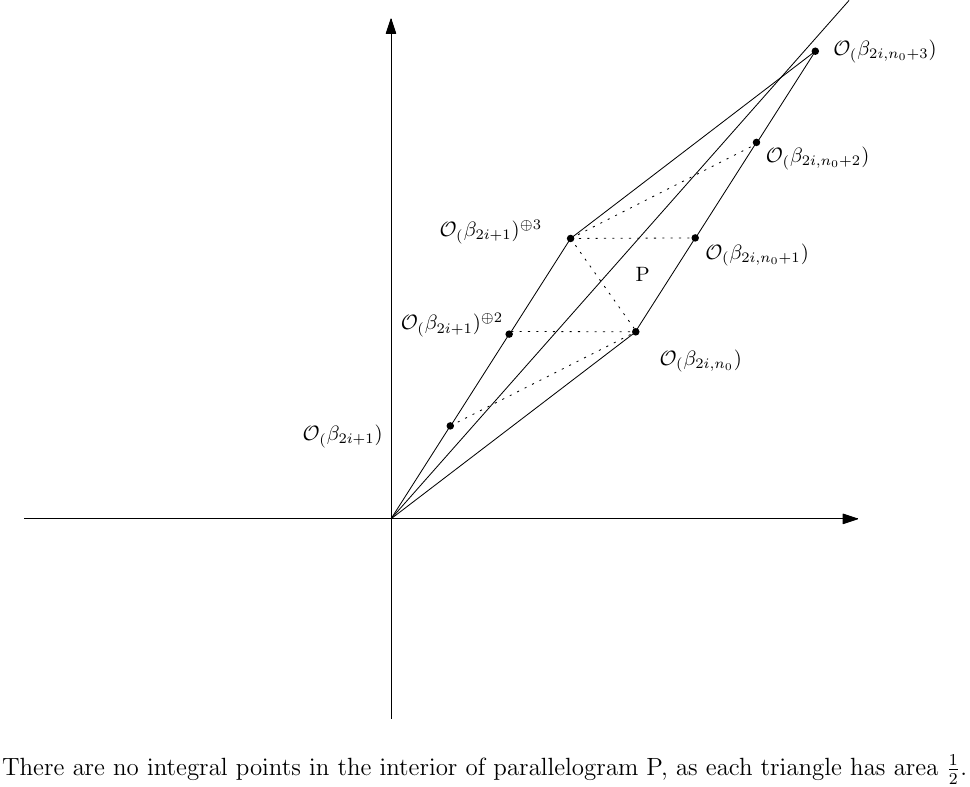}
		\caption{}
        \label{Figure 9}
	\end{figure}

Consider the composition of $$p\circ h: \mathcal{O}(\beta_{2i,n_1}) \xrightarrow{h}\mathcal{O}(\beta_{2i+1})^{k+1}\xrightarrow{p} \mathcal{O}(\beta_{2i+1}),$$ where $p$ is the projection of $\mathcal{O}(\beta_{2i+1})^{k+1}$ to its first factor. As $$\beta_{2i+1}, \beta_{2i,n_1},\beta_{2i,n_1-1}$$ are the vertices of a Farey triangle, by Theorem \ref{Theorem between minimal triangles and Farey triangles}, we get the following commutative diagram \[\begin{tikzcd}
    0\arrow[r] & \mathcal{O}(\beta_{2i,n_0}) \arrow[r,"f"] \arrow[d,"g"]& \mathcal{O}(\beta_{2i,n_1}) \arrow[r,"h"] \arrow[d,"id"]& \mathcal{O}(\beta_{2i+1})^{k+1} \arrow[r] \arrow[d,"p"] & 0 \\  0\arrow[r] & \mathcal{O}(\beta_{2i,n_1-1}) \arrow[r,"f'"] & \mathcal{O}(\beta_{2i,n_1}) \arrow[r,"p\circ h"] & \mathcal{O}(\beta_{2i+1}) \arrow[r] & 0
\end{tikzcd}\]

Hence, we proved that the composition map \[\begin{tikzcd}
   \mathrm{Pr}(\Hom(\mathcal{O}(\beta_{2i,n_1-1}),\mathcal{O}(\beta_{2i,n_1}))\otimes \Hom( \mathcal{O}(\beta_{2i,n_0}), \mathcal{O}(\beta_{2i,n_1-1}))) \\ \longrightarrow \Hom(\mathcal{O}(\beta_{2i,n_0}), \mathcal{O}(\beta_{2i,n_1})) 
\end{tikzcd}\] is surjective. Combining this with the inductive assumption, we complete the proof.
 \end{proof}
\begin{remark}\label{convergents and semiconvergents}
The result \cite[Theorem 5.6.29]{farguesfontaine-courbes} can be viewed as a special case of this proposition.
    This proposition implies that regardless of whether we consider convergents only or all semi-convergents, we get the same class of colimit and limit objects in $\QCoh(X_{FF})$. 
    
  Although this proposition does not hold verbatim for elliptic curves,  there is an analogous result in that setting. The main difference is that there exist indecomposable vector bundles on elliptic curves whose character class is not primitive (see \cite{Atiyahclassification}).

\end{remark}

\begin{prop}\label{prop: focus on convergents ellpitic curve case}
    Let $E$ be an elliptic curve.  For any two vector bundles $\mathcal{L}_{2i}, \mathcal{L}_{2i+2}$ on $E$ with Chern characters $$Ch(\mathcal{L}_{2i})=(p_{2i},q_{2i}), Ch(\mathcal{L}_{2i+2})=(p_{2i+2}, q_{2i+2}),$$ and any nonzero morphism $$f\in \Hom(\mathcal{L}_{2i}, \mathcal{L}_{2i+2}),$$ there exists a vector bundle $\mathcal{L}_{2i,m}$ with Chern character equal to $(p_{2i,m},q_{2i,m})$ for each $0\leq m\leq a_{2i+2}$, and $f$ factors as $$\mathcal{L}_{2i}\rightarrow \mathcal{L}_{2i,1}\rightarrow \cdots \rightarrow \mathcal{L}_{2i, a_{2i+2}-1}\rightarrow \mathcal{L}_{2i+2}.$$ 
\end{prop}

\begin{proof}
    The proof is similar to the proof of the previous proposition. Indeed, we can prove that $f$ is an injection in $\Coh(E)$ and its cokernel is a semistable object of slope $\frac{p_{2i+2}-p_{2i}}{q_{2i+2}-q_{2i}}$ by the same argument in the proof of Proposition \ref{prop:focus on convergents}. By \cite[Lemma 2]{FMtransformsandvectorbundlesonellipticcurve}, we know that if $a_{2i+2}>1$, we have the following commutative diagram. \[\begin{tikzcd}
        0\arrow[r] & \mathcal{L}_{2i}\arrow[r] \arrow[d,"id"] & \mathcal{L}_{2i,m} \arrow[d] \arrow[r] & V\arrow[d] \arrow[r] & 0 \\ 0\arrow[r] & \mathcal{L}_{2i}\arrow[r,"f"] & \mathcal{L}_{2i+2}\arrow[r] & \coker(f)\arrow[r] & 0
    \end{tikzcd}\] where $V$ is a semistable  vector subbundle of $\coker(f)$, and  $Ch(V)=(dp_{2i+1}, dq_{2i+1})$, where $d<a_{2i+2}$. And $\mathcal{L}_{2i,m}$ is the pullback of $V$ with respect to the surjection $\mathcal{L}_{2i+2}\rightarrow \coker(f)$. The proposition follows by induction.
\end{proof}

Hence, we only need to consider convergents in the definition of $\QCoh_{\mathbb{R}}(E)$, the continuum envelope of coherent sheaves on an elliptic curve $E$.

\subsection{Homological algebra on the roller coaster $F_{\theta,\frac{1}{0}}$} We end this section by studying some homological algebra on the roller coaster in Figure \ref{Figure 8}, to simplify the statements, we need the following definitions.

\begin{definition}\label{definition: vector bundle associated with the shortest path}
    We orient the Farey graph $F_{\theta,\frac{1}{0}}$ (the roller coaster in Figure \ref{Figure 8}) as follows: an edge $e$ joining $\lambda<\mu$ is oriented from $\lambda$ to $\mu$. To each such edge $e$, we associate it with a vector bundle $\mathcal{O}(e)$ on $X_{FF}$, which corresponds to the third vertex of the Farey triangle in $F_{\theta,\frac{1}{0}}$ that contains $e$.

    Hence, for any two different vertices $\lambda_1<\lambda_2$ in the directed graph $F_{\theta,\frac{1}{0}}$, there is a unique shortest directed path $P_{\lambda_1,\lambda_2}$ from $\lambda_1$ to $\lambda_2$.  We associate to the path $P_{\lambda_1,\lambda_2}$ with the vector bundle $$\mathcal{O}(P_{\lambda_1,\lambda_2})\coloneqq \oplus_{e\in P_{\lambda_1,\lambda_2}} \mathcal{O}(e), $$ where $e\in P_{\lambda_1,\lambda_2}$ means that $e$ is an edge in the path $P_{\lambda_1,\lambda_2}$.
\end{definition}

\begin{definition}
     Let $\theta$, $\beta_i=\frac{p_i}{q_i},\beta_{i,m}=\frac{p_{i,m}}{q_{i,m}}$ be as in Definition \ref{definition of limit and colimit objects}, and let $F_{\theta,\frac{1}{0}}$ be the directed Farey graph, which we also call a roller coaster. We call an edge $e\in F_{\theta,\frac{1}{0}}$ an \textit{exterior edge} if $e$ is on the boundary of $F_{\theta,\frac{1}{0}}$, i.e., $e$ either connects $\beta_0$ and $\beta_{-1}$ or $e$ is any edge that lies on one side of the line $x=-\theta y$ in Figure \ref{Figure 8}.
     
     Similarly, we call an edge $e\in F_{\theta,\frac{1}{0}}$ an \textit{interior edge} if $e$ lies in the interior of $F_{\theta,\frac{1}{0}}$, i.e., $e$ either connects $\beta_{2i,m}$ with $\beta_{2i+1}$ for some $i\geq 0$ and $0< m\leq a_{2i+2}$, or connects $\beta_{2j}$ with $\beta_{2j-1,n}$ for some $j\geq 0$ and $0<n\leq a_{2j+1}$.
\end{definition}
\begin{theorem}\label{theorem about the roller coaster}
     Let $\theta$, $\beta_i=\frac{p_i}{q_i},\beta_{i,m}=\frac{p_{i,m}}{q_{i,m}}$ be as in Definition \ref{definition of limit and colimit objects}. Then:

     \begin{enumerate}
         \item For any rational number $\lambda\in\mathbb{Q}$, the functors $\Hom(\mathcal{O}(\lambda),-)$ and $\Hom(-,\mathcal{O}(\lambda))$ take a minimal triangle into a short exact sequence of vector spaces over $\mathbb{Q}_p$.

         \item For any nontrivial morphism $e\in \Hom(\mathcal{O}(\beta_{i,m}),\mathcal{O}(\beta_{j,n}))$, it is either an injection (if $q_{i,m}\leq q_{j,n}$) or a surjection (if $q_{i,m}\geq q_{j,n}$) in $\Coh(X_{FF})$. 

         \item Fix nonzero morphisms $\{f_{2i+1,n}\}$ and $\{g_{2i,m}\}$, and $h_{-1,0}$, where $$f_{2i+1,n}\in \Hom(\mathcal{O}(\beta_{2i+1,n}),\mathcal{O}(\beta_{2i+1,n-1}))$$ $$g_{2i,m}\in \Hom(\mathcal{O}(\beta_{2i,m}),\mathcal{O}(\beta_{2i,m+1})),$$ and $h_{-1,0}\in \Hom(\mathcal{O}(\beta_{0}),\mathcal{O}(\beta_{-1}))$ for any $0<n\leq a_{2i+3}$ and $0\leq m\leq a_{2i+2}-1$. These morphisms correspond to the set of exterior edges. Then there exists nonzero morphisms $$h_{2i,m}\in \Hom(\mathcal{O}(\beta_{2i,m}),\mathcal{O}(\beta_{2i+1})), $$ and $$ h_{2i-1,n}\in \Hom(\mathcal{O}(\beta_{2i}),\mathcal{O}(\beta_{2i-1,n}))$$ for any 
$i\geq 0$, $0< m\leq a_{2i+2}-1$ and $0<n\leq a_{2i+1}-1$, corresponding to the interior edges such that every face of the resulting diagram commutes. 

Moreover, if $f'\in \Hom(\mathcal{O}(\beta_{i,m}),\mathcal{O}(\beta_{j,n}))$ is a composition of the morphisms belonging to the families $\{f_{2i+1,n}\}$, $\{g_{2i,m}\}$ and $\{h_{k,\ell}\}$, we have the following short exact sequences: \begin{itemize}
             \item $$0\rightarrow\mathcal{O}(\beta_{i,m})\xrightarrow{f'} \mathcal{O}(\beta_{j,n})\rightarrow \mathcal{O}(P_{\beta_{i,m},\beta_{j,n}}) \rightarrow 0$$ if $q_{j,n}>q_{i,m}$;

             \item $$0\rightarrow \mathcal{O}(P_{\beta_{i,m},\beta_{j,n}}) \rightarrow \mathcal{O}(\beta_{i,m})\xrightarrow{f'} \mathcal{O}(\beta_{j,n})\rightarrow 0$$ if $q_{j,n}<q_{i,m}$.
             
         \end{itemize}
         \item For any set  of nontrivial morphisms $$h_{2i,m}\in \Hom(\mathcal{O}(\beta_{2i,m}),\mathcal{O}(\beta_{2i+1})), $$ and $$ h_{2i+1,n}\in \Hom(\mathcal{O}(\beta_{2i}),\mathcal{O}(\beta_{2i+1,n})), $$ for any 
$i\geq 0$, $0< m\leq a_{2i+2}-1$ and $0<n\leq a_{2i+1}-1$, there exists a set of nontrivial morphisms  $$f_{2i+1,n}\in \Hom(\mathcal{O}(\beta_{2i+1,n}),\mathcal{O}(\beta_{2i+1,n-1}))$$ and $$g_{2i,m}\in \Hom(\mathcal{O}(\beta_{2i,m}),\mathcal{O}(\beta_{2i,m+1}))$$ for any $0<n\leq a_{2i+3}$ and $0\leq m\leq a_{2i+2}-1$, such that the directed Farey graph $F_{\theta,\frac{1}{0}}$ is commutative. 
     \end{enumerate}
\end{theorem}

\begin{proof}
    Part (1) follows easily from Proposition \ref{prop:homandext1}, as  we either have $$\Hom(\mathcal{O}(\lambda),\mathcal{O}(\mu))=0$$ or $$\Ext^1(\mathcal{O}(\lambda),\mathcal{O}(\mu))=0$$ for any $\mu\in\mathbb{Q}_{\infty}$ (as $\lambda\in\mathbb{Q}$ not equal to $\frac{1}{0}$). 

    Part (2) follows from the fact there are no integral points in the interior of the roller coaster and the argument in the proof of Theorem \ref{Theorem between minimal triangles and Farey triangles}.

    Part (4) and the first half of part (3) can be proved similarly, so we only need to prove part (3).

    We prove the first half of part (3) by induction. We present one case of the induction step; the base case and the other case are similar.

    Suppose we have already constructed a set of nontrivial morphisms up to $h_{2i,m}$, where $0< m\leq a_{2i+2}$, such that the edges below $h_{2i,m}$ form a commutative diagram. It suffices to show that there exists $h_{2i,m+1}$ such that the following diagram is commutative \[\begin{tikzcd}
        & \mathcal{O}(\beta_{2i,m+1}) \arrow[ld, dotted, "h_{2i,m+1}"]& \\  \mathcal{O}(\beta_{2i+1})  & & \mathcal{O}(\beta_{2i,m}) \arrow [ll,"h_{2i,m}"] \arrow[lu,"g_{2i,m}"]
    \end{tikzcd}\]

    This follows from applying the functor $\Hom(-,\mathcal{O}(\beta_{2i+1}))$ to the following minimal triangle $$0\rightarrow\mathcal{O}(\beta_{2i,m})\xrightarrow[]{g_{2i,m}}\mathcal{O}(\beta_{2i,m+1})\rightarrow \mathcal{O}(\beta_{2i+1})\rightarrow 0, $$ we obtain the following surjective map 
    
     \[\begin{tikzcd}[column sep=large] 
     \Hom(\mathcal{O}(\beta_{2i,m+1}),\mathcal{O}(\beta_{2i+1})) \arrow[r,twoheadrightarrow, "\circ g_{2i,m}"] & \Hom(\mathcal{O}(\beta_{2i,m}),\mathcal{O}(\beta_{2i+1}))
     \end{tikzcd}\]
    as $\Ext^1(\mathcal{O}(\beta_{2i+1}),\mathcal{O}(\beta_{2i+1}))=0$ ($\beta_{2i+1}\neq \frac{1}{0}$ as $i\geq 0$). Hence, we have proved the existence of $h_{2i,m+1}$ and the first half of part (3).

   For the second half of part (3),  one can easily show that the unique shortest path $P_{\beta_{i,m},\beta_{j,n}}$ is either concave up if $q_{i,m}>q_{j,n}$ or concave down if $q_{i,m}<q_{j,n}$. By the concavity of the path $P_{\beta_{i,m},\beta_{j,n}}$ and Proposition \ref{prop:homandext1}, an induction along the path shows that one direct summand of $\ker(f')$, respectively $\coker(f')$, splits off at each edge of the path $P_{\beta_{i,m},\beta_{j,n}}$. Hence we know that either $\ker(f')\simeq \mathcal{O}(P_{\beta_{i,m},\beta_{j,n}})$ or $\coker(f')\simeq \mathcal{O}(P_{\beta_{i,m},\beta_{j,n}})$.
\end{proof}

\begin{remark}
  Statements (3) and (4) basically mean that any assignment of nonzero morphisms to the exterior edges extends to the interior edges, and conversely any assignment to the interior edges extends to the exterior edges. Such extensions are far from unique.

  Although this theorem does not hold in the elliptic curve case verbatim, there are similar results which hold in the elliptic curve case. Indeed, for generic choices of vector bundles on the elliptic curve $E$ for each semi-convergents $\beta_{i,m}$, the theorem still holds, since $\Ext^1(L,L')=0$ for two different holomorphic line bundles $L,L'\in \mathrm{Pic}^0(E)$.
\end{remark}

\section{Homological algebra of continuum envelopes}\label{section: homological algebra}
In this section, we will study the homological algebra of the continuum envelope $\QCoh_{\mathbb{R}}(X_{FF})$. 

\subsection{Hom-spaces and endomorphism algebras} Theorem \ref{theorem about the roller coaster}.(3) yields the following proposition concerning $$\Hom(\mathcal{O}(\theta^-,\{g_{2i}\}), \mathcal{O}(\theta^+,\{f_{2j+1}\})).$$

\begin{prop}\label{prop:C_p dimension 1}
    Let $\{g_{2i}\}$ and $\{f_{2j+1}\}$ be any two sequences of nonzero morphisms. Then there is the following short exact sequence  $$0\rightarrow \Pi_{i=0}^{\infty} \Hom(\mathcal{O}(\beta_i),\mathcal{O}({\beta_i}))^{a_{i+1}}\rightarrow \Hom(\mathcal{O}(\theta^-,\{g_{2i}\}), \mathcal{O}(\theta^+,\{f_{2j+1}\}))\xrightarrow{\vartheta} C\rightarrow 0$$
\end{prop}

\begin{proof}
    This follows from the proof of Theorem \ref{theorem about the roller coaster}.(3). Indeed, for any choices of the families $\{f_{2j+1}\}$ and $\{g_{2i}\}$, consider the corresponding colimit and limit objects \[\begin{tikzcd}
    \varinjlim(\mathcal{O}(\beta_0)\arrow[r,hookrightarrow,"g_0"]     &\mathcal{O}(\beta_2)\arrow[r,hookrightarrow,"g_2"]  & \cdots \mathcal{O}(\beta_{2i})\arrow[r,hookrightarrow,"g_{2i}"] & \cdots) = \mathcal{O}(\theta^-,\{g_{2i}\}), \\ 
   \varprojlim(\mathcal{O}(\beta_{-1}) & \mathcal{O}(\beta_1)\arrow[l,twoheadrightarrow,"f_1"] & \cdots \mathcal{O}(\beta_{2j+1})\arrow[l,twoheadrightarrow,"f_3"] &\cdots \arrow[l,twoheadrightarrow,"f_{2j+1}"] )=\mathcal{O}(\theta^+,\{f_{2j+1}\}),
    \end{tikzcd}\]
     and any morphism  $\alpha$ in $$\Hom(\mathcal{O}(\theta^-,\{g_{2i}\}), \mathcal{O}(\theta^+,\{f_{2j+1}\}))=\varprojlim_{f_{2j+1}\circ}\varprojlim_{\circ g_{2i}}\Hom(\mathcal{O}(\beta_{2i}),\mathcal{O}(\beta_{2j+1})),$$ we obtain a morphism $\alpha_{0,-1}$ in $\Hom(\mathcal{O}(\beta_0),\mathcal{O}(\beta_{-1}))$, this defines the map $\vartheta$ in the short exact sequence. 

    Moreover, for any given morphism $\alpha_{0,-1} \in \Hom(\mathcal{O}(\beta_0),\mathcal{O}(\beta_{-1}))$, by Theorem \ref{theorem about the roller coaster}.(3), there exist morphisms $\alpha_{2i,2i-1}\in \Hom(\mathcal{O}(\beta_{2i}),\mathcal{O}(\beta_{2i-1}))$ and $\alpha_{2i,2i+1}\in \Hom(\mathcal{O}(\beta_{2i}),\mathcal{O}(\beta_{2i+1}))$ such that the following diagram is commutative. 

    \[\begin{tikzcd}
   \mathcal{O}(\beta_0)\arrow[r,hookrightarrow,"g_0"]    \arrow[d,"\alpha_{0,-1}"]\arrow[rd,"\alpha_{0,1}"] &\mathcal{O}(\beta_2)\arrow[r,hookrightarrow,"g_2"]  \arrow[d,"\alpha_{2,1}"] \arrow[rd,"\alpha_{2,3}"]& \cdots \mathcal{O}(\beta_{2i})\arrow[r,hookrightarrow,"g_{2i}"] \arrow[d,"\alpha_{2i,2i-1}"] \arrow[rd,"\alpha_{2i,2i+1}"] & \cdots &= \mathcal{O}(\theta^-,\{g_{2i}\}) \\ \mathcal{O}(\beta_{-1}) & \mathcal{O}(\beta_1)\arrow[l,twoheadrightarrow,"f_1"] & \cdots \mathcal{O}(\beta_{2i-1})\arrow[l,twoheadrightarrow,"f_3"] &\cdots \arrow[l,twoheadrightarrow,"f_{2i-1}"] &=\mathcal{O}(\theta^+,\{f_{2j+1}))
    \end{tikzcd}\]
The two families of morphisms $\{\alpha_{2i,2i-1}\}$ and $\{\alpha_{2i,2i+1}\}$ determine a morphism $$\alpha \in \Hom(\mathcal{O}(\theta^-,\{g_{2i}\}), \mathcal{O}(\theta^+,\{f_{2j+1}))$$ with $\vartheta(\alpha)= \alpha_{0,-1}$. Hence, this proves the surjectivity of $\vartheta$.

For the kernel of $\vartheta$, consider the following short exact sequence provided by Theorem \ref{theorem about the roller coaster}(3). $$0\rightarrow\mathcal{O}({\beta_0})^{a_1}\rightarrow \mathcal{O}(\beta_1)\xrightarrow{f_1} \mathcal{O}(\beta_{-1})\rightarrow 0.$$ Applying the functor $\Hom(-,\mathcal{O}(\beta_0))$ to this short exact sequence gives $$0\rightarrow \Hom(\mathcal{O}(\beta_0),\mathcal{O}(\beta_0))^{a_1}\rightarrow \Hom(\mathcal{O}(\beta_0),\mathcal{O}(\beta_{1}))\xrightarrow{f_1\circ} \Hom(\mathcal{O}(\beta_0),\mathcal{O}(\beta_{-1}))\rightarrow 0,$$ i.e., $\Hom(\mathcal{O}(\beta_0),\mathcal{O}(\beta_0))^{a_1}$ is isomorphic to the subspace whose composition with $f_1$ is zero. 

Then we argue inductively as in the proof of (3) in Theorem \ref{theorem about the roller coaster}. We get the short exact sequence.
\end{proof}

\begin{remark}
    The short exact sequence in the proposition can be viewed as a limit of short exact sequences of $p$-adic Banach spaces, and $\vartheta$ is the Fontaine $\vartheta$ map defined in \S\ref{subsection:Fargues--Fontaine curve}.

    Any nontrivial morphism in $\Hom(\mathcal{O}(\theta^-,\{g_{2i}\}), \mathcal{O}(\theta^+,\{f_{2j+1}\}))$ is an injection in $\QCoh(X_{FF})$.
\end{remark}

A direct calculation shows that, for any two families of nonzero morphisms $\{g_{2i}\}$ and $\{f'_{2j+1}\}$, and $\theta>\theta',$ the Hom-space $$\Hom(\mathcal{O}(\theta^-,\{g_{2i}\}), \mathcal{O}(\theta'^+,\{f'_{2j+1}))=\varprojlim_{f'_{2j+1}\circ}\varprojlim_{\circ g_{2i}}\Hom(\mathcal{O}(\beta_{2i}),\mathcal{O}(\beta'_{2j+1}))=0$$ because $\beta_{2i}>\beta'_{2j+1}$ for all $i,j \gg 0$.

Next, for two possibly different families of nonzero morphisms $\{g_{2i}\}$ and $\{g'_{2i}\}$, consider the Hom-space $$\Hom(\mathcal{O}(\theta^-,\{g_{2i}\}), \mathcal{O}(\theta^-,\{g'_{2i}\})).$$ 

\begin{theorem}\label{thm about colimits}
    
 Fix an irrational number $\theta$. Let $\{g_{2i}\}$ and $\{g'_{2i}\}$ be two families of nonzero morphisms, and let $$\mathcal{O}(\theta^-,\{g_{2i}\}), \mathcal{O}(\theta^-,\{g'_{2i}\}))$$ be the corresponding colimit objects, respectively. Then the following conditions are equivalent. 
\begin{enumerate}
    \item There exists a nontrivial morphism $$\alpha\in \Hom(\mathcal{O}(\theta^-,\{g_{2i}\}), \mathcal{O}(\theta^-,\{g'_{2i}\})).$$
    \item There exist an integer $N>0$ and nonzero morphisms $h_{2i}$ for all $i\geq N$, such that the following diagram is commutative  \[\begin{tikzcd}
   \mathcal{O}(\beta_{2N})\arrow[r,hookrightarrow,"g_{2N}"]  \arrow[d,"h_{2N}"]   &\mathcal{O}(\beta_{2N+2})\arrow[r,hookrightarrow,"g_{2N+2}"]  \arrow[d,"h_{2N+2}"]& \cdots \mathcal{O}(\beta_{2i})\arrow[r,hookrightarrow,"g_{2i}"] \arrow[d,"h_{2i}"]& \cdots = \mathcal{O}(\theta^-,\{g_{2i}\}) \\ \mathcal{O}(\beta_{2N})\arrow[r,hookrightarrow,"g_{2N}'"]     &\mathcal{O}(\beta_{2N+2})\arrow[r,hookrightarrow,"g_{2N+2}'"]  & \cdots \mathcal{O}(\beta_{2i})\arrow[r,hookrightarrow,"g_{2i}'"] & \cdots = \mathcal{O}(\theta^-,\{g_{2i}'\})
    \end{tikzcd}\]

    \item The colimit objects $\mathcal{O}(\theta^-,\{g_{2i}\}), \mathcal{O}(\theta^-,\{g'_{2i}\})$ are isomorphic.

\end{enumerate}
\end{theorem}
\begin{proof}
    To show that (1) implies (2), let $\alpha\in \Hom(\mathcal{O}(\theta^-,\{g_{2i}\}), \mathcal{O}(\theta^-,\{g'_{2i}\}))$ be a nonzero morphism. Since 
    \begin{equation*}
    \begin{split}
        \Hom(\mathcal{O}(\theta^-,\{g_{2i}\}), \mathcal{O}(\theta^-,\{g'_{2i}\}))&= \varprojlim_{g_{2i}\circ} \Hom(\mathcal{O}(\beta_{2i}),\mathcal{O}(\theta^-,\{g'_{2i}\})) \\ &=\varprojlim_{g_{2i}\circ}\varinjlim_{\circ g_{2j}'}\Hom(\mathcal{O}(\beta_{2i}),\mathcal{O}(\beta_{2j}))
         \end{split}
    \end{equation*}
    The last equality follows from the fact that the vector bundle $\mathcal{O}(\beta_{2i})$ is a compact object in $\QCoh(X_{FF})$. And change the index of $\{g_{2i}'\}$ to $\{g_{2j}'\}$.

    Therefore, a morphism $\alpha\in \Hom(\mathcal{O}(\theta^-,\{g_{2i}\}), \mathcal{O}(\theta^-,\{g'_{2i}\}))$ determines a compatible family of morphisms $\{\alpha_{i}\}$, where $$\alpha_i\in \Hom(\mathcal{O}(\beta_{2i}),\mathcal{O}(\theta^-,\{g'_{2i}\})).$$

 Since $\alpha\neq 0$, there exists $N_0$ such that its component $\alpha_{N_0}\in \Hom(\mathcal{O}(\beta_{2N_0}), \mathcal{O}(\theta^-,\{g_{2i}'\}))$ is nonzero. Because $$\Hom(\mathcal{O}(\beta_{2N_0}), \mathcal{O}(\theta^-,\{g_{2i}'\}))=\varinjlim_{\circ g_{2i}'}\Hom(\mathcal{O}(\beta_{2N_0}),\mathcal{O}(\beta_{2i})),$$ we know that there exists an integer $N_1$ and a nonzero morphism $$\alpha_{N_0,N_1}\in \Hom(\mathcal{O}(\beta_{2N_0}),\mathcal{O}(\beta_{2N_1}))$$ such that the following diagram is commutative.\[\begin{tikzcd}
\cdots \mathcal{O}(\beta_{2N_0})\arrow[r,hookrightarrow,"g_{2N_0}"]\arrow[rrd,"\alpha_{N_0}"] \arrow[d,"\alpha_{N_0,N_1}"]& \cdots &= \mathcal{O}(\theta^-,\{g_{2i}\}) \arrow[d,"\alpha"] \\  \cdots \mathcal{O}(\beta_{2N_1})\arrow[r,hookrightarrow,"g_{2N_1}'"] & \cdots &= \mathcal{O}(\theta^-,\{g_{2i}'\})
 \end{tikzcd}\]

 For the morphism $\alpha_{N_0+1}\in \Hom(\mathcal{O}(\beta_{2N_0+2}),\mathcal{O}(\theta^-,\{g_{2i}'\}))$, we can find an integer $N_2$ such that the following diagram is commutative (without loss of generality, we assume that $N_2>N_1>N_0$). \[\begin{tikzcd}
     \cdots \mathcal{O}(\beta_{2N_0})\arrow[r,hookrightarrow,"g_{2N_0}"] \arrow[d,"\alpha_{N_0,N_1}"]& \mathcal{O}(\beta_{2N_0+2}) \arrow[r,hookrightarrow, "g_{2N_0+2}"] \arrow[rd,"\alpha_{N_0+1,N_2}"]&\cdots & & = \mathcal{O}(\theta^-,\{g_{2i}\}) \arrow[d,"\alpha"] \\  \cdots \mathcal{O}(\beta_{2N_1})\arrow[r,hookrightarrow,"g_{2N_1}'"] & \cdots\arrow[r,hookrightarrow, "g_{2N_2-2}'"] & \mathcal{O}(\beta_{2N_2}) \arrow[r,hookrightarrow] & \cdots &= \mathcal{O}(\theta^-,\{g_{2i}'\})
 \end{tikzcd}\]

        The left-hand part of the diagram is commutative because the two relevant composites with the natural morphism $$\cdots\circ g_{2N_2+2}'\circ g_{2N_2}':\mathcal{O}(\beta_{2N_2})\rightarrow \mathcal{O}(\theta^-,\{g_{2i}'\})$$ are equal and this morphism is a monomorphism in $\QCoh(X_{FF})$.
        Hence we have the following commutative diagram \[\begin{tikzcd}
            \mathcal{O}(\beta_{2N_0})\arrow[r,hookrightarrow,"g_{2N_0}"] \arrow[d,"g_{2N_2-4}'\circ\cdots\circ\alpha_{N_0,N_1}"]& \mathcal{O}(\beta_{2N_0+2})\arrow[d,"\alpha_{N_0+1,N_2}"] \\ 
            \mathcal{O}(\beta_{2N_2-2}) \arrow[r,hookrightarrow, "g_{2N_2-2}' "] & \mathcal{O}(\beta_{2N_2}) 
        \end{tikzcd}\]
    This commutative diagram extends to the following commutative diagram of short exact sequences by Theorem \ref{theorem about the roller coaster}(3). \[\begin{tikzcd}
         0\arrow[r] &   \mathcal{O}(\beta_{2N_0})\arrow[r,hookrightarrow,"g_{2N_0}"] \arrow[d,"g_{2N_2-4}'\circ\cdots\circ\alpha_{N_0,N_1}"]& \mathcal{O}(\beta_{2N_0+2})\arrow[d,"\alpha_{N_0+1,N_2}"] \arrow[r] & \mathcal{O}(\beta_{2N_0+1})^{a_{2N_0+2}}\arrow {r} \arrow[d,"h"]& 0\\ 0\arrow[r] &
            \mathcal{O}(\beta_{2N_2-2}) \arrow[r,hookrightarrow, "g_{2N_2-2}' "] & \mathcal{O}(\beta_{2N_2}) \arrow[r] & \mathcal{O}(\beta_{2N_2-1})^{a_{2N_2}}\arrow {r} & 0
        \end{tikzcd}\]

        As $N_2>N_0+1$, we know that $\beta_{2N_0+1}>\beta_{2N_2-1}$ and hence $$\Hom(\mathcal{O}(\beta_{2N_0+1}),\mathcal{O}(\beta_{2N_2-1}))=0,$$  which shows that the morphism $h$ in the previous diagram is $0$. Thus we obtain a nonzero morphism $\alpha_{N_0+1,N_2-1}$ and the  following commutative diagram. 
        \[\begin{tikzcd}[column sep=large]
            \mathcal{O}(\beta_{2N_0})\arrow[r,hookrightarrow,"g_{2N_0}"] \arrow[d,swap,"g_{2N_2-4}'\circ\cdots\circ\alpha_{N_0,N_1}"]& \mathcal{O}(\beta_{2N_0+2})\arrow[d,"\alpha_{N_0+1,N_2}"] \arrow[ld,swap,"\alpha_{N_0+1,N_2-1}"]\\ 
            \mathcal{O}(\beta_{2N_2-2}) \arrow[r,hookrightarrow, "g_{2N_2-2}' "] & \mathcal{O}(\beta_{2N_2}) 
        \end{tikzcd}\]
      Inductively using the same argument, we get a nontrivial morphism $\alpha_{N_1,N_1}: \mathcal{O}(\beta_{N_1})\rightarrow \mathcal{O}(\beta_{N_1})$ which is compatible with all these nontrivial morphisms.   Continuing to argue as above, we get the commutative diagram in (2). This completes the proof of (1) implies (2). 

      It is obvious that (3) implies (1). We only need to prove that (2) implies (3). As the Hom-space  $\Hom(\mathcal{O}(\beta_{2i}), \mathcal{O}(\beta_{2i}))$ is the central division algebra associated with the element $\beta_{2i}\in \mathbb{Q}/\mathbb{Z}=\mathrm{Br}(\mathbb{Q}_p)$ (see \cite[Proposition 8.2.8]{farguesfontaine-courbes}), we know that every nontrivial morphism $h_{2i}$ has an inverse morphism $h_{2i}^{-1}$, and it is easy to see the following diagram is commutative. 
      \[\begin{tikzcd}
   \mathcal{O}(\beta_{2N})\arrow[r,hookrightarrow,"g_{2N}'"]  \arrow[d,"h_{2N}^{-1}"]   &\mathcal{O}(\beta_{2N+2})\arrow[r,hookrightarrow,"g_{2N+2}'"]  \arrow[d,"h_{2N+2}^{-1}"]& \cdots \mathcal{O}(\beta_{2i})\arrow[r,hookrightarrow,"g_{2i}'"] \arrow[d,"h_{2i}^{-1}"]& \cdots = \mathcal{O}(\theta^-,\{g_{2i}'\}) \\ \mathcal{O}(\beta_{2N})\arrow[r,hookrightarrow,"g_{2N}"]     &\mathcal{O}(\beta_{2N+2})\arrow[r,hookrightarrow,"g_{2N+2}"]  & \cdots \mathcal{O}(\beta_{2i})\arrow[r,hookrightarrow,"g_{2i}"] & \cdots = \mathcal{O}(\theta^-,\{g_{2i}\})
    \end{tikzcd}\]

     These two sets of nontrivial morphisms $\{h_{2i}\}$ and $\{h_{2i}^{-1}\}$ give us two inverse nontrivial morphisms $$\alpha\in \Hom(\mathcal{O}(\theta^-,\{g_{2i}\}), \mathcal{O}(\theta^-,\{g'_{2i}\}))$$ and $$\alpha^{-1}\in \Hom(\mathcal{O}(\theta^-,\{g'_{2i}\}), \mathcal{O}(\theta^-,\{g_{2i}\})).$$ This proves (3).
\end{proof}

\begin{remark}
    We can not prove the similar result  for $\mathcal{O}(\theta^{+}, \{f_{2i+1}\})$ in $\QCoh(X_{FF})$ as we do not know the cocompactness of $\mathcal{O}(\beta_{2i+1})$. The dual statement is, however, true for formal limit object in $\Pro(\Coh(X_{FF}))$.
\end{remark}
As a corollary of this theorem, we have the following result. 
\begin{corollary}\label{corollary:infinite rank vector bundles}
       For any set of nontrivial morphisms $\{g_{2i}\}$, the object $\mathcal{O}(\theta^-, \{g_{2i}\})$ is an indecomposable vector bundle of infinte rank, and the Hom-space $$\Hom(\mathcal{O}(\theta^-, \{g_{2i}\}),\mathcal{O}(\theta^-, \{g_{2i}\}))$$ is a nontrivial division algebra over $\mathbb{Q}_p$.
\end{corollary}

\begin{proof}
Firstly, we will show that the colimit object $\mathcal{O}(\theta^-, \{g_{2i}\})$ is an infinite-dimensional vector bundle on $X_{FF}$ (see \cite{Infinitedimensionalvectorbundles} for the definition of infinite-dimensional vector bundles). Indeed, let $U=\mathrm{Spec}(R)\subset X_{FF}$ be an open affine subscheme of $X_{FF}$, and $M_i=\mathrm{H}^0(U, \mathcal{O}(\beta_{2i}))$ be the corresponding projective modules over $R$ ($M_i$ is a free module because $R$ is a PID). As the functor $\mathrm{H}^0(U,-)$ commutes with colimits (by the affiness of $U$), we get that $$M\coloneqq \mathrm{H}^{0}(U, \mathcal{O}(\theta^-, \{g_{2i}\}))= \varinjlim_{g_{2i}} M_i.$$ Hence, $M$ is a countably generated flat module. By \cite[Theorem 2.2]{Infinitedimensionalvectorbundles}, we only need to prove that $M$ is a Mittag--Leffler module. The functor $\Hom(-,R)$ sends the injective morphisms $\{g_{2i}\}$ to surjective morphisms as $\mathrm{H}^0(U, \coker(g_{2i}))$ is a projective module for any $i$, hence the dual system of $M$ satisfies the Mittag-Leffler condition. Therefore, $M$ is a projective module, and in fact, a free module over $R$ by Kaplansky's theorem.

    By the proof of Theorem \ref{thm about colimits}, one can show that any nontrivial morphism $$\alpha\in \Hom(\mathcal{O}(\theta^-, \{g_{2i}\}),\mathcal{O}(\theta^-, \{g_{2i}\}))$$ is invertible. Hence the statement is true.
\end{proof}

\begin{remark}
Theorem \ref{thm about colimits} holds for elliptic curve case with a slightly different proof. In fact,  the same argument will show that  $$\Hom(\mathcal{O}(\theta^{-},\{\mathcal{L}_{2i}\}, \{g_{2i}\}), \mathcal{O}(\theta^{-},\{\mathcal{L}_{2i}\}, \{g_{2i}\}))=\mathbb{C}$$ for any choice of vector bundles $\mathcal{L}_{2i}$ and morphisms $\{g_{2i}\}$. The only difference is that the cokernel of the morphism $f: \mathcal{L}_{2i}\xrightarrow{g_{2i}} \mathcal{L}_{2i+2}$ does not split as a direct sum of copies of $\mathcal{L}_{2i+1}$. However, we still have that $\coker(f)$ is semistable of slope $\beta_{2i+1}$ (see Proposition \ref{prop: focus on convergents ellpitic curve case}), which is sufficient for the proof.
\end{remark}

\begin{remark}
    According to Theorem \ref{thm about colimits}, one would expect the moduli space of $\mathcal{O}(\theta^-, \{g_{2i}\})$ to be $$ \varinjlim_{j\to\infty}(\Pi_{i=j}^{\infty} \Hom^*(\mathcal{O}(\beta_{2i}),\mathcal{O}(\beta_{2i+2}))/\Pi_{i=j}^{\infty}\mathrm{Aut}(\mathcal{O}(\beta_{2i})) $$ where  $\Hom^*(\mathcal{O}(\beta_{2i}),\mathcal{O}(\beta_{2i+2}))$ denotes the set of nontrivial morphisms in $\Hom(\mathcal{O}(\beta_{2i}),\mathcal{O}(\beta_{2i+2}))$, and $\Pi_{i=j}^{\infty}\mathrm{Aut}(\mathcal{O}(\beta_{2i}),\mathcal{O}(\beta_{2i})))$ acts on $$\Pi_{i=j}^{\infty} \Hom^*(\mathcal{O}(\beta_{2i}),\mathcal{O}(\beta_{2i+2}))$$ in a natural way.  This should be thought of as an Ind-Pro object of $v$-stacks, and it is infinite-dimensional. For an elliptic curve $E$, one can show that the moduli space of quasi-coherent sheaves $$\mathcal{O}(\theta^{-},\{\mathcal{L}_{2i}\}, \{g_{2i}\})$$ is a $\mathbb{C}^{\times}$-gerbe over the Ind-Pro scheme (and is not a scheme) $$\varinjlim_{j\to\infty} E\times \mathbb{P}^{a_{2j}-1}\times E\times \mathbb{P}^{a_{2j+2}-1}\times\cdots,$$ where $\mathbb{P}^{a_{2j}-1}$ is the projectivization of $\Hom(\mathcal{L}_{2j-2},\mathcal{L}_{2j})$. 
\end{remark}
The next lemma shows that the $K$-class of $\mathcal{O}(\theta^-,\{g_{2i}\})$ is independent of the choice of $\{g_{2i}\}$. 
\begin{lemma}\label{K class lemma}
    For any set of nontrivial morphisms  $\{g_{2i}\}$, there exists a short exact sequence $$0\rightarrow \mathcal{O}(\beta_0)\rightarrow \mathcal{O}(\theta^-,\{g_{2i}\})\rightarrow \bigoplus_{j=0}^{\infty}\mathcal{O}(\beta_{2j+1})^{a_{2j+2}}\rightarrow 0.$$
\end{lemma}

\begin{proof}
This short exact sequence is the colimit of a direct system of short exact sequences from Theorem \ref{theorem about the roller coaster}(3). Indeed, we have the following commutative diagram, \[\begin{tikzcd} & & 0 \arrow[d]& 0\arrow[d] & \\
    0\arrow[r] & \mathcal{O}(\beta_0)\arrow[d,"id"]\arrow[r,"g_{2i-2}\circ\cdots\circ g_0"] &\mathcal{O}(\beta_{2i})\arrow[d,"g_{2i}"] \arrow[r]  & \bigoplus_{j=1}^{i-1}\mathcal{O}(\beta_{2j+1})^{a_{2j+2}} \arrow[d,"i"] \arrow[r] & 0\\  0\arrow[r] & \mathcal{O}(\beta_0) \arrow[r,"g_{2i}\circ\cdots\circ g_0"] &\mathcal{O}(\beta_{2i+2})\arrow[d] \arrow[r]  & \coker(g_{2i}\circ\cdots\circ g_0) \arrow[r] \arrow[d] &0 \\ & & \mathcal{O}(\beta_{2i+1})^{a_{2i+2}} \arrow[r,"id"] \arrow[d]& \mathcal{O}(\beta_{2i+1})^{a_{2i+2}}\arrow[d] & \\ & & 0& 0& 
\end{tikzcd}\]
where each row and each column is a short exact sequence in $\Coh(X_{FF})$. As $$\Ext^1(\mathcal{O}(\beta_{2i+1})^{a_{2i+2}}, \bigoplus_{j=1}^{i-1}\mathcal{O}(\beta_{2j+1})^{a_{2j+2}})=0$$ by Theorem \ref{prop:homandext1}, we get that $\coker(g_{2i}\circ\cdots\circ g_0)\simeq \bigoplus_{j=1}^{i}\mathcal{O}(\beta_{2j+1})^{a_{2j+2}}$. The morphism $i$ is the natural inclusion, up to isomorphism. Hence, passing to the limit as $i\to\infty$, we get the short exact sequence in the lemma because the filtered colimit is exact.
\end{proof}

\begin{remark}
    This result does not hold for elliptic curves. In fact, the $K$-class of $$\mathcal{O}(\theta^{-},\{\mathcal{L}_{2i,m}\}, \{g_{2i,m}\})$$ does depend on the choice of vector bundles $\{\mathcal{L}_{2i,m}\}$. 
\end{remark}

\begin{remark}
    We obtain an analogous result for $\mathcal{O}(\theta^+,\{f_{2i+1}\})$ by the same argument, using the Mittag--Leffler condition.
\end{remark}

Let $\theta$ be an irrational number, and let $\theta=[a_0;a_1,a_2,\cdots]$ be its continued fraction representation, $\beta_{i}=\frac{p_i}{q_i}$, $\beta_{i,m}=\frac{p_{i,m}}{q_{i,m}}$ be its associated convergents and semi-convergents. We need the following definition of $c(\theta)$.
\begin{definition}\label{defn: c(theta)}
   Let $c_i(\theta)$ denote the greatest common divisor of $\{q_{2i},q_{2i+2}, q_{2i+4},\cdots \}$, hence $c_i(\theta)$ divides $c_j(\theta)$ for any $j\geq i$. We let $$c(\theta)\coloneqq \begin{cases}
       \lim\limits_{i\to\infty} c_i(\theta), & \text{if the sequence converges}; \\ \infty, & \text{otherwise}. 
   \end{cases}$$
\end{definition}
The following lemma shows that we can replace the denominators $q_{2j}$ of the convergent $\beta_{2j}$ in the definition of $c_i(\theta)$  by $a_{2j}$ for any $j>i$.
\begin{lemma}\label{another expression of gcd}
    Let $\theta$ be an irrational number, and $\theta=[a_0;a_1,a_2,\cdots]$ be its continued fraction representation. Then $$c_i(\theta)=\gcd(q_{2i}, a_{2i+2},a_{2i+4},\cdots ).$$
\end{lemma}
\begin{proof}
    If we use $c_{i,j}(\theta)$ to denote $\gcd(q_{2i},q_{2i+2},\cdots,q_{2j})$ for any $j>i$, then for a fixed integer $i$, the sequence $\{c_{i,j}(\theta)\}_{j=i+1}^{\infty}$ is a decreasing sequence that converges to $c_i(\theta)$. 

    It suffices to show that $\gcd(q_{2i},q_{2i+2},\cdots,q_{2j})=\gcd(q_{2i},a_{2i+2},\cdots,a_{2j})$ for any $j>i$. This follows directly from the simple fact $$\gcd(q_{2i}, q_{2i+2})=\gcd(q_{2i}, q_{2i}+a_{2i+2}q_{2i+1})=\gcd(q_{2i}, a_{2i+2})$$ since $\gcd(q_{2i},q_{2i+1})=1$ (as we have $p_{2i+1}q_{2i}-q_{2i+1}p_{2i}=1$).
\end{proof}

In fact, we can show that irrational numbers $\theta$ with $c(\theta)>1$ are rare.

\begin{lemma}\label{lemma:rare case with bigger dimension}
    The Lebesgue measure of the following set $$\mathcal{L}\coloneqq \{\theta\in\mathbb{R}|\theta\notin\mathbb{Q}\ and \ c(\theta)>1\}$$ is $0$.
\end{lemma}
\begin{proof}

     Consider the Gauss transformation $$T:(0,1)\rightarrow (0,1).$$ There exists a $T$-invariant measure $\mu$ on $(0,1)$ that is equivalent to the Lebesgue measure, and the Gauss transformation $T$ is mixing with respect to $\mu$ (see \cite[Theorem 7.3]{Ergodictheory}), and hence $T^2$ is ergodic with respect to this measure. 

On the other hand, $a_{2n}(\theta)=1$ is equivalent to $T^{2n}(\theta)\in(\frac{1}{2},1)$. As $T^2$ is ergodic, the pointwise ergodic theorem implies that for almost every $\theta\in(0,1)$, there exist infinitely many $n$ such that $T^{2n}(\theta)\in(\frac{1}{2},1)$. By Lemma \ref{another expression of gcd}, no point of $\mathcal{L}_{0}\coloneqq \mathcal{L}\cap (0,1)$ has this property. Hence, we have $\mathrm{Leb}(\mathcal{L})=\mathrm{Leb}(\mathcal{L}_0)=0$.
\end{proof}
\begin{theorem}\label{vanishing theorem}
   Let $\theta$, $\theta'$ be irrational numbers and $\gamma\in\mathbb{Q}_{\infty}$. Let $\{g_{2i}\}$, $\{g'_{2i}\}$, $\{f_{2i+1}\}$, $\{f_{2i+1}'\}$ be corresponding sequences of nonzero morphisms. Then the following results hold. 
\begin{enumerate}
    \item If $c(\theta)$ is finite,  the $\mathbb{Q}_p$-dimension of the division algebra $$\Hom(\mathcal{O}(\theta^-,\{g_{2i}\}),\mathcal{O}(\theta^-,\{g_{2i}\})$$ is finite, and divides $c(\theta)^2$.

    \item We have $$\Hom(\mathcal{O}(\theta^-,\{g_{2i}\}), \mathcal{O}(\theta'^+,\{f_{2i+1}\})=0,$$ if $\theta'<\theta$; and $$\Hom(\mathcal{O}(\theta^-,\{g_{2i}\}), \mathcal{O}(\theta'^-,\{g_{2i}'\})=0,$$ if $\theta'<\theta$.

    \item  We have $$\Ext^i(\mathcal{O}(\theta^-,\{g_{2i}\}), \mathcal{O}(\theta'^+,\{f_{2i+1}\})=0$$ for any $i\geq 1$, if $\theta'> \theta$; and $$\Ext^{i}(\mathcal{O}(\theta^-,\{g_{2i}\}),\mathcal{O}(\theta'^-,\{g_{2i}'\})=0$$ for any $i\geq 1$, if $\theta'\geq \theta$.

    \item We have $$\Hom(\mathcal{O}(\gamma),\mathcal{O}(\theta^-,\{g_{2i}\}))=0$$ and $$\Hom(\mathcal{O}(\gamma),\mathcal{O}(\theta^+,\{f_{2i+1}\}))=0$$ if $\gamma>\theta$;  similarly, $$\Ext^i(\mathcal{O}(\gamma),\mathcal{O}(\theta^-,\{g_{2i}\}))=0,$$ and $$\Ext^i(\mathcal{O}(\gamma),\mathcal{O}(\theta^+,\{f_{2i+1}\}))=0$$ for any $i\geq 1$, if $\gamma<\theta$.

    \item We have $$\Hom(\mathcal{O}(\theta^-,\{g_{2i}\}), \mathcal{O}(\gamma))=0$$  if $\gamma<\theta$; and $$\Ext^i(\mathcal{O}(\theta^-,\{g_{2i}\}), \mathcal{O}(\gamma))=0$$  for any $i\geq 1$, if $\gamma>\theta$.
\end{enumerate}
\end{theorem}
\begin{proof}
    All statements except (1) follow directly from the fact that the Hom functor commutes with the relevant (co)limits and from the compactness of vector bundles. It remains to prove (1). 

    Consider the following diagram \[\begin{tikzcd}
         0\arrow[r] &   \mathcal{O}(\beta_{2N_0})\arrow[r,hookrightarrow,"g_{2N_0}"] \arrow[d,dotted,"h_{2N_0}"]& \mathcal{O}(\beta_{2N_0+2})\arrow[d,dotted,"h_{2N_0+2}"] \arrow[r] & \mathcal{O}(\beta_{2N_0+1})^{a_{2N_0+2}}\arrow {r} \arrow[d,dotted,"h_{2N_0+1}"]& 0\\ 0\arrow[r] &
            \mathcal{O}(\beta_{2N_0}) \arrow[r,hookrightarrow, "g_{2N_0} "] & \mathcal{O}(\beta_{2N_0+2}) \arrow[r] & \mathcal{O}(\beta_{2N_0+1})^{a_{2N_0+2}}\arrow {r} & 0
        \end{tikzcd}\] and we claim that if there exists a nontrivial triple $(h_{2N_0}, h_{2N_0+2}, h_{2N_0+1})$ such that the diagram above is commutative, then any one morphism in the triple uniquely determines the other two. We show only that $h_{2N_0}$ determines $h_{2N_0+2}$ and $h_{2N_0+1}$; the other two cases are similar. 

        In fact, we apply the functor $\Hom(-,\mathcal{O}(\beta_{2N_0+2}))$ to the first row and get the following exact sequence. \begin{equation*}
           \begin{split}
               0&\rightarrow  \Hom(\mathcal{O}(\beta_{2N_0+2}),\mathcal{O}(\beta_{2N_0+2})) \xrightarrow{\circ g_{2N_0}} \Hom(\mathcal{O}(\beta_{2N_0}),\mathcal{O}(\beta_{2N_0+2}))\\ & \rightarrow \Hom(\mathcal{O}(\beta_{2N_0+1})^{a_{2N_0+2}}[-1], \mathcal{O}(\beta_{2N_0+2})) \rightarrow 0. \end{split} \end{equation*} 
               
               Since $\Hom( \mathcal{O}(\beta_{2N_0+1})^{a_{2N_0+2}}, \mathcal{O}(\beta_{2N_0+2}))=0$ (by the inequality $\beta_{2N_0+1}> \beta_{2N_0+2}$) and $\Ext^1(\mathcal{O}(\beta_{2N_0+2}),\mathcal{O}(\beta_{2N_0+2}))=0$, the morphism $$\circ g_{2N_0}:\Hom(\mathcal{O}(\beta_{2N_0+2}),\mathcal{O}(\beta_{2N_0+2}))\rightarrow \Hom(\mathcal{O}(\beta_{2N_0}),\mathcal{O}(\beta_{2N_0+2})) $$ is injective. Hence, this proves the uniqueness of $h_{2N_0+2}$. If we apply the functor $$\Hom(-, \mathcal{O}(\beta_{2N_0+1})^{a_{2N_0+2}}),$$ we get the uniqueness of $h_{2N_0+1}$ in a similar way.

        Hence we know that if there exist an integer $N$ and a sequence of nontrivial morphisms $\{h_{2i}\}$ for any $i\geq N$, such that the following diagram is commutative.  \[\begin{tikzcd}
   \mathcal{O}(\beta_{2N})\arrow[r,hookrightarrow,"g_{2N}"]  \arrow[d,"h_{2N}"]   &\mathcal{O}(\beta_{2N+2})\arrow[r,hookrightarrow,"g_{2N+2}"]  \arrow[d,"h_{2N+2}"]& \cdots \mathcal{O}(\beta_{2i})\arrow[r,hookrightarrow,"g_{2i}"] \arrow[d,"h_{2i}"]& \cdots = \mathcal{O}(\theta^-,\{g_{2i}\}) \\ \mathcal{O}(\beta_{2N})\arrow[r,hookrightarrow,"g_{2N}"]     &\mathcal{O}(\beta_{2N+2})\arrow[r,hookrightarrow,"g_{2N+2}"]  & \cdots \mathcal{O}(\beta_{2i})\arrow[r,hookrightarrow,"g_{2i}"] & \cdots = \mathcal{O}(\theta^-,\{g_{2j}\})
    \end{tikzcd}\]
    Then any term $h_{2i}$ uniquely determines this entire sequence.

    If we use $H_N$ to denote the set of such sequences of nontrivial morphisms $\{h_{2i}\}_{i\geq N}$, it is easy to see that $H_N$ has a structure of division algebra over $\mathbb{Q}_p$. By the previous claim, we know that the natural forgetful map $H_N\rightarrow \Hom(\mathcal{O}(\beta_{2i}),\mathcal{O}(\beta_{2i}))$ is an injection of division algebras over $\mathbb{Q}_p$ for any $i\geq N$. By \cite[Theorem 12.7]{Pierceassociativealgebra}, we know that the dimension of $H_n$ divides the dimension of $\Hom(\mathcal{O}(\beta_{2i}),\mathcal{O}(\beta_{2i}))$, which is $q_{2i}^2$ for every $i\geq N$. Hence, the dimension of $H_N$ divides $c_N(\theta)^2$ by definition.
  By Theorem \ref{thm about colimits}, we know that $$\Hom(\mathcal{O}(\theta^-,\{g_{2i}\}),\mathcal{O}(\theta^-,\{g_{2i}\}))=\varinjlim_{N\to \infty} H_N, $$ where the transition map $H_i\rightarrow H_{i+1}$ is the natural forgetful map, which is an inclusion. By the definition of $c(\theta)$, (1) is proved.
\end{proof}
\begin{example}\label{golden ratio}
    If we take $\theta$ to be $\frac{1+\sqrt{5}}{2}$, the continued fraction is $[1;1,1,\cdots]$, and the convergents $\beta_i=\frac{F_{i+1}}{F_i}$, where $F_i$ is the $i$-th Fibonacci number with $F_{-1}=0, F_0=1$. The minimal triangles in the associated Farey graph are $$0\rightarrow \mathcal{O}(\frac{F_{2i+1}}{F_{2i}})\rightarrow \mathcal{O}(\frac{F_{2i+3}}{F_{2i+2}})\rightarrow\mathcal{O}(\frac{F_{2i+2}}{F_{2i+1}})\rightarrow 0$$ and $$0\rightarrow \mathcal{O}(\frac{F_{2i+1}}{F_{2i}})\rightarrow \mathcal{O}(\frac{F_{2i+2}}{F_{2i+1}})\rightarrow\mathcal{O}(\frac{F_{2i}}{F_{2i-1}})\rightarrow 0.$$

    By Lemma \ref{another expression of gcd}, we know that $c(\frac{1+\sqrt{5}}{2})=1$, and hence $$\Hom(\mathcal{O}(\frac{1+\sqrt{5}}{2}^-,\{g_{2i}\}),\mathcal{O}(\frac{1+\sqrt{5}}{2}^-,\{g_{2i}\}))=\mathbb{Q}_p,$$ for any sequence of nontrivial morphisms $\{g_{2i}\}$ by Theorem \ref{vanishing theorem}.
\end{example}
\begin{remark}
    Example \ref{golden ratio} and the previous theorem show that there is no Riemann-Roch type theorem for $\QCoh_{\mathbb{R}}(X_{FF})$. 

Indeed, for any two sequences of nontrivial morphisms $\{g_{2i}\}$ and $\{g'_{2i}\}$, the $K$-classes  of $\mathcal{O}(\frac{1+\sqrt{5}}{2}^-,\{g_{2i}\})$ and $\mathcal{O}(\frac{1+\sqrt{5}}{2}^-,\{g_{2i}'\})$ are the same by Lemma \ref{K class lemma}. 

However, by Theorem \ref{vanishing theorem}, we know that $$\Hom(\mathcal{O}(\frac{1+\sqrt{5}}{2}^-,\{g_{2i}\}),\mathcal{O}(\frac{1+\sqrt{5}}{2}^-,\{g_{2i}\}))=\mathbb{Q}_p,$$ $$\Ext^i(\mathcal{O}(\frac{1+\sqrt{5}}{2}^-,\{g_{2i}\}),\mathcal{O}(\frac{1+\sqrt{5}}{2}^-,\{g_{2i}\}))=0$$ for any $i\geq 1$ and $$\Ext^i(\mathcal{O}(\frac{1+\sqrt{5}}{2}^-,\{g_{2i}\}),\mathcal{O}(\frac{1+\sqrt{5}}{2}^-,\{g_{2i}'\}))=0$$ for any $i\geq 0$ for any two distinct objects $\mathcal{O}(\frac{1+\sqrt{5}}{2}^-,\{g_{2i}\}),\mathcal{O}(\frac{1+\sqrt{5}}{2}^-,\{g_{2i}'\})$ in their infinite-dimensional moduli v-stack. 

Therefore, the Riemann--Roch theorem could not hold for $\QCoh_{\mathbb{R}}(X_{FF})$ unless we ignore the finite $\mathbb{Q}_p$-dimensional part.
    
\end{remark}

%\begin{corollary}\label{irrational number with bounded numbers in continued fractions}
    %Let $\theta$ be an irrational number, and let $\theta=[a_0;a_1,a_2,\cdots]$ be its continued fraction representation, if the set of $\{a_i\}$ is bounded above, then the division algebra $$\Hom(\mathcal{O}(\theta^-,\{g_{2i}\}),\mathcal{O}(\theta^-,\{g_{2i}\}))$$ is finite-dimensional.
%\end{corollary}
%\begin{proof}
    %This follows directly from Theorem \ref{vanishing theorem}.(1) and Lemma \ref{another expression of gcd}.
%\end{proof}
%\begin{corollary}
    %If $\theta$ is a real quadratic irrational number, then the division algebra $$\Hom(\mathcal{O}(\theta^-,\{g_{2i}\}),\mathcal{O}(\theta^-,\{g_{2i}\})$$ is finite-dimensional.
%\end{corollary}

%\begin{proof}
    %By Lagrange's theorem, we know that the continued fraction of a real quadratic number is eventually periodic. Hence, this corollary follows directly from the previous one.
%\end{proof}

%\begin{example}
    %If we let
     %$$\mathcal{C}_2\coloneqq \left\{[0;a_1,\cdots,a_i,\cdots] \mid \text{$1\leq a_i\leq 2$ for any $i\geq 1$}\right\}\subset [0,1],$$ 
    %it is well known that $\mathcal{C}_2$ is a Cantor set, with Hausdorff dimension around $0.53128\cdots$ (see \cite{computingthehausdorffdimensionofE_2}).

    %For any irrational number $\theta\in\mathcal{C}_2$, the dimension of $$\Hom(\mathcal{O}(\theta^-,\{g_{2i}\}),\mathcal{O}(\theta^-,\{g_{2i}\})$$ could only be $1$, $2$, or $4$.
%\end{example}

By Lemma \ref{lemma:rare case with bigger dimension} and Theorem \ref{vanishing theorem}, we have the following corollary.

\begin{corollary}    For almost every (with respect to Lebesgue measure) irrational number $\theta$, we have that $$\Hom(\mathcal{O}(\theta^-,\{g_{2i}\}),\mathcal{O}(\theta^-,\{g_{2i}\}))\simeq \mathbb{Q}_p$$ for any choice of $\{g_{2i}\}$. 
\end{corollary}

On the other hand, we will show that there are infinitely many irrational numbers $\theta$ for which $c(\theta)=\infty $ and $ \Hom(\mathcal{O}(\theta^-,\{g_{2i}\}),\mathcal{O}(\theta^-,\{g_{2i}\}))$ is an infinite-dimensional $\mathbb{Q}_p$-division algebra.

\begin{theorem}\label{thm: infinite-dimensional division algebra}
There exists an infinite set $S$ of irrational numbers such that for any irrational number $\theta\in S$, there exists a sequence of nonzero morphisms $\{g_{2i}\}$ such that $$ \Hom(\mathcal{O}(\theta^-,\{g_{2i}\}),\mathcal{O}(\theta^-,\{g_{2i}\}))$$ is an infinite-dimensional $\mathbb{Q}_p$-division algebra.
\end{theorem}
\begin{proof}
    Let us construct a special class of continued fractions $$\theta=[a_0;a_1,a_2,\cdots]$$ with convergents $\beta_i\coloneqq \frac{p_i}{q_i}$ such that $\theta$ satisfies the following conditions: \begin{enumerate}

        \item for any prime $P$ dividing $q_{2i-2}$, the prime $P$ divides $a_{2i}$ but $P^2$ does not;

        \item for every $i$, there exists a prime $P_i$ that divides $q_{2i}$, but not $q_{2i-2}$.
    \end{enumerate}

    We construct $a_i$ inductively. We can start with an arbitrary $a_0,a_1,a_2\in \mathbb{N}$ with $a_2$ square free. Assume inductively that we have already constructed $a_{i}$ for all $0\leq i\leq 2k+1$ such that conditions (1) and (2) hold up to $k$. To construct $a_{2k+2}$, factor $q_{2k}$ as $$q_{2k}=P_1^{b_1}\cdots P_n^{b_n}$$ where $P_j$ is a prime divisor of $q_{2k}$ and $b_j\geq 1$ for any $1\leq j\leq n$. We let $$a_{2k+2}=P_1\cdots P_nx$$ where $x$ is a positive integer to be determined. By condition (1), we need $\gcd(P_1\cdots P_n, x)=1$. For condition (2), as $$q_{2k+2}=q_{2k}+a_{2k+2}q_{2k+1}=P_1\cdots P_n(P_1^{b_1-1}\cdots P_n^{b_n-1}+xq_{2k+1}),$$ we let $P$ be a prime number which does not divide $q_{2k}q_{2k+1}$, then if $$x\equiv-q_{2k+1}^{-1}P_1^{b_1-1}\cdots P_n^{b_n-1} \ (mod\ P),$$ then $P$ divides $q_{2k+2}$. 

    By the Chinese remainder theorem, one can find such $x$ with $\gcd(P_1\cdots P_n, x)=1$. We may take each $a_{2i+1}$ to be any positive integer. Continuing this construction then yields infinitely many irrational numbers that satisfy conditions (1) and (2).

    We claim that for any continued fraction which satisfies conditions (1) and (2), there exists a sequence of nontrivial morphisms $\{g_{2i}\}$ such that the division algebra $$ \Hom(\mathcal{O}(\theta^-,\{g_{2i}\}),\mathcal{O}(\theta^-,\{g_{2i}\}))$$ is an infinite-dimensional $\mathbb{Q}_p$-division algebra.

    To prove the claim, let us consider the $$\Hom(\mathcal{O}(\beta_{2i}),\mathcal{O}(\beta_{2i}))-\Hom(\mathcal{O}(\beta_{2i+2}),\mathcal{O}(\beta_{2i+2}))$$ bimodule structure of $\Hom(\mathcal{O}(\beta_{2i}),\mathcal{O}(\beta_{2i+2}))$. 

    We use $K_n$ to denote the cyclic field extension of $\mathbb{Q}_p$ of degree $n$. Recall that $\Hom(\mathcal{O}(\beta_{2i}),\mathcal{O}(\beta_{2i}))$ is a  central simple division algebra over $\mathbb{Q}_p$ with $K_{q_{2i}}$ as a maximal subfield; its associated Brauer class is $\frac{p_{2i}}{q_{2i}}\in \mathbb{Q}/\mathbb{Z}$. 

    Set $$d_i\coloneqq \gcd(q_{2i}, q_{2i+2})= \gcd(q_{2i},a_{2i+2})=\gcd(a_{2i+2}, q_{2i+2}).$$ We study its Hom-space  $\Hom(\mathcal{O}(\beta_{2i}),\mathcal{O}(\beta_{2i+2}))$ via the equivalence between the category of $\psi$-modules over $B^+=\cap_{n\geq 0}\psi^n(B^+_{cris})$ and the category of vector bundles on $X_{FF}$ (see \cite[\S 11.4]{farguesfontaine-courbes}). 

    Under this equivalence, $\mathcal{O}(\beta_{2i})$ corresponds to a $\psi$-module $B_{2i}$ whose underlying space is $\oplus_{i=1}^{q_{2i}} B^+e_i$, where $\{e_i\}$ is a basis of $B_{2i}$, and the semilinear morphism $\psi$ acts as $\psi(e_i)=e_{i+1}$ for $1\leq i<q_{2i}$, $\psi(e_{q_{2i}})=p^{p_{2i}}e_1$. We have a similar description for $B_{2i+2}$. 

    Moreover, we have $$\Hom(\mathcal{O}(\beta_{2i}),\mathcal{O}(\beta_{2i+2}))= \Hom_{B^+,\psi}(B_{2i},B_{2i+2}),$$ and the latter space can be described as a subspace of the matrix space $M_{q_{2i+2}\times q_{2i}}(B^+)$ consisting of the matrices $$\begin{pmatrix} b_{11} & \cdots & b_{1q_{2i}} \\ b_{21} & \cdots & b_{2q_{2i}}\\ \cdots & \cdots & \cdots \\ b_{q_{2i+2}1} & \cdots & b_{q_{2i+2}q_{2i}}
        
    \end{pmatrix}$$ such that \begin{equation*}
        \begin{array}{cc}

    \begin{pmatrix} b_{11} & \cdots & b_{1q_{2i}} \\ b_{21} & \cdots & b_{2q_{2i}}\\ \cdots & \cdots & \cdots \\ b_{q_{2i+2}1} & \cdots & b_{q_{2i+2}q_{2i}} \end{pmatrix}\begin{pmatrix}
        0 & 0& \cdots & p^{p_{2i}} \\ 1& 0& \cdots & 0\\ \cdots & \cdots & \cdots & \cdots \\ 0 &\cdots & 1 &0
    \end{pmatrix} \\ =\begin{pmatrix}
        0 & 0& \cdots & p^{p_{2i+2}} \\ 1& 0& \cdots & 0\\ \cdots & \cdots & \cdots & \cdots \\ 0 &\cdots & 1 &0
    \end{pmatrix} \begin{pmatrix} \psi(b_{11}) & \cdots & \psi(b_{1q_{2i}}) \\ \psi(b_{21}) & \cdots & \psi(b_{2q_{2i}})\\ \cdots & \cdots & \cdots \\ \psi(b_{q_{2i+2}1}) & \cdots & \psi(b_{q_{2i+2}q_{2i}})
        
    \end{pmatrix} \end{array} \end{equation*} Hence, a matrix in  $\Hom_{B^+,\psi}(B_{2i},B_{2i+2})$ can be written as $$\begin{pmatrix}
        b_{11} & p^{p_{2i+2}}\psi(b_{q_{2i+2},1}) &p^{p_{2i+2}}\psi^2(b_{q_{2i+2}-1, 1 }) &\cdots  & p^{p_{2i+2}}\psi^{q_{2i}-1}(b_{q_{2i+2}-q_{2i}+2,1})\\ b_{21} & \psi(b_{11}) & p^{p_{2i+2}}\psi^2(b_{q_{2i+2}-1,1}) & \cdots & p^{p_{2i+2}}\psi^{q_{2i}-1}(b_{q_{2i+2}-q_{2i}+3,1}) \\ b_{31} & \psi(b_{21}) & \psi^2(b_{11}) &\cdots &\cdots \\ \cdots & \cdots &\cdots &\cdots & \cdots \\ b_{q_{2i+2}+1,1} &\cdots &\cdots &\cdots & \psi^{q_{2i}-1}(b_{11} ) \\ \cdots & \cdots &\cdots &\cdots & \cdots \\ b_{q_{2i+2},1} &\psi(b_{q_{2i+2}-1,1}) & \psi^2(b_{q_{2i+2}-2,1}) &\cdots & \psi^{q_{2i}-1}(b_{q_{2i+2}-q_{2i}+1,1})
    \end{pmatrix}$$ where the following equations hold: $$p^{p_{2i}}b_{k,1}=p^{p_{2i+2}}\psi^{q_{2i}}(b_{q_{2i+2}-q_{2i}+k,1})$$ for $1\leq k\leq q_{2i}$ and $$p^{p_{2i}}b_{k,1}=\psi^{q_{2i}}(b_{k-q_{2i},1})$$ for  $q_{2i}<k\leq q_{2i+2}$. This implies that $$b_{11},b_{21},\cdots b_{d_i,1}\in (B^+)^{\psi^{q_{2i}q_{2i+2}/{d_i}}=p^{a_{2i+2}/d_i}},$$ and all other entries in this matrix are determined by these elements. This implies the following well-known isomorphism. $$\Hom(\mathcal{O}(\beta_{2i}),\mathcal{O}(\beta_{2i+2}))\simeq \mathrm{H}^0(X_{FF}, \mathcal{O}(\frac{a_{2i+2}/d_i}{q_{2i}q_{2i+2}/d_i}))^{\oplus d_i}$$ as $\mathbb{Q}_p$-Banach spaces. A morphism $g\in \Hom(\mathcal{O}(\beta_{2i}),\mathcal{O}(\beta_{2i+2}))$ is said to be of pure type if its matrix representation has entries $b_{21}=b_{31}=\cdots =b_{d_i,1}=0$.  In fact, the matrix representations provide a more suitable way to study the $$\Hom(\mathcal{O}(\beta_{2i}),\mathcal{O}(\beta_{2i}))-\Hom(\mathcal{O}(\beta_{2i+2}),\mathcal{O}(\beta_{2i+2}))$$ bimodule structure of $\Hom(\mathcal{O}(\beta_{2i}),\mathcal{O}(\beta_{2i+2}))$.

     One can show that the central division algebra $$D_{\beta_{2i}}\coloneqq  \Hom(\mathcal{O}(\beta_{2i}),\mathcal{O}(\beta_{2i}))$$ is isomorphic to the following matrix algebra $$\{\begin{pmatrix}
         d_{11} & p^{p_{2i}}\psi(d_{q_{2i},1}) & \cdots & p^{p_{2i}}\psi^{q_{2i}-1}(d_{21}) \\ d_{21} & \psi(d_{11}) &\cdots & p^{p_{2i}}\psi^{q_{2i}-1}(d_{31}) \\ \cdots &\cdots &\cdots &\cdots \\ d_{q_{2i}1} & \cdots &\cdots & \psi^{q_{2i}-1}(d_{11})
     \end{pmatrix}| d_{k,1}\in(B^+)^{\psi^{q_{2i}}=1}\text{for any $1\leq k\leq q_{2i}$} \}.$$

     We have a similar matrix representation for $$D_{\beta_{2i+2}}\coloneqq  \Hom(\mathcal{O}(\beta_{2i+2}),\mathcal{O}(\beta_{2i+2})).$$

     Hence for any nonzero elements $b_{11}\in (B^+)^{\psi^{q_{2i}q_{2i+2}/{d_i}}=p^{a_{2i+2}/d_i}}$ and $a\in (B^+)^{\psi^{d_i}=1}\simeq K_{d_i}$, where $K_{d_i}$ is the unramified field extension of $\mathbb{Q}_p$ of degree $d_i$,  we have \begin{equation*}
         \begin{array}{cc}
           
     \begin{pmatrix}
         b_{11} & 0 &\cdots \\ 0& \psi(b_{11}) & \cdots \\ \cdots &\cdots &\cdots 
     \end{pmatrix}\begin{pmatrix}
         a & 0&\cdots &0 \\ 0&\psi(a) &\cdots & 0\\ \cdots &\cdots & \cdots &\cdots \\ 0 & 0& \cdots & \psi^{q_{2i}-1}(a)
     \end{pmatrix}\\ =\begin{pmatrix}
          a & 0&\cdots &0 \\ 0&\psi(a) &\cdots & 0\\ \cdots &\cdots & \cdots &\cdots \\ 0 & 0& \cdots & \psi^{q_{2i+2}-1}(a)
     \end{pmatrix}\begin{pmatrix}
           b_{11} & 0 &\cdots \\ 0& \psi(b_{11}) & \cdots \\ \cdots &\cdots &\cdots 
     \end{pmatrix}   
         \end{array}\end{equation*}
    as $d_i=gcd(q_{2i}, q_{2i+2})$. Hence, if $g_{2i}\in \Hom(\mathcal{O}(\beta_{2i}),\mathcal{O}(\beta_{2i+2}))$ is of pure type, i.e. the associated entries $b_{k,1}$ are $0$ for all $2\leq k\leq d_i$, then the set of the pairs $(h_{2i}, h_{2i+2})$ that make the following diagram commutative \[\begin{tikzcd}
   \mathcal{O}(\beta_{2i})\arrow[r,hookrightarrow,"g_{2i}"]  \arrow[d,"h_{2i}"]   &\mathcal{O}(\beta_{2i+2})  \arrow[d,"h_{2i+2}"] \\ \mathcal{O}(\beta_{2i})\arrow[r,hookrightarrow,"g_{2i}"]     &\mathcal{O}(\beta_{2i+2})
    \end{tikzcd}\] has the structure of division algebra over $K_{d_i}$.
    %as it is in fact from the push-forward of the finite map (see \cite[Section 8.2]{farguesfontaine-courbes})  $$\pi_{q_{2i}q_{2i+2}/d_i}: X_{C^{\flat}, K_{q_{2i}q_{2i+2}/d_i}}\rightarrow X_{FF}. $$

    Therefore, if we take $g_{2N}$ to be a morphism of pure type for all $N\geq i$,  the division algebra $H_i$ defined in the proof of Theorem \ref{vanishing theorem} is a division algebra over the field $K_{d_i}$ by conditions (1). By conditions (2), we know that the sequence $\{d_i\}$ tends to infinity. Hence the division algebra $$H=\varinjlim_{i \to \infty} H_i$$ is infinite-dimensional.  
\end{proof}
\begin{remark}
    We have shown that the division algebra $$ \Hom(\mathcal{O}(\theta^-,\{g_{2i}\}),\mathcal{O}(\theta^-,\{g_{2i}\}))$$ depends on the Diophantine properties of the irrational number $\theta$ and the choice of morphisms $\{g_{2i}\}$. We expect that different algebraic structures provide a stratification of the infinite-dimensional moduli stack, and the dimension of the endomorphism algebra is upper semi-continuous. The authors hope to further study these questions in future work.

    Note that the division algebra $H$ is a colimit of finite-dimensional $\mathbb{Q}_p$-division algebras, hence, the $\mathbb{Q}_p$-dimension of $H$ is countable. By comparison with Proposition \ref{prop:C_p dimension 1}, we conclude that $\mathcal{O}(\theta^-,\{g_{2i}\})$ is not isomorphic to $\mathcal{O}(\theta^+,\{f_{2i+1}\})$ for any choices of nontrivial morphisms $\{g_{2i}\}$ and $\{f_{2i+1}\}$.
\end{remark}

\subsection{Quotient categories $\mathcal{P}(\phi,\phi+1]/\mathcal{P}(\phi+1)$}\label{Subsection: quotient categories}

In this subsection, we will study the abelian quotient categories $\mathcal{P}(\phi,\phi+1]/\mathcal{P}(\phi+1)$ for $X_{FF}$, where $-\cot(\pi\phi)\eqqcolon \lambda\in \mathbb{Q}\cup \{-\infty\}$. We denote the quotient category by $\mathcal{Q}_{\lambda}$. We will show that it is semi-simple with a unique simple object $S_{\lambda}$, and then study $\End_{\mathcal{Q}_{\lambda}}(S_{\lambda})$, which is a generalization of the Colmez-Fontaine division algebra $\mathscr{C}$ (see \cite[\S 5, \S 9] {BanachColmezspaces}, \cite[\S 7.3]{LeBrasresult}).

\begin{prop}\label{prop: quotient categories}
   Assume that $\lambda\coloneqq-\cot(\pi\phi)\in\mathbb{Q}\cup \{-\infty\}$. Then the quotient category $\mathcal{P}(\phi,\phi+1]/\mathcal{P}(\phi+1)$ is isomorphic to the quotient category $\mathcal{P}[\phi,\phi+1)/\mathcal{P}(\phi)$, and they are semisimple with a single simple object. 
\end{prop}

\begin{proof}
Up to a shift functor $[n]$, we can assume that $\phi\in[0,1)$. And we write $\lambda=\frac{p}{q}$. As in the proof of Corollary \ref{Hearts of FF curves}, we can find a Farey triangle with vertices $\infty>r>s>\lambda$. We use $S_{\lambda}$ to denote the image of $\mathcal{O}(r)$ in the quotient category $\mathcal{P}(\phi,\phi+1]/\mathcal{P}(\phi+1)$. We will prove that this object is the only simple object in the semi-simple quotient category $\mathcal{P}(\phi,\phi+1]/\mathcal{P}(\phi+1)$. A similar argument applies to $\mathcal{P}[\phi,\phi+1)/\mathcal{P}(\phi)$.
    
    As the Harder--Narasimhan filtration splits in $D^b(X_{FF})$, we only need to show that every indecomposable object in $\mathcal{P}(\phi,\phi+1]$ is isomorphic to a direct sums of copies of $S_{\lambda}$ in the quotient category. The indecomposable objects in $\mathcal{P}(\phi,\phi+1]$ consist of $\mathcal{O}(\mu)$ for any $\lambda<\mu\in\mathbb{Q}$, $\mathcal{O}(\eta)[1]$ for any $\lambda\geq \eta\in\mathbb{Q}$, and indecomposable torsion sheaves. For any $\mu \in\mathbb{Q}$, we use $\mathcal{O}_{\mu}$ to denote the corresponding indecomposable object in $\mathcal{P}(\phi,\phi+1]$. Note that this does not include torsion sheaves.

The two vectors $u_1$ and $u_2$ corresponding to $s$ and $\lambda$ form a $\mathbb{Z}$-basis of $\mathbb{Z}^2$. We claim that for any indecomposable object $\mathcal{E}\in \mathcal{P}(\phi,\phi+1]$, if $v(\mathcal{E})=bu_1+au_2$ for $b\in\mathbb{N}$,$a\in\mathbb{Z}$, the image of $\mathcal{E}$ in $\mathcal{Q}_{\lambda}$ is isomorphic to $S_{\lambda}^{\bigoplus b}.$

We prove the claim by induction on $b$; the case $b=0$ is trivial. If $b=1$ and $q\neq 1$, for any such $\mathcal{E}$, we can find a sequence of minimal triangles $$0\rightarrow \mathcal{O}_{\mu_i} \rightarrow\mathcal{O}_{\mu_{i+1}}\rightarrow \mathcal{O(\lambda})[1]\rightarrow 0 $$ in $\mathcal{P}(\phi,\phi+1]$ such that $\mu_0=r$ and $\mu_n=v(\mathcal{E})$. Hence, this proves that $\mathcal{E}\simeq S_{\lambda}$ in $\mathcal{Q}_{\lambda}$ under the assumption $q\neq 1$.

If $q=1$, for any closed point $x\in X_{FF}$, we have the following short exact sequence $$0\rightarrow \mathcal{O}(p)\rightarrow \mathcal{O}(p+1)\rightarrow \kappa_x\rightarrow 0$$ in $\Coh(X_{FF})$. This implies that the skyscraper sheaf $\kappa_x\simeq S_{\lambda}$ in the quotient category $\mathcal{Q}_{\lambda}$. For every other indecomposable object $\mathcal{E}$ with $b=1$, the argument for $q\neq 1$ applies.

For general $b>1$, consider the vectors $$\mathcal{V}\coloneqq \{v_k=v(\mathcal{E})+k(p,q)| k\in\mathbb{Z}\}.$$ Assume that this set does not contain any vector $(n,0)$ with $n\in\mathbb{Z}_{>0}$. For each vector $v\in \mathcal{V}$, we associate an object $\mathcal{O}_v\in\mathcal{P}(\phi,\phi+1]$ as follows.  For any $v=(a_0,a_1)\in\mathcal{V}$ in the upper half-plane, i.e. $a_1>0$, $$\mathcal{O}_v \coloneqq \mathcal{O}(a_0,a_1)=\mathcal{O}(\frac{a_0}{a_1})^{\oplus gcd(a_0,a_1)}.$$ For any $v=(a_0,a_1)\in \mathcal{V}$ in the lower half-plane, i.e. $a_1<0$, $$\mathcal{O}_v\coloneqq \mathcal{O}(-a_0,-a_1)[1]=\mathcal{O}(\frac{-a_0}{-a_1})^{\oplus gcd(a_0,a_1)}[1].$$ 

To prove the claim, it suffices to show that the objects $\mathcal{O}_{v_k}$ are isomorphic to each other in $\mathcal{Q}_{\lambda}$ for any $k\in\mathbb{Z}$, as $\mathcal{O}_{bu_1}\simeq \mathcal{O}_{u_1}^{\oplus b}\simeq S_{\lambda}^{\oplus b}$ in $\mathcal{Q}_{\lambda}$ by induction. By \cite[Theorem 1.1.2]{Extensionsofvectorbundlesonffcurve}, we have the following short exact sequence $$ 0\rightarrow \mathcal{O}_{v_k}\rightarrow \mathcal{O}_{v_{k+1}}\rightarrow \mathcal{O}(\lambda)[1]\rightarrow 0$$ in $\mathcal{P}(\phi,\phi+1]$ for any $k\in\mathbb{Z}$. This proves the claim under the assumption that $\mathcal{V}$ does not contain any integral point on the $x$-axis.

If this set $\mathcal{V}$ contains a vector $(n,0),$ for any indecomposable torsion sheaf $\mathcal{E}$ with $v(\mathcal{E})=(n,0),$ consider the degree $q$ \'etale cover $\pi: X_{q}\rightarrow X_{FF}$ (see \cite[\S 5.6.4]{farguesfontaine-courbes}). There exists an indecomposable torsion sheaf $\Tilde{\mathcal{E}}$ with $\pi_*(\Tilde{\mathcal{E}})\simeq \mathcal{E}$. And there is a short exact sequence $$0\rightarrow \mathcal{O}_{X_q}(p)\rightarrow \mathcal{O}_{X_q}(p+n)\rightarrow \Tilde{\mathcal{E}}\rightarrow 0$$ in $\Coh(X_q)$.  Applying $\pi_*$ to it, we obtain $$0\rightarrow \mathcal{O}(\frac{p}{q})\rightarrow \mathcal{O}(p+n,q)\rightarrow {\mathcal{E}}\rightarrow 0$$ in $\Coh(X_{FF})$. This shows that $\mathcal{E}\simeq \mathcal{O}(p+n,q)$ in the quotient category $\mathcal{Q}_{\lambda}$. By the previous argument, one can prove that $\mathcal{O}(p+n,q)\simeq  S_{\lambda}^{\oplus b}$. Hence, the claim is proved. 

Hence, we have proved that $\mathcal{P}(\phi,\phi+1]/\mathcal{P}(\phi+1)$ and $\mathcal{P}[\phi,\phi+1)/\mathcal{P}(\phi)$ are semi-simple with a single simple object $S_{\lambda}$. One can directly show that $$\mathrm{End}_{\mathcal{P}(\phi,\phi+1]/\mathcal{P}(\phi+1)}(S_{\lambda})\simeq  \mathrm{End}_{\mathcal{P}[\phi,\phi+1)/\mathcal{P}(\phi)}(S_{\lambda})$$ by the definition of the quotient category. Hence, the proof is complete.
\end{proof}

For any $\lambda=\frac{p}{q}\in\mathbb{Q}\cup\{-\infty\}$, we have the following result about $\End_{\mathcal{Q}_{\lambda}}(S_{\lambda})$.
\begin{corollary}
    The endomorphism algebra $\mathrm{End}_{\mathcal{Q}_{\lambda}}(S_{\lambda})$ is a division algebra, and for any $\frac{r}{s}\in \mathbb{Q}$ with $|det\begin{pmatrix}
        r & s \\ p & q
    \end{pmatrix}|=n>0,$ we have a monomorphism of $\mathbb{Q}_p$ algebras $$D_{\frac{r}{s}}\hookrightarrow \mathrm{M}_n(\mathrm{End}_{\mathcal{Q}_{\lambda}}(S_{\lambda})) ,$$ where $D_{\frac{r}{s}}$ is the central simple division algebra corresponding to $\frac{r}{s}$ in the Brauer group $\mathbb{Q}/\mathbb{Z}=Br(\mathbb{Q}_p)$, and $\mathrm{M}_n(\mathrm{End}_{\mathcal{Q}_{\lambda}}(S_{\lambda}))$ denotes the algebra of $n\times n$ matrices over $\mathrm{End}_{\mathcal{Q}_{\lambda}}(S_{\lambda})$.

Moreover, for  $(r,s)=(1,0)$, if $det \begin{pmatrix}
        1 &0 \\ p & q
    \end{pmatrix}=q>0,$ we have  monomorphisms of $\mathbb{Q}_p$ algebras $$\mathcal{C}\hookrightarrow \mathrm{M}_q(\mathrm{End}_{\mathcal{Q}_{\lambda}}(S_{\lambda})) ,$$ for any untilt $\mathcal{C}$ of $C^{\flat}$ in characteristic 0.
\end{corollary}
\begin{proof}
    Consider the simple objects $\mathcal{O}_{\frac{r}{s}}$ in $\mathcal{P}(\phi,\phi+1]$; its image is isomorphic to $S_{\lambda}^{\oplus n}$ in $\mathcal{Q}_{\lambda}$ by the proof of Proposition \ref{prop: quotient categories}.  The quotient functor $q:\mathcal{P}(\phi,\phi+1]\rightarrow \mathcal{Q}_{\lambda}$ induces a nonzero algebras homomoprhism $$D_{\frac{r}{s}}\hookrightarrow \mathrm{M}_n(\mathrm{End}_{\mathcal{Q}_{\lambda}}(S_{\lambda})).$$ This is a monomorphism because $D_{\frac{r}{s}}$ is simple. A similar argument can be applied to all untilts of $C^{\flat}$ in characteristic 0.
\end{proof}
\begin{remark}
    In the special case when $\phi=\frac{1}{2}$, the division algebra is `le corps' studied by Colmez (see \cite[\S 5, \S 9] {BanachColmezspaces}, \cite[\S 7.3]{LeBrasresult}). The authors do not know whether these large $\mathbb{Q}_p$-division algebras are central simple in general.
\end{remark}

 Our method of studying the continuum envelopes of $\mathcal{A}$ is applicable to many examples of stability conditions $\sigma=(\mathcal{A},Z)$. In fact, we define a special type, called an arithmetic stability condition, to which our method applies.
    
\begin{definition}
    A stability condition $\sigma=(\mathcal{A},Z)$ is called an \textit{arithmetic stability condition} if the following conditions are satisfied. \begin{enumerate}
        
        \item The image of the central charge lies in $\mathbb{Z}+\sqrt{-1}\cdot\mathbb{Z}\subset\mathbb{C}$.
        \item For any primitive vector $(p,q)\in\mathbb{Z}^2$, there exists a stable object $E$ with $Z(E)=-p+\sqrt{-1}\cdot q$.
        \item For any two primitive vectors $(p_1,q_1), (p_2,q_2)\in\mathbb{Z}^2$ which are connected by a Farey geodesic with $\frac{p_1}{q_1}<\frac{p_2}{q_2}$ and any stable object $E$ with $Z(E)=-p_1+\sqrt{-1}\cdot q_1$, there exists a stable object $F$ with $Z(F)=-p_2+\sqrt{-1}\cdot q_2$ and $\Hom(E,F)\neq 0$.
    \end{enumerate}
\end{definition}

Many interesting stability conditions satisfy this definition. In fact, any stability condition whose central charge lies in $\mathbb{Q}+\sqrt{-1}\cdot\mathbb{Q}\subset\mathbb{C}$ is $\widetilde{GL}_2^+(\mathbb{R})$-equivalent to a stability condition that satisfies condition (1). Moreover, the set of stability conditions whose central charges lie in $\mathbb{Q}+\sqrt{-1}\cdot\mathbb{Q}\subset\mathbb{C}$ is dense in the rational component of the space of stability conditions (see \cite[Corollary 5.0.5]{APsheaves}).
\section{Final remarks}\label{section:final remarks}

\subsection{Stability conditions on twistor projective line}

There is an Archimedean analog of the Fargues--Fontaine curve, which is called the twistor projective line (Historically, one may instead view the Fargues–-Fontaine curve as a non-Archimedean analogue of the twistor projective line.). We briefly recall some basic facts about the twistor projective line here (see \cite{simpson1997mixed}, \cite{noteontwistorP1}).

The twistor projective line is another real structure on $\mathbb{P}_{\mathbb{C}}^1$, motivated by Penrose's idea of relating particles spinors to the light cone in Minkowski's space. The descent datum is given by the conjugate-linear involution: $z\mapsto -\frac{1}{\bar{z}}$. It is denoted by $\mathbb{P}_{tw}^1$. We have a natural map $\pi:\mathbb{P}_{\mathbb{C}}^1\rightarrow \mathbb{P}_{tw}^1$.

Consider the vector bundles on $\mathbb{P}_{tw}^1$. For $\lambda\in\frac{1}{2}\mathbb{Z}$, define $$\mathcal{O}_{\mathbb{P}_{tw}^1}(\lambda)=\begin{cases}
    \pi_*\mathcal{O}_{\mathbb{P}_{\mathbb{C}}^1}(2\lambda), & \text{if $\lambda\notin \mathbb{Z}$};  \\ \mathcal{L} \text{ such that $\pi^*\mathcal{L}=\mathcal{O}_{\mathbb{P}_{\mathbb{C}}^1}(2\lambda)$}, & \text{if $\lambda\in\mathbb{Z}$}.
\end{cases}$$ The slope of $\mathcal{O}_{\mathbb{P}_{tw}^1}(\lambda)$ is $\lambda$. 

The Harder--Narasimhan filtration of a vector bundle $\mathcal{E}$ on $\mathbb{P}_{tw}^1$ is given by $$\mathcal{E}=\mathcal{O}_{\mathbb{P}_{tw}^1}(\lambda_1)\oplus \mathcal{O}_{\mathbb{P}_{tw}^1}(\lambda_2)\oplus\cdots\oplus \mathcal{O}_{\mathbb{P}_{tw}^1}(\lambda_n),$$ where $\lambda_1\geq \lambda_2\geq\cdots\geq\lambda_n$ are half-integers in $\frac{1}{2}\mathbb{Z}$. There is also a well-known correspondence between the $U(1)$-equivariant vector bundles on $\mathbb{P}_{tw}^1$ and real Hodge structures (see \cite{simpson1997mixed}, \cite[\S 5.2]{noteontwistorP1}). 

Note that there exists a Bridgeland stability condition $\sigma=(\Coh(\mathbb{P}_{tw}^1), Z)$ on $D^b(\mathbb{P}_{tw}^1)$, where the central charge $Z$ is determined by the following formula $$Z(\mathcal{O}_{\mathbb{P}_{tw}^1}(\lambda))=\begin{cases}
-2\lambda+2\sqrt{-1}, &\text{if $\lambda\notin \mathbb{Z}$};  \\ -\lambda+\sqrt{-1}, & \text{if $\lambda\in\mathbb{Z}$}.
\end{cases}$$
  \begin{figure}[ht]
		\centering
		\includegraphics[scale=0.8]{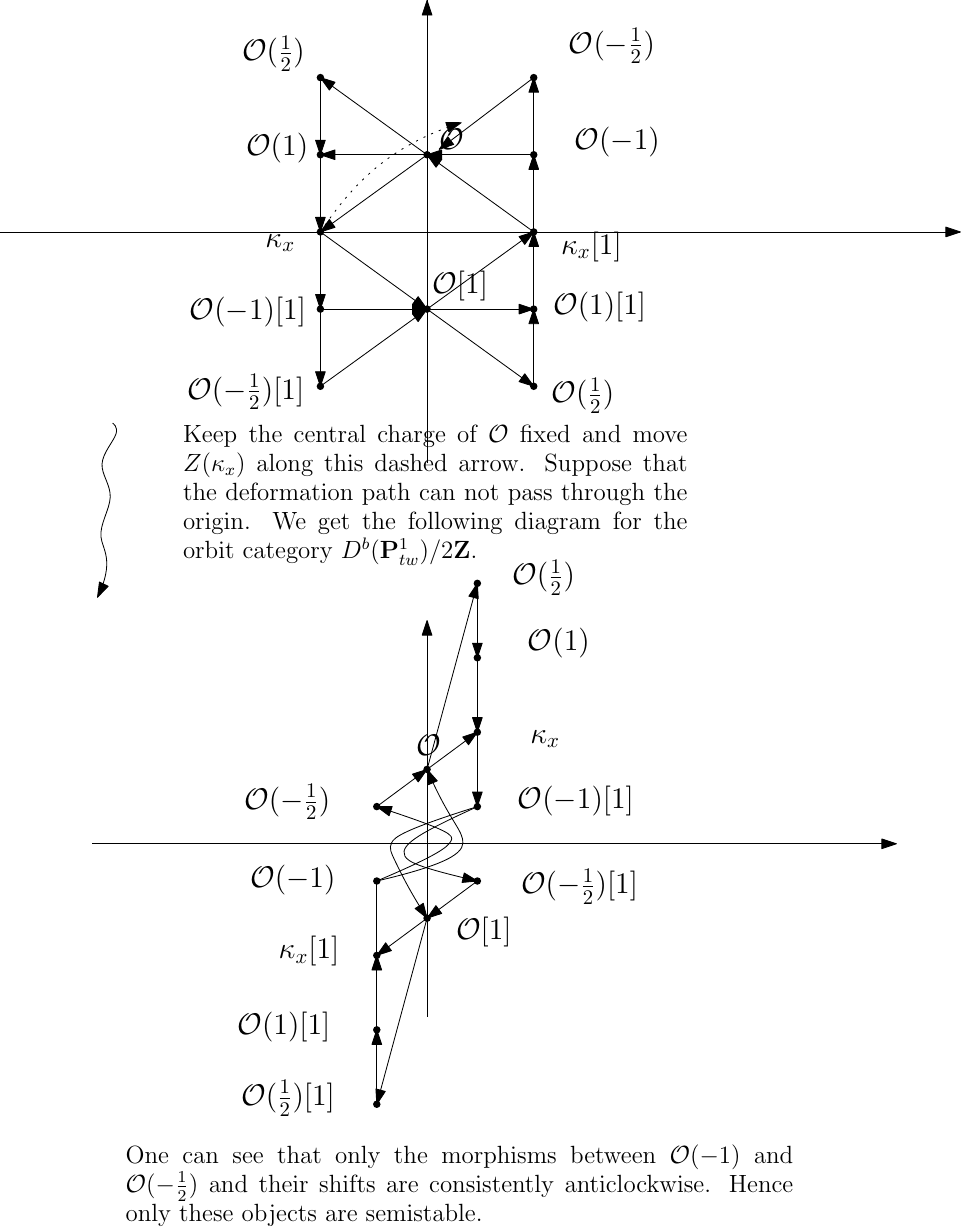}
		\caption{}
        \label{Figure 10}
	\end{figure} 
In fact, there are many other stability conditions on $D^b(\mathbb{P}_{tw}^1)$.
\begin{prop}
    For any stability condition in $\Stab(\mathbb{P}_{tw}^1)$, there exist $p>0$ and $\lambda\in\frac{1}{2}\mathbb{Z}$ such that $\mathcal{O}(\lambda)[p]$ and $\mathcal{O}(\lambda+\frac{1}{2})$ are semistable and $\phi(\mathcal{O}(\lambda)[p]), \phi(\mathcal{O}(\lambda+\frac{1}{2}))\in(r,r+1]$ for some $r\in\mathbb{R}$. 

    If $\phi(\mathcal{O}(\lambda)[1])< \phi(\mathcal{O}(\lambda+\frac{1}{2}))$, the multiples of the shifts of $\mathcal{O}(\lambda)[1], \mathcal{O}(\lambda+\frac{1}{2})$ are the only semistable objects. If $\phi(\mathcal{O}(\lambda)[1])\geq \phi(\mathcal{O}(\lambda+\frac{1}{2}))$, all indecomposable objects are semistable.
\end{prop}
\begin{proof}
    This proof is similar to that of \cite[Theorem 1.2]{StabilitymanifoldofP1}. We include a pictorial illustration (Figure \ref{Figure 10}) from the point of view of \cite{stabilityconditionsoncycliccategories}.

\end{proof}

Note that this is different from the case of the p-adic Fargues--Fontaine curve $X_{FF}$ (see Proposition \ref{Prop: stabilit conditions on FF curves}).

There are Bridgeland stability conditions on the product of twistor projective lines.
\begin{prop}
    There exist stability conditions on the product of twistor projective lines $\mathbb{P}_{tw}^1\times_{Spec(\mathbb{R})} \cdots \times_{Spec(\mathbb{R})} \mathbb{P}_{tw}^1$.
\end{prop}
\begin{proof}
    This follows from the construction in \cite{Stabilityconditionsonproductvarieties}. Note that although we assumed that the scheme is defined over an algebraically closed field in \cite{Stabilityconditionsonproductvarieties}, this assumption is not necessary. In fact, the only place where we used this assumption was in \cite[Theorem 3.3]{Stabilityconditionsonproductvarieties}. We used it to apply Bertini's theorem, which holds over an infinite field. 
    
    In fact, when the stability condition $\sigma$ is a numerical stability condition (see \cite[\S 2]{bayer2014projectivity} for the definition of numerical stability conditions), we can prove \cite[Theorem 3.3]{Stabilityconditionsonproductvarieties} by Grothendieck--Riemann--Roch for any projective scheme over any field $k$. So, the construction in \cite{Stabilityconditionsonproductvarieties} works for any field $k$ in that case.
\end{proof}
\begin{remark}
    Unfortunately, as the schematic Fargues--Fontaine curve is not of finite type, we can not apply the Abramovich--Polishchuk's construction to it, and hence the construction in \cite{Stabilityconditionsonproductvarieties} can not be applied directly to the product of Fargues--Fontaine curves (as diamonds). However, we still expect that many results in \cite{Stabilityconditionsinfamilies} are applicable to family of Fargues--Fontaine curves over a perfectoid space.
\end{remark}

\subsection{Similarity to elliptic curves and noncommutative tori}\label{noncommutative tori}

As we have mentioned before, there is strong similarity between elliptic curves and the Fargues--Fontaine curve. In particular, Corollary \ref{proposition of noetherian and aritinian} holds for elliptic curves. We will recall the notion of noncommutative tori and its relation to $\mathcal{P}(\phi,\phi+1]$, where $\mathcal{P}$ is the standard slicing in $D^b(E)$, the bounded derived category of coherent sheaves on an elliptic curve $E$ (see \cite{Quasicoherentsheaves}).

Noncommutative geometry was introduced by Alain Connes in the 1980s (see \cite{Noncommutativedifferentialgeometry}). Noncommutative tori are fundamental objects in this theory. They arise from Kronecker foliations. We briefly recall their construction.

Let $T$ be the two-dimensional torus $T=\mathbb{R}^2/\mathbb{Z}^2,$ with the Kronecker foliation associated to an irrational number $\theta'$, i.e., defines by the differential equation $$dy=\theta'dx.$$ 

We use $X$ to denote the quotient space of $T$ by this foliation. It is easy to see that the classical measure theory can not distinguish between $X$ and a one-point space. However, in the theory of noncommutative geometry, we can associate a smooth groupoid to a foliation and then associate a $C^*$-algebra to that groupoid (see \cite[Chapter 2.9]{connes1994noncommutative}. Following these standard procedures, one can show that the $C^*$-algebra $A_{\theta}$ associated to the Kronecker foliation is generated by two unitaries $U$ and $V$ such that $$VU=exp(2\pi i\theta)UV,$$ where $\theta$ is another irrational number determined by $\theta'$. Moreover, one can equip $A_{\theta}$ with a complex structure given by some $\tau\in\mathbb{C}\backslash \mathbb{R}$. The same $\tau$ also determines a complex elliptic curve $X_{\tau}$. Polishchuk studied the category of holomorphic vector bundles with respect to such a holomorphic structure in \cite{Classficationofholomorphicvectorbundles}. It is related to the Bridgeland stability condition  on $D^b(X_{\tau})$ in the following way.

\begin{theorem}[{\cite[Corollary 2.1]{Classficationofholomorphicvectorbundles}}]
    Let $X_{\tau}$ be the associated complex elliptic curve, and let $\sigma=(\mathcal{P},Z)$ be the standard Mumford stability condition on $D^b(X_{\tau})$. Then the category of holomorphic vector bundles on $A_{\theta,\tau}$ is equivalent to $\mathcal{P}(\phi,\phi+1]$, where $\phi$ and $\theta$ are related in the usual way.
\end{theorem}

Motivated by this result, we ask the following question.

\textbf{Question}: Does the Fargues--Fontaine curve $X_{FF}$ have noncommutative deformations $X_{FF}^{\theta}$ such that the category of vector bundles on $X_{FF}^{\theta}$ is equivalent to $\mathcal{P}(\phi,\phi+1]$?
\bibliographystyle{alpha}
	\bibliography{bibfile}

@Book{farguesfontaine-courbes,
 Author = {Fargues, Laurent and Fontaine, Jean-Marc},
 Title = {Courbes et fibr{\'e}s vectoriels en th{\'e}orie de {Hodge} {{\(p\)}}-adique},
 FSeries = {Ast{\'e}risque},
 Series = {Ast{\'e}risque},
 ISSN = {0303-1179},
 Volume = {406},
 ISBN = {978-2-85629-896-1},
 Year = {2018},
 Publisher = {Paris: Soci{\'e}t{\'e} Math{\'e}matique de France (SMF)},
 Language = {French},
 zbMATH = {7005651},
 Zbl = {1470.14001}
}

@Article{scholze-perfectoidspacesIHES,
 Author = {Scholze, Peter},
 Title = {Perfectoid spaces},
 FJournal = {Publications Math{\'e}matiques},
 Journal = {Publ. Math., Inst. Hautes {\'E}tud. Sci.},
 ISSN = {0073-8301},
 Volume = {116},
 Pages = {245--313},
 Year = {2012},
 Language = {English},
 DOI = {10.1007/s10240-012-0042-x},
 zbMATH = {6120994},
 Zbl = {1263.14022}
}

@article{bridgeland2007stability,
 AUTHOR = {Bridgeland, Tom},
 TITLE = {Stability conditions on triangulated categories},
 JOURNAL = {Ann. of Math. (2)},
 FJOURNAL = {Annals of Mathematics. Second Series},
 VOLUME = {166},
 YEAR = {2007},
 NUMBER = {2},
 PAGES = {317--345}
}

@article{polishchuk2007constant,
    AUTHOR = {Polishchuk, A.},
    TITLE = {Constant families of {$t$}-structures on derived categories of
    coherent sheaves},
    JOURNAL = {Mosc. Math. J.},
    FJOURNAL = {Moscow Mathematical Journal},
    VOLUME = {7},
    YEAR = {2007},
    NUMBER = {1},
    PAGES = {109--134, 167},
}

@article {APsheaves,
AUTHOR = {Abramovich, Dan and Polishchuk, Alexander},
TITLE = {Sheaves of {$t$}-structures and valuative criteria for stable
complexes},
JOURNAL = {J. Reine Angew. Math.},
FJOURNAL = {Journal f\"{u}r die Reine und Angewandte Mathematik. [Crelle's
Journal]},
VOLUME = {590},
YEAR = {2006},
PAGES = {89--130}
}

@book{hotta2007d,
	title={D-modules, perverse sheaves, and representation theory},
	author={Hotta, Ryoshi and Tanisaki, Toshiyuki},
	volume={236},
	year={2007},
	publisher={Springer Science \& Business Media}
}

@article{bayer2014projectivity,
 AUTHOR = {Bayer, Arend and Macr\`\i , Emanuele},
 TITLE = {Projectivity and birational geometry of {B}ridgeland moduli
 spaces},
 JOURNAL = {J. Amer. Math. Soc.},
 FJOURNAL = {Journal of the American Mathematical Society},
 VOLUME = {27},
 YEAR = {2014},
 NUMBER = {3},
 PAGES = {707--752}
}

@article{beilinson1982faisceaux,
	title={Faisceaux pervers, Analyse et topologie sur les espaces singuliers {(I) CIRM}, 6--10 juillet 1981},
	author={Beilinson, A.A. and Bernstein, J. and Deligne, P.},
	journal={Ast{\'e}risque},
	volume={100},
	year={1982}
}

@article{bayer2014mmp,
AUTHOR = {Bayer, Arend and Macr\`\i , Emanuele},
TITLE = {M{MP} for moduli of sheaves on {K}3s via wall-crossing: nef
and movable cones, {L}agrangian fibrations},
JOURNAL = {Invent. Math.},
FJOURNAL = {Inventiones Mathematicae},
VOLUME = {198},
YEAR = {2014},
NUMBER = {3},
PAGES = {505--590}
}

@article{bridgeland2008stability,
AUTHOR = {Bridgeland, Tom},
TITLE = {Stability conditions on {$K3$} surfaces},
JOURNAL = {Duke Math. J.},
FJOURNAL = {Duke Mathematical Journal},
VOLUME = {141},
YEAR = {2008},
NUMBER = {2},
PAGES = {241--291}
}

@article{baye2011bridgeland,
 AUTHOR = {Bayer, Arend and Macr\`\i , Emanuele and Toda, Yukinobu},
 TITLE = {Bridgeland stability conditions on threefolds {I}:
 {B}ogomolov-{G}ieseker type inequalities},
 JOURNAL = {J. Algebraic Geom.},
 FJOURNAL = {Journal of Algebraic Geometry},
 VOLUME = {23},
 YEAR = {2014},
 NUMBER = {1},
 PAGES = {117--163}
}

@article{kontsevich2008stability,
	title={Stability structures, motivic {D}onaldson-{T}homas invariants and cluster transformations},
	author={Kontsevich, Maxim and Soibelman, Yan},
	journal={arXiv preprint arXiv:0811.2435},
	year={2008}
}

@inproceedings{douglas2002dirichlet,
		 AUTHOR = {Douglas, Michael R.},
		 TITLE = {Dirichlet branes, homological mirror symmetry, and stability},
		 BOOKTITLE = {Proceedings of the {I}nternational {C}ongress of
		 {M}athematicians, {V}ol. {III} ({B}eijing, 2002)},
		 PAGES = {395--408},
		 PUBLISHER = {Higher Ed. Press, Beijing},
		 YEAR = {2002}
	}

@article {Stabilityconditionsonproductvarieties,
	AUTHOR = {Liu, Yucheng},
	TITLE = {Stability conditions on product varieties},
	JOURNAL = {J. Reine Angew. Math.},
	FJOURNAL = {Journal f\"{u}r die Reine und Angewandte Mathematik. [Crelle's
	Journal]},
	VOLUME = {770},
	YEAR = {2021},
	PAGES = {135--157},
	ISSN = {0075-4102},
	MRCLASS = {14F08},
	MRNUMBER = {4193465},
	DOI = {10.1515/crelle-2020-0010},
	URL = {https://doi.org/10.1515/crelle-2020-0010},
}

@Article{CurvecountingtheoriesviastableojectsI,
	Author = {Toda, Yukinobu},
	Title = {Curve counting theories via stable objects. {I}: {DT}/{PT} correspondence},
	FJournal = {Journal of the American Mathematical Society},
	Journal = {J. Am. Math. Soc.},
	ISSN = {0894-0347},
	Volume = {23},
	Number = {4},
	Pages = {1119--1157},
	Year = {2010},
	Language = {English},
	DOI = {10.1090/S0894-0347-10-00670-3},
	zbMATH = {5797891},
	Zbl = {1207.14020}
}

@Article{Stabilityconditionsinfamilies,
	Author = {Bayer, Arend and Lahoz, Mart{\'{\i}} and Macr{\`{\i}}, Emanuele and Nuer, Howard and Perry, Alexander and Stellari, Paolo},
	Title = {Stability conditions in families},
	FJournal = {Publications Math{\'e}matiques},
	Journal = {Publ. Math., Inst. Hautes {\'E}tud. Sci.},
	ISSN = {0073-8301},
	Volume = {133},
	Pages = {157--325},
	Year = {2021},
	Language = {English},
	DOI = {10.1007/s10240-021-00124-6},
	zbMATH = {7395794},
	Zbl = {1481.14033}
}

@Article{Scatteringdiagrams,
	Author = {Bridgeland, Tom},
	Title = {Scattering diagrams, {Hall} algebras and stability conditions},
	FJournal = {Algebraic Geometry},
	Journal = {Algebr. Geom.},
	ISSN = {2313-1691},
	Volume = {4},
	Number = {5},
	Pages = {523--561},
	Year = {2017},
	Language = {English},
	zbMATH = {6849618},
	Zbl = {1388.16013}
}

@Article{CurvecountingtheoriesviastableobjectsII,
	Author = {Toda, Yukinobu},
	Title = {Curve counting theories via stable objects. {II}: {DT}/{ncDT} flop formula},
	FJournal = {Journal f{\"u}r die Reine und Angewandte Mathematik},
	Journal = {J. Reine Angew. Math.},
	ISSN = {0075-4102},
	Volume = {675},
	Pages = {1--51},
	Year = {2013},
	Language = {English},
	DOI = {10.1515/CRELLE.2011.176},
	zbMATH = {6146431},
	Zbl = {1267.14049}
}

@article {Noncommutativedifferentialgeometry,
	AUTHOR = {Connes, Alain},
	TITLE = {Noncommutative differential geometry},
	JOURNAL = {Inst. Hautes \'{E}tudes Sci. Publ. Math.},
	FJOURNAL = {Institut des Hautes \'{E}tudes Scientifiques. Publications
	Math\'{e}matiques},
	NUMBER = {62},
	YEAR = {1985},
	PAGES = {257--360},
	ISSN = {0073-8301},
	MRCLASS = {58G12 (16A61 18F25 19D55 19K56 46L80 53C99)},
	MRNUMBER = {823176},
	MRREVIEWER = {Jonathan M. Rosenberg},
	URL = {http://www.numdam.org/item?id=PMIHES_1985__62__257_0},
}

@Misc{Langbookondiophantineapproximations,
	Author = {Lang, Serge},
	Title = {Introduction to {D}iophantine approximations},
	Year = {1966},
	Language = {English},
	HowPublished = {Addison-{Wesley} {Series} in {Mathematics}. {Reading}, {Mass}. etc.: {Addison}-{Wesley} {Publishing} {Company}. {VIII}, 83 p. (1966).},
	zbMATH = {3232020},
	Zbl = {0144.04005}
}

@Article{BanachColmezspaces,
	Author = {Colmez, Pierre},
	Title = {Espaces de {Banach} de dimension finie. ({Finite}-dimensional {Banach} {Spaces})},
	FJournal = {Journal of the Institute of Mathematics of Jussieu},
	Journal = {J. Inst. Math. Jussieu},
	ISSN = {1474-7480},
	Volume = {1},
	Number = {3},
	Pages = {331--439},
	Year = {2002},
	Language = {French},
	DOI = {10.1017/S1474748002000099},
	zbMATH = {2011902},
	Zbl = {1044.11102}
}

@Article{LeBrasresult,
	Author = {Le Bras, Arthur-C{\'e}sar},
	Title = {Banach-{Colmez} spaces and coherent sheaves on the {Fargues}-{Fontaine} curve},
	FJournal = {Duke Mathematical Journal},
	Journal = {Duke Math. J.},
	ISSN = {0012-7094},
	Volume = {167},
	Number = {18},
	Pages = {3455--3532},
	Year = {2018},
	Language = {English},
	DOI = {10.1215/00127094-2018-0034},
	zbMATH = {7009770},
	Zbl = {1439.14085}
}

@article {Atiyahclassification,
	AUTHOR = {Atiyah, M. F.},
	TITLE = {Vector bundles over an elliptic curve},
	JOURNAL = {Proc. London Math. Soc. (3)},
	FJOURNAL = {Proceedings of the London Mathematical Society. Third Series},
	VOLUME = {7},
	YEAR = {1957},
	PAGES = {414--452},
	ISSN = {0024-6115},
	MRCLASS = {14.55 (14.20)},
	MRNUMBER = {131423},
	MRREVIEWER = {F. Hirzebruch},
	DOI = {10.1112/plms/s3-7.1.414},
	URL = {https://doi.org/10.1112/plms/s3-7.1.414},
}

@Article{Quadraticdifferentialsasstabilityconditions,
	Author = {Bridgeland, Tom and Smith, Ivan},
	Title = {Quadratic differentials as stability conditions},
	FJournal = {Publications Math{\'e}matiques},
	Journal = {Publ. Math., Inst. Hautes {\'E}tud. Sci.},
	ISSN = {0073-8301},
	Volume = {121},
	Pages = {155--278},
	Year = {2015},
	Language = {English},
	DOI = {10.1007/s10240-014-0066-5},
	URL = {eprints.whiterose.ac.uk/84399/1/BrSm_arXiv-Version2-Corrected.pdf},
	zbMATH = {6460465},
	Zbl = {1328.14025}
}

@Article{Flatsurfacesandstabilityconditions,
	Author = {Haiden, F. and Katzarkov, L. and Kontsevich, M.},
	Title = {Flat surfaces and stability structures},
	FJournal = {Publications Math{\'e}matiques},
	Journal = {Publ. Math., Inst. Hautes {\'E}tud. Sci.},
	ISSN = {0073-8301},
	Volume = {126},
	Pages = {247--318},
	Year = {2017},
	Language = {English},
	DOI = {10.1007/s10240-017-0095-y},
	zbMATH = {6827888},
	Zbl = {1390.32010}
}

@article{Filtrationsandtorsionpairs,
 author = {Liu, Yucheng},
 title = {Filtrations and torsion pairs in {Abramovich}-{Polishchuk}'s heart},
 fjournal = {Peking Mathematical Journal},
 journal = {Peking Math. J.},
 issn = {2096-6075},
 volume = {8},
 number = {4},
 pages = {767--789},
 year = {2025},
 language = {English},
 doi = {10.1007/s42543-024-00084-w},
 zbMATH = {8180066}
}

@article{stabilityconditionsoncycliccategories,
	title={Stability conditions on cyclic categories {I}: basic definitions and examples},
	author={Yucheng Liu},
	journal={arXiv preprint arXiv:2211.16939},
	year={2022}
}

@Article{StabilitymanifoldofP1,
	Author = {Okada, So},
	Title = {Stability manifold of {{\(\mathbb P^1\)}}},
	FJournal = {Journal of Algebraic Geometry},
	Journal = {J. Algebr. Geom.},
	ISSN = {1056-3911},
	Volume = {15},
	Number = {3},
	Pages = {487--505},
	Year = {2006},
	Language = {English},
	DOI = {10.1090/S1056-3911-06-00432-2},
	zbMATH = {5135138},
	Zbl = {1117.14021}
}

@Article{Quasicoherentsheaves,
 Author = {Polishchuk, A.},
 Title = {Quasicoherent sheaves on complex noncommutative two-tori},
 FJournal = {Selecta Mathematica. New Series},
 Journal = {Sel. Math., New Ser.},
 ISSN = {1022-1824},
 Volume = {13},
 Number = {1},
 Pages = {137--173},
 Year = {2007},
 Language = {English},
 DOI = {10.1007/s00029-007-0034-8},
 zbMATH = {5237949},
 Zbl = {1135.32011}
}

@Article{Classficationofholomorphicvectorbundles,
 Author = {Polishchuk, A.},
 Title = {Classification of holomorphic vector bundles on noncommutative two-tori},
 FJournal = {Documenta Mathematica},
 Journal = {Doc. Math.},
 ISSN = {1431-0635},
 Volume = {9},
 Pages = {163--181},
 Year = {2004},
 Language = {English},
 zbMATH = {2095723},
 Zbl = {1048.32012}
}

@Article{themodularsurfaceandcontinuedfractions,
 Author = {Series, Caroline},
 Title = {The modular surface and continued fractions},
 FJournal = {Journal of the London Mathematical Society. Second Series},
 Journal = {J. Lond. Math. Soc., II. Ser.},
 ISSN = {0024-6107},
 Volume = {31},
 Pages = {69--80},
 Year = {1985},
 Language = {English},
 DOI = {10.1112/jlms/s2-31.1.69},
 zbMATH = {3867655},
 Zbl = {0545.30001}
}

@Misc{Themodulargroup,
 Author = {Jones, G. A. and Singerman, D. and Wicks, K.},
 Title = {The modular group and generalized {Farey} graphs},
 Year = {1991},
 Language = {English},
 HowPublished = {Groups, {Vol}. 2, {Proc}. {Int}. {Conf}., {St}. {Andrews}/{UK} 1989, {Lond}. {Math}. {Soc}. {Lect}. {Note} {Ser}. 160, 316-338 (1991).},
 zbMATH = {4202583},
 Zbl = {0728.20040}
}

@Article{curvatureandrank,
 Author = {Brock, Jeffrey and Farb, Benson},
 Title = {Curvature and rank of {Teichm{\"u}ller} space},
 FJournal = {American Journal of Mathematics},
 Journal = {Am. J. Math.},
 ISSN = {0002-9327},
 Volume = {128},
 Number = {1},
 Pages = {1--22},
 Year = {2006},
 Language = {English},
 DOI = {10.1353/ajm.2006.0003},
 zbMATH = {5015886},
 Zbl = {1092.32008}
}

@article{Circlehomeomorphisms,
  title={Circle homeomorphisms with square summable diamond shears},
  author={{\v{S}}ari{\'c}, Dragomir and Wang, Yilin and Wolfram, Catherine},
  journal={arXiv preprint arXiv:2211.11497},
  year={2022}
}

@article{Differentheartsonellipticcurves,
 author = {Liu, Yucheng},
 title = {Different hearts on elliptic curves},
 fjournal = {Acta Mathematica Sinica. English Series},
 journal = {Acta Math. Sin., Engl. Ser.},
 issn = {1439-8516},
 volume = {41},
 number = {3},
 pages = {847--853},
 year = {2025},
 language = {English},
 doi = {10.1007/s10114-025-3286-3},
 zbMATH = {8019809},
 Zbl = {1573.14061}
}

@book {HTT,
    AUTHOR = {Lurie, Jacob},
     TITLE = {Higher topos theory},
    SERIES = {Annals of Mathematics Studies},
    VOLUME = {170},
 PUBLISHER = {Princeton University Press, Princeton, NJ},
      YEAR = {2009},
     PAGES = {xviii+925},
      ISBN = {978-0-691-14049-0; 0-691-14049-9},
   MRCLASS = {18-02 (18B25 18E35 18G30 18G55 55U40)},
  MRNUMBER = {2522659},
MRREVIEWER = {Mark\ Hovey},
       DOI = {10.1515/9781400830558},
       URL = {https://doi.org/10.1515/9781400830558},
}

@Book{Pierceassociativealgebra,
 Author = {Pierce, Richard S.},
 Title = {Associative algebras},
 FSeries = {Graduate Texts in Mathematics},
 Series = {Grad. Texts Math.},
 ISSN = {0072-5285},
 Volume = {88},
 Year = {1982},
 Publisher = {Springer, Cham},
 Language = {English},
 zbMATH = {3783206},
 Zbl = {0497.16001}
}

@Article{Extensionsofvectorbundlesonffcurve,
 Author = {Birkbeck, Christopher and Feng, Tony and Hansen, David and Hong, Serin and Li, Qirui and Wang, Anthony and Ye, Lynnelle},
 Title = {Extensions of vector bundles on the {Fargues}-{Fontaine} curve},
 FJournal = {Journal of the Institute of Mathematics of Jussieu},
 Journal = {J. Inst. Math. Jussieu},
 ISSN = {1474-7480},
 Volume = {21},
 Number = {2},
 Pages = {487--532},
 Year = {2022},
 Language = {English},
 DOI = {10.1017/S1474748020000183},
 zbMATH = {7482201},
 Zbl = {1497.14024}
}

@article{fargues2021geometrization,
  title={Geometrization of the local {L}anglands correspondence},
  author={Fargues, Laurent and Scholze, Peter},
  journal={arXiv preprint arXiv:2102.13459},
  year={2021}
}

@book{connes1994noncommutative,
  title={Noncommutative geometry},
  author={Connes, Alain},
  year={1994},
  publisher={Springer}
}

@article{noteontwistorP1,
  title={Notes on the Twistor $\mathbb{P}^1 $},
  author={Woit, Peter},
  journal={arXiv preprint arXiv:2202.02657},
  year={2022}
}

@article{simpson1997mixed,
  title={Mixed twistor structures},
  author={Simpson, Carlos},
  journal={arXiv preprint alg-geom/9705006},
  year={1997}
}

@article {FMtransformsandvectorbundlesonellipticcurve,
	AUTHOR = {Hein, Georg and Ploog, David},
	TITLE = {Fourier-{M}ukai transforms and stable bundles on elliptic
	curves},
	JOURNAL = {Beitr\"{a}ge Algebra Geom.},
	FJOURNAL = {Beitr\"{a}ge zur Algebra und Geometrie. Contributions to Algebra
	and Geometry},
	VOLUME = {46},
	YEAR = {2005},
	NUMBER = {2},
	PAGES = {423--434},
	ISSN = {0138-4821},
	MRCLASS = {14H60 (14H52)},
	MRNUMBER = {2196927},
	MRREVIEWER = {Ana Cristina L\'{o}pez-Mart\'{\i}n},
}

@article{scholze2017etale,
  title={{\'E}tale cohomology of diamonds},
  author={Scholze, Peter},
  journal={arXiv preprint arXiv:1709.07343},
  year={2017}
}

@book{scholze2020berkeley,
  title={Berkeley Lectures on p-adic Geometry (AMS-207)},
  author={Scholze, Peter and Weinstein, Jared},
  year={2020},
  publisher={Princeton University Press}
}

@article{anschutz2021fourier,
  title={A {F}ourier Transform for {B}anach-{C}olmez spaces},
  author={Ansch{\"u}tz, Johannes and Bras, Arthur-C{\'e}sar Le},
  journal={arXiv preprint arXiv:2111.11116},
  year={2021}
}

@book{kedlaya2015relative,
 author = {Kedlaya, Kiran S. and Liu, Ruochuan},
 title = {Relative {{\(p\)}}-adic {Hodge} theory: foundations},
 fseries = {Ast{\'e}risque},
 series = {Ast{\'e}risque},
 issn = {0303-1179},
 volume = {371},
 isbn = {978-2-85629-807-7},
 year = {2015},
 publisher = {Paris: Soci{\'e}t{\'e} Math{\'e}matique de France (SMF)},
 language = {English},
 zbMATH = {6456685},
 Zbl = {1370.14025}
}

@book{aigner2015markov,
  title={Markov's theorem and 100 years of the uniqueness conjecture},
  author={Aigner, Martin},
  year={2015},
  publisher={Springer}
}

@InCollection{Infinitedimensionalvectorbundles,
 Author = {Drinfeld, Vladimir},
 Title = {Infinite-dimensional vector bundles in algebraic geometry: an introduction},
 BookTitle = {The unity of mathematics. In honor of the ninetieth birthday of I. M. Gelfand. Papers from the conference held in Cambridge, MA, USA, August 31--September 4, 2003.},
 ISBN = {0-8176-4076-2},
 Pages = {263--304},
 Year = {2006},
 Publisher = {Boston, MA: Birkh{\"a}user},
 Language = {English},
 zbMATH = {5050079},
 Zbl = {1108.14012}
}

@Book{Ergodictheory,
 Author = {Cornfeld, I. P. and Fomin, S. V. and Sinai, Ya. G.},
 Title = {Ergodic theory. {Transl}. from the {Russian} by {A}. {B}. {Sossinskii}},
 FSeries = {Grundlehren der Mathematischen Wissenschaften},
 Series = {Grundlehren Math. Wiss.},
 ISSN = {0072-7830},
 Volume = {245},
 Year = {1982},
 Publisher = {Springer, Cham},
 Language = {English},
 zbMATH = {3775877},
 Zbl = {0493.28007}
}

@Article{Moduliofp-divisiblegroups,
 Author = {Scholze, Peter and Weinstein, Jared},
 Title = {Moduli of {{\(p\)}}-divisible groups},
 FJournal = {Cambridge Journal of Mathematics},
 Journal = {Camb. J. Math.},
 ISSN = {2168-0930},
 Volume = {1},
 Number = {2},
 Pages = {145--237},
 Year = {2013},
 Language = {English},
 DOI = {10.4310/CJM.2013.v1.n2.a1},
 zbMATH = {6292858},
 Zbl = {1349.14149}
}

@article{Modularvariants,
  title={Modular variants of p-adic fundamental sequence},
  author={Du, Heng and Jiang, Qingyuan and Liu, Yucheng},
  journal={arXiv preprint arXiv:2606.12220},
  year={2026}
}
\end{document}